\documentclass[reqno,12pt]{amsart}

\numberwithin{equation}{section}

\newcommand{\1}{\mathbf {1}}
\newcommand{\N}{{\mathbb N}}
\newcommand{\Z}{{\mathbb Z}}
\newcommand{\Q}{{\mathbb Q}}

\newcommand{\C}{{\mathbb C}}
\newcommand{\F}{{\mathcal F}}
\newcommand{\G}{{\mathcal G}}
\newcommand{\I}{{\mathcal I}}
\newcommand{\J}{{\mathcal J}}
\newcommand{\M}{{\mathcal M}}
\newcommand{\W}{{\mathcal W}}

\newcommand{\al}{\alpha}
\newcommand{\be}{\beta}

\newcommand{\om}{\omega}
\newcommand{\sg}{\sigma}
\newcommand{\vep}{\varepsilon}
\newcommand{\tom}{\tilde{\omega}}
\newcommand{\ty}{\tilde{Y}}
\newcommand{\hM}{\hat{M}}
\newcommand{\hW}{\hat{W}}
\newcommand{\la}{\langle}
\newcommand{\op}{\oplus}
\newcommand{\ot}{\otimes}
\newcommand{\ra}{\rangle}
\newcommand{\bv}{\mathbf{v}}

\DeclareMathOperator{\End}{End}
\DeclareMathOperator{\Hom}{Hom}
\DeclareMathOperator{\Ker}{Ker}

\DeclareMathOperator{\Vir}{Vir}
\DeclareMathOperator{\ch}{ch}
\DeclareMathOperator{\spn}{span}

\DeclareMathOperator{\wt}{wt}

\newtheorem{thm}{Theorem}[section]
\newtheorem{prop}[thm]{Proposition}
\newtheorem{lem}[thm]{Lemma}
\newtheorem{cor}[thm]{Corollary}
\newtheorem{rmk}[thm]{Remark}
\newtheorem{definition}[thm]{Definition}
\begin{document}

\title[fixed point subalgebra]
{The fixed point subalgebra of a lattice vertex
operator algebra by an automorphism of order three}

\author[K. Tanabe]{Kenichiro Tanabe $^\dagger$}
\address[K. Tanabe]{Institute of Mathematics,
University of Tsukuba,
Tsukuba 305-8571, Japan}
\email{tanabe@math.tsukuba.ac.jp}

\author[H. Yamada]{Hiromichi Yamada $^\ddagger$}
\address[H. Yamada]{Department of Mathematics,
Hitotsubashi University,
Kunitachi, Tokyo 186-8601, Japan}
\email{yamada@math.hit-u.ac.jp}

\thanks{$^\dagger$ Partially supported by JSPS Grant-in-Aid for
Scientific Research No. 17740002
\\
\hspace*{3mm} $^\ddagger$ Partially supported by JSPS Grant-in-Aid for
Scientific Research No. 17540016
}


\begin{abstract}
We study the subalgebra of the lattice vertex operator algebra
$V_{\sqrt{2}A_2}$ consisting of the fixed points of an automorphism
which is induced from an order $3$ isometry of the root lattice
$A_2$. We classify the simple modules for the subalgebra. The
rationality and the $C_2$-cofiniteness are also established.
\end{abstract}

\maketitle

\section{Introduction}
The space of fixed points of an automorphism group of finite order
in a vertex operator algebra is a vertex operator subalgebra. The
study of such subalgebras and their modules is called the orbifold
theory. It is a rich field both in the conformal field theory and in
the theory of vertex operator algebras. However, the orbifold theory
is difficult to study in general. One of the reason is that the
subalgebra of fixed points usually has more complicated structure
than the original vertex operator algebra.

The first example of the orbifold theory in vertex operator algebras
is the moonshine module $V^\natural$ by Frenkel, Lepowsky, and
Meurman \cite{FLM}. In their book $V^\natural$ was constructed as an
extension of $V_\Lambda^+$ by its simple module $V_\Lambda^{T,+}$,
where $V_\Lambda^+$ is the space of fixed points of an automorphism
$\theta$ of order two in the Leech lattice vertex operator algebra
$V_\Lambda$. This construction is called a $2B$-orbifold
construction because $\theta$ corresponds to a $2B$ involution of
the monster simple group. More generally, a vertex operator algebra
$V_L$ associated with an arbitrary positive definite even lattice
$L$ was defined in \cite{FLM}. Those lattice vertex operator
algebras provide a large family of vertex operator algebras. Such a
lattice vertex operator algebra admits an automorphism $\theta$ of
order two, which is a lift of the isometry $\alpha \mapsto -\alpha$
of the underlying lattice $L$. The orbifold theory for the fixed
point subalgebra $V_L^+$ of $\theta$ has been developed extensively.
In fact, the simple $V_L^+$-modules were classified \cite{AD} and
the fusion rules were determined \cite{ADL}. Furthermore, the
$C_2$-cofiniteness of $V_L^+$ was established \cite{ABD, Yamskulna}.

In this paper we study the fixed point subalgebra by an automorphism
of order three for a certain lattice vertex operator algebra.
Namely, let $L = \sqrt{2}A_2$ be $\sqrt{2}$ times an ordinary root
lattice of type $A_2$ and $\tau$ be an isometry of the root lattice
of type $A_2$ induced from an order three permutation on the set of
positive roots. We classify the simple modules for the subalgebra
$V_L^\tau$ of fixed points by $\tau$. Moreover, we show that
$V_L^\tau$ is rational and $C_2$-cofinite.

In previous papers \cite{DLTYY, KLY2} we have already discussed the
vertex operator algebra $V_L^\tau$. It was shown that $V_L^\tau =
M^0 \op W^0$ is a direct some of a subalgebra $M^0$ and its simple
highest weight module $W^0$. Actually, $M^0$ is a tensor product of
a $W_3$ algebra of central charge $6/5$ and a $W_3$ algebra of
central charge $4/5$. The property of a $W_3$ algebra of central
charge $6/5$ as the first component of the tensor product $M^0$ was
investigated in \cite{DLTYY}. It is generated by the Virasoro
element $\tom^1$ and a weight $3$ vector $J$. On the other hand, the
second component of $M^0$, which is a $W_3$ algebra of central
charge $4/5$, was studied in \cite{KMY}. It is generated by the
Virasoro element $\tom^2$ and a weight $3$ vector $K$. Each of these
$W_3$ algebras possesses a symmetry of order three. The order three
symmetry of the second $W_3$ algebra is related to the $\Z_3$ part
of $L^\perp/L \cong \Z_2 \times \Z_2 \times \Z_3$, where $L^\perp$
denotes the dual lattice of $L$. As an $M^0$-module, $W^0$ is
generated by a highest weight vector $P$ of weight $2$. Thus the
vertex operator algebra $V_L^\tau$ is generated by the five elements
$\tom^1$, $\tom^2$, $J$, $K$, and $P$.

There are $12$ inequivalent simple $V_L$-modules, which correspond
to the cosets of $L$ in its dual lattice $L^\perp$ (cf. \cite{D1}).
Let $(U,Y_U)$ be a simple $V_L$-module. One can define a new simple
$V_L$-module $(U \circ \tau, Y_{U \circ \tau})$ by $U \circ \tau =
U$ as vector spaces and $Y_{U \circ \tau}(v,z) = Y_U(\tau v,z)$ for
$v \in V_L$. Then $U \mapsto U \circ \tau$ is a permutation on the
set of simple $V_L$-modules. In the case where $U$ and $U \circ
\tau$ are equivalent $V_L$-modules, $U$ is said to be $\tau$-stable.
If $U$ is $\tau$-stable, then the eigenspace $U(\vep)$ of $\tau$
with eigenvalue $\xi^\vep$, where $\xi = \exp(2\pi\sqrt{-1}/3)$,
$\vep = 0,1,2$, is a simple $V_L^\tau$-module, while if $U$ belongs
to a $\tau$-orbit of length $3$, then $U$ itself is a simple
$V_L^\tau$-module and the three members in the $\tau$-orbit are
equivalent each other (cf. \cite[Theorem 6.14]{DY}). Among the $12$
inequivalent simple $V_L$-modules, three of them are $\tau$-stable
and the remaining nine simple $V_L$-modules are divided into three
$\tau$-orbits. In this way we obtain $12$ simple $V_L^\tau$-modules.
It is known that there are three inequivalent simple $\tau$-twisted
(resp. $\tau^2$-twisted) $V_L$-modules $V_L^{T_{\chi_j}}(\tau)$
(resp. $V_L^{T_{\chi'_j}}(\tau^2)$), $j=0,1,2$. The automorphism
$\tau$ acts on these $\tau$-twisted or $\tau^2$-twisted
$V_L$-modules and each eigenspace of $\tau$ is a simple
$V_L^\tau$-module (cf. \cite[Theorem 2]{MT}). There are $18$ such
simple $V_L^\tau$-modules. Furthermore, all of these simple
$V_L^\tau$-modules are inequivalent. Hence there are at least $30$
inequivalent simple $V_L^\tau$-modules.

The main part of our argument is to show that every simple
$V_L^\tau$-module is isomorphic to one of the $30$ above mentioned
simple $V_L^\tau$-modules. Recall that $V_L^\tau = M^0 \op W^0$ and
$M^0$ is a tensor product of two $W_3$ algebras. The $W_3$ algebra
of central charge $6/5$ (resp. $4/5$) possesses $20$ (resp. $6$)
inequivalent simple modules. Thus there are $120$ inequivalent
simple $M^0$-modules. It turns out that among these simple
$M^0$-modules, $60$ of them can not appear as an $M^0$-submodule in
any simple $V_L^\tau$-module and that each simple $V_L^\tau$-module
is a direct sum of two of the remaining $60$ simple $M^0$-modules.
We note that $W^0$ is not a simple current $M^0$-module, namely,
$V_L^\tau$ is a nonsimple current extension of $M^0$. A discussion
on simple modules for another nonsimple current extension of a
certain vertex operator algebra can be found in \cite[Appendix
C]{LYY}.

The organization of this paper is as follows. In Section 2 we review
various notions about untwisted or twisted modules for vertex
operator algebras, together with some basic tools which will be used
in later sections. In Section 3 we fix notation for the vertex
operator algebra $V_L^\tau$ and collect its properties. We clarify
an argument on the simplicity of $M^0_T(\tau^i)$ and
$W^0_T(\tau^i)$, $i=1,2$, in \cite[Proposition 6.8]{KLY2}.
Furthermore, we correct some misprints in \cite[(6.46)]{KLY2} and in
an equation of \cite[page 265]{DLTYY} concerning a decomposition of
a simple $\tau$-twisted $V_L$-module $V_L^{T_{\chi_j}}(\tau)$,
$j=1,2$ as a $\tau$-twisted $M_k^0 \ot M_t^0$-module (see Remark
\ref{rmk:correction}). In Section 4 we discuss the structure of the
$30$ known simple $V_L^\tau$-modules. In particular, we calculate
the action of $o(\tom^1)$, $o(\tom^2)$, $o(J)$, $o(K)$, and $o(P)$
on the top level of these simple modules. Finally, in Section 5 we
complete the classification of simple $V_L^\tau$-modules. We also
show the rationality of $V_L^\tau$.

The authors would like to thank Ching Hung Lam, Masahiko Miyamoto,
and Hiroshi Yamauchi for valuable discussions. The proof of Lemma
5.7 is essentially the same as that of \cite[Lemma C.3]{LYY}. Part
of our calculation was done by a computer algebra system
Risa/Asir. The authors are grateful to Kazuhiro Yokoyama for
helpful advice on computer programs.

\section{Preliminaries}
We recall some notation for untwisted or twisted modules for a
vertex operator algebra. We also review the twisted version of Zhu's
theory. A basic reference to twisted modules is \cite{DLM1}. For
untwisted modules, see also \cite{LL}. Let $(V,Y,\1,\om)$ be a
vertex operator algebra and $g$ be an automorphism of $V$ of finite
order $T$. Set $V^r = \{ v \in V\,|\, gv = e^{2\pi\sqrt{-1}r/T}
v\}$, so that $V = \op_{r \in \Z/T\Z} V^r$.

\begin{definition}\label{def:weak-twisted}
A weak g-twisted $V$-module $M$ is a vector space equipped with a
linear map
\begin{equation*}
Y_M(\,\cdot\,,z) : v \in V \mapsto Y_M(v,z) = \sum_{n \in \Q}
v_nz^{-n-1} \in (\End M)\{z\}
\end{equation*}
which satisfies the following conditions.

$(1)$ $Y_M(v,z) = \sum_{n \in r/T+\Z}v_n z^{-n-1}$ for $v \in V^r$.

$(2)$ $v_n w = 0$ if $n \gg 0$, where $v \in V$ and $w \in M$.

$(3)$ $Y_M(\1,z) = \mathrm{id}_M$.

$(4)$ For $u \in V^r$ and $v \in V$, the $g$-twisted Jacobi identity
holds.
\begin{equation}\label{eq:twisted-Jacobi}
\begin{split}
& z_0^{-1}\delta\big(\frac{z_1-z_2}{z_0}\big) Y_M(u,z_1)Y_M(v,z_2)-
z_0^{-1}\delta\big(\frac{z_2-z_1}{-z_0}\big)
Y_M(v,z_2)Y_M(u,z_1)\\
& \qquad = z_2^{-1}\big(\frac{z_1-z_0}{z_2}\big)^{-r/T}
\delta\big(\frac{z_1-z_0}{z_2}\big) Y_M(Y(u,z_0)v, z_2).
\end{split}
\end{equation}
\end{definition}

Compare the coefficients of $z_0^{-l-1}z_1^{-m-1}z_2^{-n-1}$ in both
sides of \eqref{eq:twisted-Jacobi} for $u \in V^r$, $v \in V^s$, $l
\in \Z$, $m\in \frac{r}{T}+\Z$, and  $n \in \frac{s}{T}+\Z$. Then we
obtain the following identity.
\begin{equation}\label{eq:twisted-Borcherds}
\sum_{i=0}^{\infty}\binom{m}{i} (u_{l+i}v)_{m+n-i} =
\sum_{i=0}^{\infty} (-1)^i \binom{l}{i}
\big(u_{l+m-i}v_{n+i}-(-1)^l v_{l+n-i}u_{m+i}\big).
\end{equation}
In the case $l=0$, \eqref{eq:twisted-Borcherds} reduces to
\begin{equation}\label{eq:twisted-commutator2}
[u_m,v_n] = \sum_{i=0}^\infty \binom{m}{i} (u_i v)_{m+n-i}.
\end{equation}

The Virasoro element $\om$ is contained in $V^0$. Let $L(n) =
\om_{n+1}$ for $n \in \Z$. Then
\begin{gather*}
[L(m), L(n)] = (m-n)L(m+n) + \frac{1}{12}(m^3-m) \delta_{m+n, 0}
(\mathrm{rank }V),\\
\frac{d}{dz}Y_M(v,z) = Y_M(L(-1)v,z)
\end{gather*}
for $v \in V$ \cite[(3.8), (3.9)]{DLM1}.

An important consequence of \eqref{eq:twisted-Jacobi} is the
associativity formula \cite[(3.5)]{DLM1}
\begin{equation}\label{eq:twisted-associativity}
(z_0 + z_2)^{k+r/T} Y_M(u, z_0+z_2)Y_M(v,z_2)w = (z_2 +
z_0)^{k+r/T}Y_M(Y(u,z_0)v,z_2)w,
\end{equation}
where $u \in V^r$, $v \in V$, $w \in M$, and $k$ is a nonnegative
integer such that $z^{k+r/T}Y_M(u,z)w \in M[[z]]$.

Let $(M, Y_M)$ and $(N, Y_N)$ be weak $g$-twisted $V$-modules. A
homomorphism of $M$ to $N$ is a linear map $f : M \rightarrow N$
such that $f Y_M(v,z) = Y_N(v,z) f$ for all $v \in V$.

Let $\N$ be the set of nonnegative integers.

\begin{definition}\label{def:admissible-twisted}
A $\frac{1}{T}\N$-graded weak $g$-twisted $V$-module $M$ is a weak
$g$-twisted $V$-module with a $\frac{1}{T}\N$-grading $M = \op_{n
\in \frac{1}{T}\N} M_{(n)}$ such that
\begin{equation}\label{eq:admissible-grading}
v_m M_{(n)} \subset M_{(n + \wt(v) - m -1)}
\end{equation}
for any homogeneous vectors $v \in V$.
\end{definition}

A $\frac{1}{T}\N$-graded weak $g$-twisted $V$-module here is called
an admissible $g$-twisted $V$-module in \cite{DLM1}. Without loss we
can shift the grading of a $\frac{1}{T}\N$-graded weak $g$-twisted
$V$-module $M$ so that $M_{(0)} \ne 0$ if $M \ne 0$. We call such an
$M_{(0)}$ the top level of $M$.

\begin{definition}\label{def:twisted-module}
A $g$-twisted $V$-module $M$ is a weak $g$-twisted $V$-module with a
$\C$-grading $M = \op_{\lambda \in \C} M_\lambda$, where $M_\lambda
= \{ w \in M\,|\,L(0) w = \lambda w\}$. Moreover, each $M_\lambda$
is a finite dimensional space and for any fixed $\lambda$,
$M_{\lambda + n/T} = 0$ for all sufficiently small integers $n$.
\end{definition}

A $g$-twisted $V$-module is sometimes called an ordinary $g$-twisted
$V$-module. By \cite[Lemma 3.4]{DLM1}, any $g$-twisted $V$-module is
a $\frac{1}{T}\N$-graded weak $g$-twisted $V$-module. Indeed, assume
that $M$ is a $g$-twisted $V$-module. For each $\lambda \in \C$ with
$M_\lambda \ne 0$, let $\lambda_0=\lambda + m/T$ be such that $m \in
\Z$ is minimal subject to $M_{\lambda_0} \ne 0$. Let $\Lambda$ be
the set of all such $\lambda_0$ and let $M_{(n)} = \op_{\lambda \in
\Lambda} M_{n+\lambda}$. Then $M_{(n)}$ satisfies the condition in
Definition \ref{def:admissible-twisted}. Thus we have the following
inclusions.
\begin{align*}
\{ g\text{-twisted }V\text{-modules} \} & \subset \{
\textstyle{\frac{1}{T}\N}\text{-graded
weak }g\text{-twisted }V\text{-modules}\}\\
& \subset \{\text{weak }g\text{-twisted }V\text{-modules}\}
\end{align*}

\begin{definition}\label{def:g-rational}
A vertex operator algebra $V$ is said to be $g$-rational if every
$\frac{1}{T}\N$-graded weak $g$-twisted $V$-module is semisimple,
that is, a direct sum of simple $\frac{1}{T}\N$-graded weak
$g$-twisted $V$-modules.
\end{definition}

Let $M$ be a weak $g$-twisted $V$-module. The next lemma is a
twisted version of \cite[Lemma 3.12]{Li3}. In fact, using the
associativity formula \eqref{eq:twisted-associativity} we can prove
it by essentially the same argument as in \cite{Li3}.

\begin{lem}\label{lem:Formula1}
Let $u \in V^r$, $v \in V^s$, $w \in M$, and $k$ be a nonnegative
integer such that $z^{k+r/T}Y_M(u,z)w \in M[[z]]$. Let $p \in
\frac{r}{T}+\Z$, $q \in \frac{s}{T}+\Z$, and $N$ be a nonnegative
integer such that $z^{N+1+q}Y_M(v,z)w \in M[[z]]$. Then
\begin{equation}\label{eq:Formula1}
u_pv_qw = \sum_{i=0}^N \sum_{j=0}^\infty \binom{p-k-r/T}{i}
\binom{k+r/T}{j} (u_{p-k-r/T-i+j} v)_{q+k+r/T+i-j}w.
\end{equation}
\end{lem}

Conversely, $(u_p v)_q w$ can be written as a linear combination
of some vectors of the form $u_i v_j w$.

\begin{lem}\label{lem:Formula2}
Let $u \in V^r$, $v \in V^s$, $w \in M$. Then for $p \in \Z$ and
$q \in \frac{r+s}{T}+\Z$, the vector $(u_p v)_q w$ is a linear
combination of $u_i v_j w$ with $i \in \frac{r}{T}+\Z$ and $j \in
\frac{s}{T}+\Z$.
\end{lem}

\begin{proof}
Let $X = \spn \{ u_i v_j w\,|\,i \in \frac{r}{T}+\Z, j \in
\frac{s}{T}+\Z\}$. We use \eqref{eq:twisted-Borcherds}. Take $m \in
\frac{r}{T}+\Z$ such that $u_{m+i}w = 0$ for $i \ge 0$. Let $N \in
\Z$ be such that $u_{N+i} v = 0$ for $i > 0$. If $p > N$, then $u_p
v = 0$ and the assertion is trivial. Assume that $p \le N$. For
$j=0,1, \ldots, N-p$, let $l = p+j$ and $n=q-m-j$ in
\eqref{eq:twisted-Borcherds}. Then
\begin{equation*}
\sum_{i=0}^\infty \binom{m}{i} (u_{p+j+i} v)_{q-j-i} w =
\sum_{i=0}^\infty (-1)^i \binom{p+j}{i} u_{p+m+j-i} v_{q-m-j+i} w.
\end{equation*}
The right hand side of this equation is contained in $X$. Consider
the left hand side for each of $j=N-p, N-p-1, \ldots, 1,0$. Then we
see that $(u_N v)_{q-N+p}w \in X$, $(u_{N-1} v)_{q-N+p+1}w \in X$,
$\ldots$ , and $(u_p v)_q w \in X$.
\end{proof}

For subsets $A$, $B$ of $V$ and a subset $X$ of $M$, set $A \cdot X
= \spn\{ u_n w\,|\,u \in A, w \in X, n \in \frac{1}{T}\Z\}$ and $A
\cdot B = \spn\{ u_n v\,|\,u \in A, v \in B, n \in \Z\}$. Then it
follows from \eqref{eq:Formula1} that $A \cdot (B \cdot X) \subset
(A \cdot B)\cdot X$ (see also \cite[(2.2)]{Yamauchi}). For a vector
$w \in M$, this in particular implies that $V \cdot w$ is a weak
$g$-twisted $V$-submodule of $M$. If $w$ is an eigenvector for
$L(0)$, then $V \cdot w$ is a direct sum of eigenspaces for $L(0)$.
Each eigenspace is not necessarily of finite dimension. Thus $V
\cdot w$ is not a $g$-twisted module in general. This subject was
discussed in \cite{ABD, Buhl, Yamauchi}. We will review it later in
this section.

In \cite{Z}, Zhu introduced an associative algebra $A(V)$ called the
Zhu algebra for a vertex operator algebra $V$, which plays a crucial
role in the study of representations for $V$. Later, Dong, Li and
Mason \cite{DLM1} constructed an associative algebra $A_g(V)$ called
the $g$-twisted Zhu algebra in order to generalize Zhu's theory to
$g$-twisted representations for $V$. The definition of $A_g(V)$ is
similar to that of $A(V)$. Let $V$, $g$, $T$, and $V^r$ be as
before. Roughly speaking, $A_g(V)=V/O_g(V)$ is a quotient space of
$V$ with a binary operation $\ast_g$. It is in fact an associative
algebra with respect to $\ast_g$. If $r \ne 0$, then $V^r \subset
O_g(V)$. Thus $A_g(V) = (V^0 + O_g(V))/O_g(V)$. For the case $g=1$,
see \eqref{eq:Zhu-Operation} in Section \ref{sec:Classification}.

A certain Lie algebra $V[g]$ was considered in \cite{DLM1}. Any weak
$g$-twisted $V$-module is a module for the Lie algebra $V[g]$ (cf.
\cite[Lemma 5.1]{DLM1}). Moreover, for a $V[g]$-module $M$, the
space $\Omega(M)$ of lowest weight vectors with respect to $V[g]$
was defined. If $M$ is a weak $g$-twisted $V$-modules, then
$\Omega(M)$ is the set of $w \in M$ such that $v_{\wt(v) -1+n}w=0$
for all homogeneous vectors $v \in V$ and $0 < n \in \frac{1}{T}\Z$.
The map $v \mapsto o(v)$ for homogeneous vectors $v \in V^0$ induces
a representation of the associative algebra $A_g(V)$ on $\Omega(M)$,
where $o(v)= v_{\wt(v) -1}$. If $M$ is a $\frac{1}{T}\N$-graded weak
$g$-twisted $V$-module, then the top level $M_{(0)}$ is contained in
$\Omega(M)$. In the case where $M$ is a simple
$\frac{1}{T}\N$-graded weak $g$-twisted $V$-module, $M_{(0)} =
\Omega (M)$ and $M_{(0)}$ is a simple $A_g(V)$-module (cf.
\cite[Proposition 5.4]{DLM1}).

For any $A_g(V)$-module $U$, a certain $\frac{1}{T}\N$-graded
$V[g]$-module $M(U)$ such that $M(U)_{(0)} = U$ was defined (cf.
\cite[(6.1)]{DLM1}). Let $W$ be the subspace of $M(U)$ spanned by
the coefficients of
\begin{align*}
& (z_0+z_2)^{\wt(u) -1 + \delta_r + r/T}Y_M(u,z_0+z_2)Y_M(v,z_2)w\\
& \qquad\qquad - (z_2+z_0)^{\wt(u) -1 + \delta_r +
r/T}Y_M(Y(u,z_0)v, z_2)w
\end{align*}
for all homogeneous $u \in V^r$, $v \in V$, $w \in U$ (cf.
\cite[(6.3)]{DLM1}). Set $\bar{M}(U) = M(U)/U(V[g])W$, which is a
quotient module of $M(U)$ by the $V[g]$-submodule generated by $W$.

The following results will be necessary in Sections \ref{sec:VL-tau}
and \ref{sec:Classification}.

\begin{thm}$($\cite[Theorem 6.2]{DLM1}$)$\label{thm:Ag-universal}
$\bar{M}(U)$ is a $\frac{1}{T}\N$-graded weak $g$-twisted
$V$-module such that its top level $\bar{M}(U)_{(0)}$ is equal to
$U$ and such that it has the following universal property: for any
weak $g$-twisted $V$-module $M$ and any homomorphism $\varphi : U
\rightarrow \Omega(M)$ of $A_g(V)$-modules, there is a unique
homomorphism $\bar{\varphi} : \bar{M}(U) \rightarrow M$ of weak
$g$-twisted $V$-modules which is an extension of $\varphi$.
\end{thm}

Let $J$ be the sum of all $\frac{1}{T}\N$-graded $V[g]$-submodules
of $M(U)$ which intersect trivially with $U$. Since $M(U)_{(0)} =
U$, it is a unique $\frac{1}{T}\N$-graded $V[g]$-submodule of $M(U)$
being maximal subject to $J \cap U = 0$. The principal point is that
$U(V[g])W \subset J$. Set $L(U) = M(U)/J$.

\begin{thm}$($\cite[Theorem 6.3]{DLM1}$)$\label{thm:LU}
$L(U)$ is a $\frac{1}{T}\N$-graded weak $g$-twisted $V$-module
such that $\Omega(L(U)) \cong U$ as $A_g(V)$-modules.
\end{thm}

\begin{rmk}\label{rmk:Ag-universal}
If $M$ is a $\frac{1}{T}\N$-graded weak $g$-twisted $V$-module and
$\varphi : U \rightarrow M_{(0)}$ is a homomorphism of
$A_g(V)$-modules, then the homomorphism $\bar{\varphi} : \bar{M}(U)
\rightarrow M$ of weak $g$-twisted $V$-modules in Theorem
\ref{thm:Ag-universal} preserves the $\frac{1}{T}\N$-grading.
Indeed, $\bar{M}(U) = \spn\{ v_n U\,|\, v \in V, n \in
\frac{1}{T}\Z\}$ by \eqref{eq:Formula1}, since $\bar{M}(U)$ is
generated by $U$ as a $\frac{1}{T}\N$-graded weak $g$-twisted
$V$-module. By \eqref{eq:admissible-grading}, $v_{\wt(v)-1-n}
\bar{M}(U)_{(0)} \subset \bar{M}(U)_{(n)}$ for any homogeneous $v
\in V$ and $n \in \frac{1}{T}\Z$. Since $\bar{M}(U)_{(0)} = U$, it
follows that $\bar{M}(U)_{(n)}$ is spanned by $v_{\wt(v)-1-n}U$ for
all homogeneous $ v \in V$. Now, $\bar{\varphi}(v_{\wt(v)-1-n}U) =
v_{\wt(v)-1-n}\bar{\varphi}(U)$ is contained in
$v_{\wt(v)-1-n}M_{(0)} \subset M_{(n)}$. Hence
$\bar{\varphi}(\bar{M}(U)_{(n)}) \subset M_{(n)}$ as required. In
the case where both of $\bar{M}(U)$ and $M$ are ordinary $g$-twisted
$V$-modules, $\bar{\varphi}$ becomes a homomorphism of ordinary
$g$-twisted $V$-modules since $\bar{\varphi}$ commutes with $L(0)$.
\end{rmk}

\begin{lem}\label{lem:surjective-LU}
Let $U$ be an $A_g(V)$-module. Let $S$ be a $\frac{1}{T}\N$-graded
weak $g$-twisted $V$-module such that it is generated by its top
level $S_{(0)}$ and such that $S_{(0)}$ is isomorphic to $U$ as an
$A_g(V)$-module. Then there is a surjective homomorphism $S
\rightarrow L(U)$ of weak $g$-twisted $V$-modules which preserves
the $\frac{1}{T}\N$-grading.
\end{lem}

\begin{proof}
By Theorem \ref{thm:Ag-universal} and Remark \ref{rmk:Ag-universal},
an isomorphism $\varphi : U \rightarrow S_{(0)}$ of $A_g(U)$-modules
can be extended to a surjective homomorphism $\bar{\varphi}:
\bar{M}(U) \rightarrow S$ of weak $g$-twisted $V$-modules which
preserves the $\frac{1}{T}\N$-grading. The kernel $\Ker
\bar{\varphi}$ of $\bar{\varphi}$ intersects trivially with
$\bar{M}(U)_{(0)}$ and so it is contained in $\op_{0 < n \in
\frac{1}{T}\N} \bar{M}(U)_{(n)}$. Let $I$ be a
$\frac{1}{T}\N$-graded $V[g]$-submodule of $M(U)$ such that $\Ker
\bar{\varphi} = I/U(V[g])W$. Then $I \cap U = 0$. This implies that
$I \subset J$. Hence $L(U) = M(U)/J$ is a homomorphic image of
$M(U)/I \cong S$.
\end{proof}

\begin{thm}$($\cite[Theorem 7.2]{DLM1}$)$\label{thm:Ag-simple}
$L$ is a functor  from the category of simple $A_g(V)$-modules to
the category of simple $\frac{1}{T}\N$-graded weak $g$-twisted
$V$-modules such that $\Omega \circ L = \mathrm{id}$ and $L \circ
\Omega = \mathrm{id}$.
\end{thm}

\begin{thm}$($\cite[Theorem 8.1]{DLM1}$)$\label{thm:Ag-property}
If $V$ is a $g$-rational vertex operator algebra, then the
following assertions hold.

$(1)$ $A_g(V)$ is a finite dimensional semisimple associative
algebra.

$(2)$ $V$ has only finitely many isomorphism classes of simple
$\frac{1}{T}\N$-graded weak $g$-twisted $V$-modules.

$(3)$ Every simple $\frac{1}{T}\N$-graded weak $g$-twisted
$V$-module is an ordinary $g$-twisted $V$-modules.
\end{thm}

In case of $g=1$, the above argument reduces to the untwisted case.
In particular, $A_g(V)$ is identical with the original Zhu algebra
$A(V)$ if $g=1$.

There is an important intrinsic property of a vertex operator
algebra, namely, the $C_2$-cofiniteness. Let $C_2(V)=\spn \{ u_{-2}v
\,|\, u,v \in V \}$. More generally, we set $C_2(M) = \spn \{
u_{-2}w \,|\, u \in V, w \in M \}$ for a weak $V$-module $M$. If the
dimension of the quotient space $V/C_2(V)$ is finite, $V$ is said to
be $C_2$-cofinite. Similarly, a weak $V$-module $M$ is said to be
$C_2$-cofinite if $M/C_2(M)$ is of finite dimension. The notion of
$C_2$-cofiniteness of a vertex operator algebra was first introduced
by Zhu \cite{Z}. The subspace $C_2(M)$ of a weak $V$-module $M$ was
studied in Li \cite{Li1}. We refer the reader to \cite{NT} also.

\begin{thm}$($\cite[Proposition 3.6]{DLM2}$)$
\label{thm:Ag-finite-dim} If $V$ is $C_2$-cofinite, then $A_g(V)$
is of finite dimension.
\end{thm}

If $V=\op_{n=0}^\infty V_n$ and $V_0 = \C\1$, then $V$ is said to
be of CFT type. Here $V_n$ denotes the homogeneous subspace of
weight $n$, that is, the eigenspace of $L(0)=\om_1$ with
eigenvalue $n$.

\begin{thm}$($\cite[Lemma 3.3]{Yamauchi}$)$
\label{thm:spanninng-set} Suppose $V$ is $C_2$-cofinite and of CFT
type. Choose a finite dimensional $L(0)$-invariant and $g$-invariant
subspace $U$ of $V$ such that $V = U + C_2(V)$. Let $W$ be a weak
$g$-twisted $V$-module generated by a vector $w$. Then $W$ is
spanned by the vectors of the form $u^1_{-n_1}u^2_{-n_2} \cdots
u^k_{-n_k}w$ with $n_1 > n_2 > \cdots > n_k > -N$ and $u^i \in U$,
$i=1,2, \ldots, k$, where $N \in \frac{1}{T}\Z$ is a constant such
that $u_m w = 0$ for all $u \in U$ and $m \ge N$.
\end{thm}

\begin{thm}$($\cite[Corollaries 3.8 and 3.9]{Yamauchi}$)$
\label{thm:C2-and-CFT} Suppose $V$ is $C_2$-cofinite and of CFT
type. Then the following assertions hold.

$(1)$ Every weak $g$-twisted $V$-module is a
$\frac{1}{T}\N$-graded weak $g$-twisted $V$-module.

$(2)$ Every simple weak $g$-twisted $V$-module is a simple ordinary
$g$-twisted $V$-module.
\end{thm}

\begin{rmk}\label{rmk:ordinary-submodule}
Suppose $V$ is $C_2$-cofinite and of CFT type. Let $M$ be a weak
$g$-twisted $V$-module and $w^1$, $\ldots$, $w^k$ be eigenvectors of
$L(0)$ in $M$. Then the weak $g$-twisted $V$-submodule $W$ generated
by $w^1$, $\ldots$, $w^k$ is an ordinary $g$-twisted $V$-module.
Indeed, $W$ is a direct sum of eigenspaces for $L(0)$ and each
homogeneous subspace is of finite dimension by Theorem
\ref{thm:spanninng-set}.
\end{rmk}

For the untwisted case, that is, the case $g=1$, we refer the reader
to \cite{ABD, Buhl, DLM0, Li1}. A spanning set for a vertex operator
algebra was first studied in \cite[Proposition 8]{GN}.

\section{Fixed point subalgebra
$(V_{\sqrt{2}A_2})^\tau$}\label{sec:VL-tau} In this section we fix
notation. We tend to follow the notation in \cite{DLTYY, KLY1, KLY2}
unless otherwise specified. We also recall certain properties of the
lattice vertex operator algebra $V_{\sqrt{2}A_2}$ associated with
$\sqrt{2}$ times an ordinary root lattice of type $A_2$ and its
subalgebras (cf. \cite{DLTYY, KLY1, KLY2, KMY}).

Let $\al_1, \al_2$ be the simple roots of type $A_2$ and set
$\al_0=-(\al_1+\al_2)$. Thus $\la \al_i, \al_i\ra = 2$ and $\la
\al_i, \al_j\ra=-1$ if $i \ne j$. Set $\be_i=\sqrt{2}\al_i$ and let
$L=\Z\be_1+\Z\be_2$ be the lattice spanned by $\be_1$ and $\be_2$.
We denote the cosets of $L$ in its dual lattice $L^\perp = \{ \al
\in \Q \ot_{\Z} L\,|\,\la\al, L\ra \subset \Z\}$ as follows.
\begin{equation*}
L^0=L,\quad  L^1=\frac{-\be_1+\be_2}{3}+L ,\quad
L^2=\frac{\be_1-\be_2}{3}+L,
\end{equation*}
\begin{equation*}
L_0=L,\quad L_a=\frac{\be_2}{2}+L,\quad
L_b=\frac{\be_0}{2}+L,\quad L_c=\frac{\be_1 }{2}+L,
\end{equation*}
\begin{equation*}
L^{(i,j)} = L_i + L^j
\end{equation*}
for $i=0,a,b,c$ and $j=0,1,2$, where $\{0,a,b,c\} \cong \Z_2 \times
\Z_2$ is Klein's four-group. Note that $L^{(i,j)},i \in
\{0,a,b,c\},j\in \{0,1,2\}$ are all the cosets of $L$ in $ L^{\perp
}$ and $L^\perp/L \cong \Z_2 \times \Z_2 \times \Z_3$.

We adopt the standard notation for the vertex operator algebra
$(V_L,Y(\,\cdot\,,z))$ associated with the lattice $L$
(cf.\cite{FLM}). In particular, ${\mathfrak h}=\C\otimes_{\Z} L$
is an abelian Lie algebra, $\hat {\mathfrak h}={\mathfrak
h}\otimes \C[t,t^{-1}]\oplus \C c$ is the corresponding affine Lie
algebra, $M(1)=\C[\al(n)\,;\,\al\in {\mathfrak h}, n<0],$ where
$\al(n)=\al\otimes t^n,$ is the unique simple $\hat{\mathfrak
h}$-module such that $\alpha(n)1=0$ for all $\alpha\in {\mathfrak
h}$ and $n>0$ and $c=1$. As a vector space $V_L = M(1) \ot \C[L]$
and for each $v \in V_L$, a vertex operator $Y(v, z) = \sum_{n \in
\Z} v_n z^{-n-1} \in \End (V_L)[[z,z^{-1}]]$ is defined. The
vector $\1 = 1 \otimes 1$ is called the vacuum vector. In our case
$\la \al , \be \ra \in 2\Z$ for any $\al, \be \in L$. Thus the
twisted group algebra $\C\{L\}$ of \cite{FLM} is naturally
isomorphic to the ordinary group algebra $\C[L]$.

There are exactly $12$ inequivalent simple $V_L$-modules, which
are represented by $V_{L^{(i,j)}}$, $i=0,a,b,c$ and $j=0,1,2$ (cf.
\cite{D1}). We use the symbol $e^{\alpha}, \alpha \in L^\perp$ to
denote a basis of $\C\{L^\perp\}$.

We consider the following three isometries of
$(L,\la\cdot,\cdot\ra)$.
\begin{equation}\label{eq:3Aut}
\begin{split}
\tau &: \be_1 \to \be_2 \to \be_0 \to \be_1,\\
\sg &: \be_1 \to \be_2, \qquad \be_2 \to \be_1,\\
\theta &: \be_i \to -\be_i, \quad i=1,2.
\end{split}
\end{equation}
Note that $\tau$ is fixed-point-free and of order $3$. The
isometries $\tau, \sg$, and $\theta$ of $L$ can be extended linearly
to isometries of $L^\perp$. Moreover, the isometry $\tau$ lifts
naturally to an automorphism of $V_L$:
\begin{equation*}
\al^1(-n_1)\cdots\al^k(-n_k)e^{\be} \longmapsto
(\tau\al^1)(-n_1)\cdots(\tau\al^k)(-n_k)e^{\tau\be}.
\end{equation*}
By abuse of notation, we  denote it by $\tau$ also. We can consider
the action of $\tau$ on $V_{L^{(i,j)}}$ in a similar way. We apply
the same argument to $\sigma$ and $\theta$. Our purpose is the
classification of simple modules for the fixed point subalgebra
$V_L^\tau = \{ v \in V_L\,|\, \tau v = v\}$ of $V_L$ by the
automorphism $\tau$.

For a simple $V_L$-module $(U, Y_U)$, let $(U\circ\tau,
Y_{U\circ\tau})$ be a new $V_L$-module such that $U\circ\tau = U$ as
vector spaces and $Y_{U\circ\tau}(v,z)=Y_U(\tau v,z)$ for $v \in
V_L$ (cf. \cite{DLM2}). Then $U \mapsto U\circ\tau$ induces a
permutation on the set of simple $V_L$-modules. If $U$ and
$U\circ\tau$ are equivalent $V_L$-modules, $U$ is said to be
$\tau$-stable. The following lemma is a straightforward consequence
of the definition of $V_{L^{(i,j)}}$.

\begin{lem}\label{STABLE}
$(1)$ $V_{L^{(0,j)}}$, $j = 0,1,2$ are $\tau$-stable.

$(2)$ $V_{L^{(a,j)}}\circ\tau = V_{L^{(c,j)}}$,
$V_{L^{(c,j)}}\circ\tau = V_{L^{(b,j)}}$, and
$V_{L^{(b,j)}}\circ\tau = V_{L^{(a,j)}}$, $j = 0,1,2$.
\end{lem}

A family of simple twisted modules for lattice vertex operator
algebras was constructed in \cite{DL2, L}. Following \cite{DL2},
three inequivalent simple $\tau$-twisted $V_L$-modules
$(V_L^{T_{\chi_j}}(\tau), Y^\tau(\cdot,z))$, $j=0,1,2$ were studied
in \cite{DLTYY, KLY2}. By the above lemma and \cite[Theorem
10.2]{DLM2}, we know that $(V_L^{T_{\chi_j}}(\tau),
Y^\tau(\cdot,z))$, $j=0,1,2$, are all the inequivalent simple
$\tau$-twisted $V_L$-modules. Similarly, there are exactly three
inequivalent simple $\tau^2$-twisted $V_L$-modules
$(V_L^{T_{\chi'_j}}(\tau^2), Y^{\tau^2}(\cdot,z))$, $j=0,1,2$.

We use the same notation for $(V_L^{T_{\chi_j}}(\tau),
Y^\tau(\cdot,z))$ and $(V_L^{T_{\chi'_j}}(\tau^2),
Y^{\tau^2}(\cdot,z))$ as in \cite[Section 4]{DLTYY}. Thus
\begin{equation*}
V_L^{T_{\chi_j}}(\tau)=S[\tau]\otimes T_{\chi_j},
\end{equation*}
where $T_{\chi_j}$, $j =0,1,2$ are the one-dimensional
representations of a certain central extention of $L$
affording the character $\chi_j$. Let
\begin{equation*}
h_1 = \frac{1}{3}(\be_1 + \xi^2 \be_2 + \xi \be_0), \qquad h_2=
\frac{1}{3}(\be_1 + \xi \be_2 + \xi^2 \be_0).
\end{equation*}
Then $\tau h_i = \xi^i h_i$, $\la h_1, h_1 \ra = \la h_2, h_2 \ra =
0$, and $\la h_1, h_2 \ra=2$. Moreover,
$\be_i=\xi^{i-1}h_1+\xi^{2(i-1)}h_2$, $i=0,1,2$. As a vector space,
$S[\tau]$ is isomorphic to a polynomial algebra with variables
$h_1(1/3 + n)$, $h_2(2/3 + n)$, $n \in \Z_{< 0}$. The isometry
$\tau$ acts on $S[\tau]$ by $\tau h_j = \xi^j h_j$. We define the
action of $\tau$ on $T_{\chi_j}$ to be the identity. The weight in
$S[\tau]$ is given by $\wt h_i(i/3 + n) = -i/3 - n$, $i=1,2$ and
$\wt 1 = 1/9$. The weight of any element of $T_{\chi_j}$ is defined
to be $0$. Note that the weight in $V_L^{T_{\chi_j}}(\tau)$ is
identical with the eigenvalue for the action of the coefficient of
$z^{-2}$ in the $\tau$-twisted vertex operator $Y^\tau (\om, z)$,
where $\om$ denotes the Virasoro element of $V_L$.

The simple $\tau^2$-twisted
$V_L$-modules $(V_L^{T_{\chi_j}}(\tau^2),
Y^{\tau^2}(\cdot,z))$, $j=0,1,2$ is
\begin{equation*}
V_L^{T_{\chi'_j}}(\tau^2)=S[\tau^2]\otimes T_{\chi'_j},
\end{equation*}
where $T_{\chi'_j}$, $j =0,1,2$ are the one-dimensional
representations of a certain central extension of $L$ affording the
character $\chi'_j$. Moreover, $S[\tau^2]$ is isomorphic to a
polynomial algebra with variables $h'_1(1/3 + n)$, $h'_2(2/3 + n)$,
$n \in \Z_{< 0}$ as a vector space, where $h'_1 = h_2$ and $h'_2 =
h_1$. Thus $\tau^2 h'_i = \xi^i h'_i$, $i=1,2$. The action of $\tau$
on $S[\tau^2]$ is given by $\tau h'_i = \xi^{2i} h'_i$, $i=1,2$. The
action of $\tau$ on $T_{\chi'_j}$ is defined to be the identity. The
weight in $S[\tau^2]$ is given by $\wt h'_i(i/3 + n) = -i/3 - n$,
$i=1,2$ and $\wt 1 = 1/9$. The weight of any element of
$T_{\chi'_j}$ is defined to be $0$. The weight in
$V_L^{T_{\chi'_j}}(\tau^2)$ is identical with the eigenvalue for the
action of the coefficient of $z^{-2}$ in the $\tau^2$-twisted vertex
operator $Y^{\tau^2} (\om, z)$.

By Lemma \ref{STABLE}, \cite[Theorem 4.4]{DM}, and \cite[Theorem
6.14]{DY},
\begin{equation*}
V_{L^{(0,j)}} (\vep) = \{ v \in V_{L^{(0,j)}} \,|\,
\tau v = \xi^\vep v \}, \quad j,\vep = 0, 1, 2
\end{equation*}
are inequivalent simple $V_L^\tau$-modules. For each of $j = 0, 1,
2$, we have that $V_{L^{(i,j)}}$, $i = a,b,c$ are equivalent simple
$V_L^\tau$-modules. Moreover, $V_{L^{(c,j)}}$, $j = 0,1,2$ are
inequivalent simple $V_L^\tau$-modules. From \cite[Theorem 2]{MT},
it follows that
\begin{equation*}
V_L^{T_{\chi_j}}(\tau)(\vep) = \{ v \in V_L^{T_{\chi_j}}(\tau)
\,|\, \tau v = \xi^\vep v \}, \quad j,\vep = 0, 1, 2
\end{equation*}
are inequivalent simple $V_L^\tau$-modules. Similar assertions hold
for simple $\tau^2$-twisted modules, namely,
\begin{equation*}
V_L^{T_{\chi'_j}}(\tau^2)(\vep) =
\{ v \in V_L^{T_{\chi'_j}}(\tau^2)
\,|\, \tau^2 v = \xi^\vep v \}, \quad j,\vep = 0, 1, 2
\end{equation*}
are inequivalent simple $V_L^\tau$-modules. In this way we obtain
$30$ simple $V_L^\tau$-modules. These $30$ simple $V_L^\tau$-modules
are inequivalent by \cite[Theorem 2]{MT}. We summarize the result as
follows.

\begin{lem}\label{lem:30SIMPLES}
The following $30$ simple $V_L^\tau$-modules are inequivalent.

$(1)$ $V_{L^{(0,j)}} (\vep)$, $j, \vep = 0,1,2$,

$(2)$ $V_{L^{(c,j)}}$, $j = 0,1,2$,

$(3)$ $V_L^{T_{\chi_j}}(\tau)(\vep)$, $j, \vep = 0,1,2$,

$(4)$ $V_L^{T_{\chi'_j}}(\tau^2)(\vep)$, $j, \vep = 0,1,2$.
\end{lem}

We consider the structure of $V_L^\tau$ in detail. Let
\begin{equation*}
x(\al)=e^{\sqrt{2}\al} + e^{-\sqrt{2}\al},\qquad
y(\al)=e^{\sqrt{2}\al} - e^{-\sqrt{2}\al},\qquad
w(\al)=\frac{1}{2}\al(-1)^2 - x(\al)
\end{equation*}
for $\al \in\{\pm\al_0, \pm\al_1, \pm\al_2\}$ and let
\begin{gather*}
\om = \frac{1}{6}\big( \al_1(-1)^2 + \al_2(-1)^2 +
\al_0(-1)^2\big),\\
\tom^1 = \frac{1}{5}\big( w(\al_1) + w(\al_2) + w(\al_0)\big),\qquad
\tom^2 = \om - \tom^1,\\
\om^1 = \frac{1}{4} w(\al_1), \qquad \om^2 = \tom^1 - \om^1.
\end{gather*}
Then $\om$ is the Virasoro element of $V_L$ and $\tom^1$, $\tom^2$
are mutually orthogonal conformal vectors of central charge $6/5,
4/5$ respectively. The subalgebra $\Vir(\tom^i)$ generated by
$\tom^i$ is isomorphic to the Virasoro vertex operator algebra of
given central charge, namely, $\Vir(\tom^1) \cong L(6/5,0)$ and
$\Vir(\tom^2) \cong L(4/5,0)$. Moreover, $\tom^1$ is a sum of two
conformal vectors $\om^1$ and $\om^2$ of central charge $1/2$ and
$7/10$ respectively and $\om^1$, $\om^2$ and $\tom^2$ are mutually
orthogonal. Note that $\tom^2$ was denoted by $\om^3$ in
\cite{DLTYY, KLY1, KLY2, KMY}. Such a decomposition of the Virasoro
element of a lattice vertex operator algebra into a sum of mutually
orthogonal conformal vectors was first studied in \cite{DLMN}.

Set
\begin{align*}
M_{k}^{i}&=\{ v\in V_{L_i}\,|\,(\tom^2)_1 v=0\},\\
W_{k}^{i}&=\{ v\in V_{L_i}\,|\,(\tom^2)_1 v= \frac{2}{5}v\} ,\quad
i=0,a,b,c,
\end{align*}
\begin{align*}
M_{t}^{j}&=\{ v\in V_{L^j}\,|\,(\omega^1)_1 v=
(\omega^2)_1 v=0\},\\
W_{t}^{j}&=\{ v\in V_{L^j}\,|\,(\omega^1)_1 v=0,\, (\omega^2)_1
v=\frac{3}{5}v\} ,\quad j=0,1,2.
\end{align*}
Then $M_k^0$ and $M_t^0$ are simple vertex operator algebras.
Moreover, $\{M_k^i,\,W_k^i;\, i=0,a,b,c\}$ and $\{M_t^j,\, W_t^j;\,
j=0,1,2\}$ are complete sets of representatives of isomorphism
classes of simple modules for $M_k^0$ and $M_t^0$, respectively (cf.
\cite{KLY1, KMY, LY}). As $\Vir(\om^1) \ot \Vir(\om^2)$-modules,
\begin{equation}\label{eq:MW-k}
\begin{split}
& M_{k}^{0}\cong \Big( L( \frac{1}{2},0) \otimes L( \frac{7}{10} ,0)
\Big) \oplus \Big( L( \frac{1}{2},\frac{1}{2}) \otimes L
\frac{7}{10},\frac{3}{2}) \Big),\\
& M_{k}^{a}\cong M_{k}^{b}\cong L( \frac{1}{2},\frac{1}{16})
\otimes L(\frac{7}{10},\frac{7}{16}),\\
& M_{k}^{c}\cong \Big( L( \frac{1}{2},\frac{1}{2}) \otimes
L(\frac{7}{10},0) \Big) \oplus \Big( L( \frac{1}{2},0) \otimes
L(\frac{7}{10},\frac{3}{2}) \Big),\\
& W_{k}^{0}\cong \Big( L( \frac{1}{2},0) \otimes L(
\frac{7}{10},\frac{3}{5}) \Big) \oplus \Big( L( \frac{1}{2},
\frac{1}{2}) \otimes L(\frac{7}{10},\frac{1}{10}) \Big),\\
& W_{k}^{a}\cong W_{k}^{b}\cong L( \frac{1}{2},\frac{1}{16})
\otimes L( \frac{7}{10},\frac{3}{80}),\\
& W_{k}^{c}\cong \Big( L( \frac{1}{2},\frac{1}{2}) \otimes L(
\frac{7}{10},\frac{3}{5}) \Big) \oplus \Big( L( \frac{1}{2},0)
\otimes L(\frac{7}{10},\frac{1}{10}) \Big),
\end{split}
\end{equation}
and as $\Vir(\tom^2)$-modules,
\begin{equation}\label{eq:MW-t}
\begin{split}
M_{t}^{0} \cong L( \frac{4}{5},0) \oplus L( \frac{4}{5},3), &
\qquad\qquad
M_{t}^{1} \cong M_{t}^{2}\cong L( \frac{4}{5},\frac{2}{3}),\\
W_{t}^{0} \cong L( \frac{4}{5},\frac{2}{5}) \oplus L(
\frac{4}{5},\frac{7}{5}), & \qquad\qquad  W_{t}^{1}  \cong
W_{t}^{2}\cong L(\frac{4}{5},\frac{1}{15}).
\end{split}
\end{equation}

Furthermore,
\begin{equation}\label{eq:VLij-dec}
V_{L^{(i,j)}} \cong (M_k^i\otimes M_t^j) \oplus
(W_k^i\otimes W_t^j)
\end{equation}
as $M_k^0\otimes M_t^0$-modules.
In particular,
\begin{equation}\label{eq:VL-dec}
V_L \cong ( M_k^0 \otimes M_t^0) \oplus (W_k^0 \otimes W_t^0).
\end{equation}

Note that $M_t^j = \{ v \in V_{L^j} \,|\, (\tom^1)_1 v = 0 \}$ and
that $M_k^0$, $W_k^0$ and $M_t^j$, $j=0,1,2$ are $\tau$-invariant.
However, $W_t^j$, $j=0,1,2$ are not $\tau$-invariant.

The fusion rules for $M_k^0$ and $M_t^0$ were determined in
\cite{LY} and \cite{M2}, respectively. They are as follows.
\begin{equation}\label{eq:Klein-fusion}
\begin{split}
M_k^i \times M_k^j &= M_k^{i+j},\\
M_k^i \times W_k^j &= W_k^{i+j},\\
W_k^i \times W_k^j &= M_k^{i+j} + W_k^{i+j}\\
\end{split}
\end{equation}
for $i,j = 0,a,b,c$ and
\begin{equation}\label{eq:Ternary-fusion}
\begin{split}
M_t^i \times M_t^j &= M_t^{i+j},\\
M_t^i \times W_t^j &= W_t^{i+j},\\
W_t^i \times W_t^j &= M_t^{i+j} + W_t^{i+j}\\
\end{split}
\end{equation}
for $i,j = 0,1,2$

The following two weight $3$ vectors are important.
\begin{equation*}
\begin{split}
J &= w(\al_1)_0w(\al_2) - w(\al_2)_0w(\al_1)\\
&= -\frac{1}{6}\Big(\beta_1(-2)(\beta_2 - \beta_0)(-1) +
\beta_2(-2)(\beta_0 - \beta_1)(-1) +
\beta_0(-2)(\beta_1 - \beta_2)(-1)\Big)\\
& \quad - (\beta_2 - \beta_0)(-1)y(\al_1) - (\beta_0 -
\beta_1)(-1)y(\al_2) - (\beta_1 - \beta_2)(-1)y(\al_0),
\end{split}
\end{equation*}
\begin{equation*}
\begin{split}
K &= - \frac{1}{9} (\be_1 - \be_2)(-1)(\be_2 - \be_0)(-1)
(\be_0 - \be_1)(-1)\\
& \quad + (\be_2 - \be_0)(-1)x(\al_1)
+ (\be_0 - \be_1)(-1)x(\al_2) + (\be_1 - \be_2)(-1)x(\al_0).
\end{split}
\end{equation*}

Let $M(0) = (M_k^0)^\tau = \{ u \in M_k^0 \,|\, \tau u = u \}$. The
vertex operator algebra $M(0)$ was studied in \cite{DLTYY}. Among
other things, the classification of simple modules, the rationality
and the $C_2$-cofiniteness for $M(0)$ were established. It is known
that $M(0)$ is a $W_3$ algebra of central charge $6/5$ with the
Virasoro element $\tom^1$. In fact, $M(0)$ is generated by $\tom^1$
and $J$. The following equations hold \cite[(3.1)]{DLTYY}.
\begin{equation}\label{eq:JnJ}
\begin{split}
J_5 J &= -84\cdot\1,\\
J_4 J &= 0,\\
J_3 J &= -420\tom^1,\\
J_2 J &= -210(\tom^1)_0\tom^1,\\
J_1 J &= 9(\tom^1)_0(\tom^1)_0\tom^1 - 240(\tom^1)_{-1}\tom^1,\\
J_0 J &= 22(\tom^1)_0(\tom^1)_0(\tom^1)_0\tom^1 -
120(\tom^1)_0(\tom^1)_{-1}\tom^1.
\end{split}
\end{equation}

Let $L^1(n)=(\tom^1)_{n+1}$ and $J(n)=J_{n+2}$ for $n \in \Z$, so
that the weight of these operators is $\wt L^1(n)= \wt J(n)=-n$.
Then
\begin{equation}\label{eq:L1L1}
[L^1(m),\,L^1(n)]=(m-n)L^1(m+n)+\frac{m^3-m}{12}\cdot\frac{6}{5}
\cdot\delta_{m+n,0},
\end{equation}
\begin{equation}\label{eq:L1J}
[L^1(m),\,J(n)]=(2m-n)J(m+n),
\end{equation}
\begin{equation}\label{eq:JJ}
\begin{split}
[J(m),\,J(n)]
&=(m-n)\Big(22(m+n+2)(m+n+3) + 35(m+2)(n+2)\Big)L^1(m+n)\\
&\qquad -120(m-n)\Big( \sum_{k \le -2} L^1(k)L^1(m+n-k) +
\sum_{k \ge -1}L^1(m+n-k)L^1(k) \Big)\\
&\qquad -\frac{7}{10}m(m^2-1)(m^2-4)\delta_{m+n,0}.
\end{split}
\end{equation}

The vertex operator algebra $M_t^0$ is known as a $3$-State Potts
model. It is a $W_3$ algebra of central charge $4/5$ with the
Virasoro element $\tom^2$ and is generated by $\tom^2$ and $K$. Both
of $\tom^2$ and $K$ are fixed by $\tau$, so that $\tau$ is the
identity on $M_t^0$. The rationality of $M_t^0$ was established in
\cite{KMY} and the $C_2$-cofiniteness of $M_t^0$ follows from
\cite{Buhl}. By a direct calculation, we can verify that
\begin{equation}\label{eq:KnK}
\begin{split}
K_5 K &= 104\cdot\1,\\
K_4 K &= 0,\\
K_3 K &= 780\tom^2,\\
K_2 K &= 390(\tom^2)_0 \tom^2,\\
K_1 K &= -27(\tom^2)_0 (\tom^2)_0 \tom^2 +
480(\tom^2)_{-1}\tom^2,\\
K_0 K &= -46(\tom^2)_0(\tom^2)_0(\tom^2)_0\tom^2 +
240(\tom^2)_0(\tom^2)_{-1}\tom^2.
\end{split}
\end{equation}

Let $L^2(n)=(\tom^2)_{n+1}$ and $K(n)=K_{n+2}$ for $n \in \Z$. Then
\begin{equation}\label{eq:L2L2}
[L^2(m),\,L^2(n)]=(m-n)L^2(m+n)+\frac{m^3-m}{12}\cdot\frac{4}{5}
\cdot\delta_{m+n,0},
\end{equation}
\begin{equation}\label{eq:L2K}
[L^2(m),\,K(n)]=(2m-n)K(m+n),
\end{equation}
\begin{equation}\label{eq:KK}
\begin{split}
[K(m),K(n)] &= -(m-n)\Big( 46(m+n+2)(m+n+3) + 65(m+2)(n+2)\Big)
L^2(m+n)\\
& \qquad +240(m-n)\Big(\sum_{k \le -2} L^2(k)L^2(m+n-k) + \sum_{k
\ge -1} L^2(m+n-k)L^2(k)\Big)\\
& \qquad +\frac{13}{15} m(m^2-1)(m^2-4)\delta_{m+n,0}.
\end{split}
\end{equation}

\begin{rmk}
Let $L_n=L^1(n)$, $W_n = \sqrt{-1/210}J(n)$, and $c = 6/5$. Then the
above commutation relations coincide with (2.1) and (2.2) of
\cite{BMP}. The same commutation relations also hold if we set
$L_n=L^2(n)$, $W_n=K(n)/\sqrt{390}$, and $c = 4/5$.
\end{rmk}

Let us review the $20$ inequivalent simple $M(0)$-modules studied
in \cite{DLTYY}. Among those simple $M(0)$-modules, eight of
them appear in simple $M_k^0$-modules, namely,
\begin{equation*}
M(\vep) = \{ u \in M_k^0 \,|\, \tau u = \xi^\vep u \}, \quad
W(\vep) = \{ u \in W_k^0 \,|\, \tau u = \xi^\vep u \}
\end{equation*}
for $\vep = 0,1,2$, $M_k^c$ and $W_k^c$. The remaining $12$ simple
$M(0)$-modules appear in simple $\tau$-twisted or $\tau^2$-twisted
$V_L$-modules. Let
\begin{align*}
M_T(\tau)(\vep) & = \{ u \in V_L^{T_{\chi_0}}(\tau) \,|\,
(\tom^2)_1
u = 0,\, \tau u = \xi^\vep u \},\\
W_T(\tau)(\vep) & = \{ u \in V_L^{T_{\chi_0}}(\tau) \,|\, (\tom^2)_1
u = \frac{2}{5} u,\, \tau u = \xi^\vep u \}.
\end{align*}
Then $M_T(\tau)(\vep)$, $W_T(\tau)(\vep)$, $\vep = 0,1,2$ are
inequivalent simple $M(0)$-modules. Similarly,
\begin{align*}
M_T(\tau^2)(\vep) & = \{ u \in V_L^{T_{\chi'_0}}(\tau^2) \,|\,
(\tom^2)_1 u = 0,\, \tau^2 u = \xi^\vep u \},\\
W_T(\tau^2)(\vep) & = \{ u \in V_L^{T_{\chi'_0}}(\tau^2) \,|\,
(\tom^2)_1 u = \frac{2}{5} u,\, \tau^2 u = \xi^\vep u \}
\end{align*}
for $\vep = 0,1,2$ are inequivalent simple $M(0)$-modules. In
\cite{DLTYY}, it was shown  that $M(\vep)$, $W(\vep)$, $M_k^c$,
$W_k^c$, $M_T(\tau)(\vep)$, $W_T(\tau)(\vep)$,
$M_T(\tau^2)(\vep)$, and $W_T(\tau^2)(\vep)$, $\vep = 0,1,2$ form
a complete set of representatives of isomorphism classes of simple
$M(0)$-modules.

Let us describe the structure of the fixed point subalgebra
$V_L^\tau$. By the definition of $M(0)$ and $M_t^0$, we see that
$V_L^\tau \supset M(0) \otimes M_t^0$. Since both of $M(0)$ and
$M_t^0$ are rational, $M(0) \ot M_t^0$ is also rational. Thus
$V_L(\vep) = \{ u \in V_L\,|\, \tau u = \xi^\vep u \}$, $\vep =
0,1,2$ can be decomposed into a direct sum of simple modules for
$M(0) \otimes M_t^0$. Any simple module for $M(0) \otimes M_t^0$ is
of the form $A \otimes B$, where $A$ and $B$ are simple modules for
$M(0)$ and $M_t^0$, respectively. By \eqref{eq:VL-dec}, it follows
that $B \cong M_t^0$ or $W_t^0$. Moreover, $V_L(\vep)$ contains the
simple $M(0)$-modules $M(\vep)$ and $W(\vep)$. The eigenvalues of
$(\tom^1)_1$ in $M(\vep)$ (resp. $W(\vep)$) are integers (resp. of
the form $3/5 + n$, $n \in \Z$), while the eigenvalues of
$(\tom^2)_1$ in $M_t^0$ (resp. $W_t^0$) are integers (resp. of the
form $2/5 + n$, $n \in \Z$). Since the eigenvalues of $\om_1 =
(\tom^1)_1 + (\tom^2)_1$ in $V_L$ are integers, we conclude that
\begin{equation}\label{eq:VLvep-dec}
V_L(\vep) \cong ( M(\vep) \otimes M_t^0) \oplus (W(\vep) \otimes
W_t^0)
\end{equation}
as $M(0) \otimes M_t^0$-modules, $\vep = 0,1,2$. In particular,
\begin{equation}\label{eq:VLtau-dec}
V_L^\tau \cong ( M(0) \otimes M_t^0) \oplus (W(0) \otimes W_t^0).
\end{equation}

>From now on we set $M^0 = M(0) \otimes M_t^0$ and $W^0 = W(0)
\otimes W_t^0$. Thus $V_L^\tau = V_L(0) \cong M^0 \op W^0$. Let
\begin{equation*}
P = y(\al_1) + y(\al_2) + y(\al_0).
\end{equation*}
Then  we can verify that $(\tom^1)_n P = (\tom^2)_n P = 0$ for $n
\ge 2$, $(\tom^1)_1 P = (8/5)P$, and $(\tom^2)_1 P = (2/5)P$.
Moreover, $J_n P = K_n P = 0$ for $n \ge 2$. Thus $W^0$ is a simple
$M^0$-module with $P$ a highest weight vector of weight $(8/5,2/5)$.
The vertex operator algebra $V_L^\tau$ is generated by $\tom^1$,
$\tom^2$, $J$, $K$ and $P$.

\begin{thm}\label{thm:C2}
$V_L^\tau$ is a simple $C_2$-cofinite vertex operator algebra.
\end{thm}

\begin{proof}
We know that $M(0)$ and $M_t^0$ are $C_2$-cofinite. Thus $M^0$ is
also $C_2$-cofinite. Since $W^0$ is generated by $P$ as an
$M^0$-module, it follows from \cite{Buhl} that $V_L^\tau$ is
$C_2$-cofinite. By \cite[Theorem 4.4]{DM}, $V_L^\tau$ is simple.
\end{proof}

Following the outline of the argument in \cite{DLTYY, KLY2}, we
discuss the structure of the simple $\tau$-twisted $V_L$-modules
$V_L^{T_{\chi_j}}(\tau)$, $j=0,1,2$ as $\tau$-twisted $M_k^0 \otimes
M_t^0$-modules. Furthermore, we correct an error in \cite{DLTYY,
KLY2} concerning a decomposition of $V_L^{T_{\chi_j}}(\tau)$ for
$j=1,2$. We first consider $V_L^{T_{\chi_0}}(\tau)$. Let $0 \ne v
\in T_{\chi_0}$ and $1$ be the identity of $S[\tau]$. Then $1 \ot v
\in S[\tau] \ot T_{\chi_0} = V_L^{T_{\chi_0}}(\tau)$. Since $M_t^0
\subset V_L^\tau$, we can decompose $V_L^{T_{\chi_0}}(\tau)$ into a
direct sum of simple $M_t^0$-modules. By a direct calculation, we
can verify that
\begin{equation*}
(\tom^2)_1 (1 \ot v) =0,\qquad (\tom^2)_1 (h_2(-\frac{1}{3}) \ot v)
= \frac{2}{5} h_2(-\frac{1}{3}) \ot v.
\end{equation*}
Thus we see that $M_t^0$ and $W_t^0$ appear as direct summands.
Since $V_L^{T_{\chi_0}}(\tau)$ is simple as a $\tau$-twisted
$V_L$-module, \eqref{eq:VL-dec} and the fusion rule $W_t^0 \times
W_t^0 = M_t^0 + W_t^0$ (cf. \eqref{eq:Ternary-fusion}) imply that
any simple $M_t^0$-submodule of $V_L^{T_{\chi_0}}(\tau)$ is
isomorphic to $M_t^0$ or $W_t^0$. Hence
\begin{equation}\label{eq:VLT-dec}
V_L^{T_{\chi_0}}(\tau) \cong
(M_T^0(\tau) \ot M_t^0) \oplus
(W_T^0(\tau) \ot W_t^0)
\end{equation}
as $\tau$-twisted $M_k^0 \ot M_t^0$-modules, where
\begin{align*}
M_T^0(\tau) &= \{ u \in V_L^{T_{\chi_0}}(\tau) \,|\,
(\tom^2)_1 u = 0\},\\
W_T^0(\tau) &= \{ u \in V_L^{T_{\chi_0}}(\tau) \,|\, (\tom^2)_1 u =
\frac{2}{5}u \}.
\end{align*}

The $\tau$-twisted $M_k^0$-modules $M_T^0(\tau)$ and $W_T^0(\tau)$
are simple . Indeed, if $N$ is a $\tau$-twisted $M_k^0$-submodule of
$M_T^0(\tau)$, then $N \ot M_t^0$ is a $\tau$-twisted $M_k^0 \ot
M_t^0$-submodule of $M_T^0(\tau) \ot M_t^0$. By \eqref{eq:Formula1},
$V_L \cdot (N \ot M_t^0) = \spn \{ a_n (N \ot M_t^0)\,|\, a \in V_L,
n \in \Q\}$ is a $\tau$-twisted $V_L$-submodule of
$V_L^{T_{\chi_0}}(\tau)$. The fusion rule $W_t^0 \times M_t^0 =
W_t^0$ and \eqref{eq:VL-dec} imply that $V_L \cdot (N \ot M_t^0)$ is
contained in $(N \ot M_t^0) \oplus (W_T^0(\tau) \ot W_t^0)$. Since
$V_L^{T_{\chi_0}}(\tau)$ is a simple $\tau$-twisted $V_L$-module, we
conclude that $M_T^0(\tau)$ is a simple $\tau$-twisted
$M_k^0$-module.

Because of the fusion rule $W_t^0 \times W_t^0 = M_t^0 + W_t^0$,
we can not apply a similar argument to $W_T^0(\tau)$. Note that
there are at most two inequivalent simple $\tau$-twisted
$M_k^0$-modules by \cite[Lemma 4.1]{DLTYY} and \cite[Theorem
10.2]{DLM2}. Note also that a weight in $M_T^0(\tau)$ or in
$W_T^0(\tau)$ means an eigenvalue of $(\tom^1)_1$. First several
terms of the characters of $M_T^0(\tau)$ and $W_T^0(\tau)$ can be
calculated easily from \eqref{eq:VLT-dec} (cf. \cite{DLTYY}).
\begin{align*}
\ch M_T^0(\tau) &= q^{1/9} + q^{1/9+2/3} + q^{1/9+1} + q^{1/9+4/3}
+ \cdots,\\
\ch W_T^0(\tau) &= q^{2/45} + q^{2/45+1/3} + q^{2/45+2/3} +
q^{2/45+1} + \cdots.
\end{align*}

Suppose $W_T^0(\tau)$ is not a simple $\tau$-twisted $M_k^0$-module.
Let $N$ be the $\tau$-twisted $M_k^0$-submodule of $W_T^0(\tau)$
generated by the top level of $W_T^0(\tau)$. Then the top level of
$N$ is a one dimensional space of weight $2/45$. If $N$ is not a
simple $\tau$-twisted $M_k^0$-module, then the sum $U$ of all proper
$\tau$-twisted $M_k^0$-submodules of $N$ is a unique maximal
$\tau$-twisted $M_k^0$-submodule of $N$. The quotient $N/U$ is a
simple $\tau$-twisted $M_k^0$-module whose top level is of weight
$2/45$. Denote the top level of $U$ by $U_\lambda$, where the weight
$\lambda$ is $2/45 + n/3$ for some $1 \le n \in \Z$. Consider the
$\tau$-twisted Zhu algebra $A_\tau(M_k^0)$ of $M_k^0$. Since
$U_\lambda$ is a finite dimensional $A_\tau(M_k^0)$-module, we can
choose a simple $A_\tau(M_k^0)$-submodule $S$ of $U_\lambda$. By
\cite[Proposition 5.4 and Theorem 7.2]{DLM1}, there is a simple
$\frac{1}{3}\N$-graded weak $\tau$-twisted $M_k^0$-module $R$ with
top level $R_\lambda$ being isomorphic to $S$ as an
$A_\tau(M_k^0)$-module. It follows from \cite[Corollary
3.8]{Yamauchi} that $R$ is in fact a simple $\tau$-twisted
$M_k^0$-module. Here we note that $M_k^0$ is $C_2$-cofinite and of
CFT type by its structure \eqref{eq:MW-k}. Since the top levels of
$M_T^0(\tau)$, $N/U$, and $R$ have different weight, they are
inequivalent simple $\tau$-twisted $M_k^0$-modules. If $N$ is a
simple $\tau$-twisted $M_k^0$-module, then it is not equal to
$W_T^0(\tau)$ by our assumption. The quotient $W_T^0(\tau)/N$ is a
$\tau$-twisted $M_k^0$-module and the weight of its top level, say
$\mu$ is $2/45+m/3$ for some $1 \le m \in \Z$. By a similar argument
as above, we see that there is a simple $\tau$-twisted
$M_k^0$-module whose top level is of weight $\mu$. Hence we have
three inequivalent simple $\tau$-twisted $M_k^0$-modules in both
cases. This contradicts the fact that there are at most two
inequivalent simple $\tau$-twisted $M_k^0$-modules. Thus
$W_T^0(\tau)$ is a simple $\tau$-twisted $M_k^0$-module.

Next, let $0 \ne v \in T_{\chi_j}$, $j=1,2$. From the
definition of $V_L^{T_{\chi_j}}(\tau)$ in \cite{DLTYY, KLY2},
we can calculate that
\begin{equation*}
(\tom^2)_1 (1 \ot v) = \frac{1}{15}(1 \ot v),
\qquad (\tom^2)_1 u^j = \frac{2}{3}u^j,
\end{equation*}
where $u^j = h_1(-\frac{2}{3}) \ot v - (-1)^j \sqrt{-3}
h_2(-\frac{1}{3})^2 \ot v$. Thus $M_t^1$ or $M_t^2$ and $W_t^1$ or
$W_t^2$ appear as $M_t^0$-submodules of $V_L^{T_{\chi_j}}(\tau)$. In
order to distinguish $M_t^1$ and $M_t^2$ (resp. $W_t^1$ and
$W_t^2$), we need to know the action of $K_2$ on these vectors (cf.
\cite{KMY}). By a direct calculation, we can verify that
\begin{equation*}
K_2 (1 \ot v) = -(-1)^j \frac{2}{9} (1 \ot v), \qquad
K_2 u^j = (-1)^j \frac{52}{9} u^j.
\end{equation*}
Hence $M_t^{3-j}$ and $W_t^{3-j}$ appear in
$V_L^{T_{\chi_j}}(\tau)$ for $j=1,2$. Let
\begin{align*}
M_T^j(\tau) &= \{ u \in V_L^{T_{\chi_j}}(\tau) \,|\,
(\tom^2)_1 u = \frac{2}{3}u \},\\
W_T^j(\tau) &= \{ u \in V_L^{T_{\chi_j}}(\tau) \,|\, (\tom^2)_1 u =
\frac{1}{15}u \}, \quad j=1,2.
\end{align*}
Then,
$V_L^{T_{\chi_j}}(\tau) \cong
(M_T^j(\tau) \ot M_t^{3-j}) \oplus
(W_T^j(\tau) \ot W_t^{3-j})$
as $\tau$-twisted $M_k^0 \ot M_t^0$-modules for $j=1,2$.
Moreover, $M_T^j(\tau)$ and $W_T^j(\tau)$, $j=1,2$ are
simple $\tau$-twisted $M_k^0$-modules.

Recall that there are at most two inequivalent simple
$\tau$-twisted $M_k^0$-modules. Looking at the smallest weight of
$M_T^j(\tau)$ and $W_T^j(\tau)$, we have that $M_T^j(\tau)$,
$j=0,1,2$ are equivalent and $W_T^j(\tau)$, $j=0,1,2$ are
equivalent, but $M_T^0(\tau)$ and $W_T^0(\tau)$ are not
equivalent. For simplicity, set $M_T(\tau) = M_T^0(\tau)$ and
$W_T(\tau) = W_T^0(\tau)$. Then
\begin{equation}\label{eq:Twisted-1}
\begin{split}
V_L^{T_{\chi_0}}(\tau) & \cong (M_T(\tau) \ot M_t^0)
\oplus (W_T(\tau) \ot W_t^0),\\
V_L^{T_{\chi_j}}(\tau) & \cong (M_T(\tau) \ot M_t^{3-j})
\oplus (W_T(\tau) \ot W_t^{3-j}), \quad j=1,2
\end{split}
\end{equation}
as $\tau$-twisted $M_k^0 \ot M_t^0$-modules.

The structure of the simple $\tau^2$-twisted $V_L$-module
$V_L^{T_{\chi'_j}}(\tau^2)$, $j=0,1,2$ as a $\tau^2$-twisted
$M_k^0 \ot M_t^0$-module is similar to that of the case for
$V_L^{T_{\chi_j}}(\tau)$. Let $0 \ne v \in T_{\chi'_0}$ and
$1$ be the identity of $S[\tau^2]$. Then
\begin{equation*}
(\tom^2)_1 (1 \ot v) =0,\qquad
(\tom^2)_1 (h'_2(-\frac{1}{3}) \ot v) = \frac{2}{5}
h'_2(-\frac{1}{3}) \ot v
\end{equation*}
and so
\begin{equation*}
V_L^{T_{\chi'_0}}(\tau^2) \cong
(M_T^0(\tau^2) \ot M_t^0) \oplus
(W_T^0(\tau^2) \ot W_t^0)
\end{equation*}
as $\tau^2$-twisted $M_k^0 \ot M_t^0$-modules, where
\begin{align*}
M_T^0(\tau^2) &= \{ u \in V_L^{T_{\chi'_0}}(\tau^2) \,|\,
(\tom^2)_1 u = 0\},\\
W_T^0(\tau^2) &= \{ u \in V_L^{T_{\chi'_0}}(\tau^2) \,|\,
(\tom^2)_1 u = \frac{2}{5}u \}.
\end{align*}

By a similar argument as in the $\tau$-twisted case, we can show
that $M_T^0(\tau^2)$ and $W_T^0(\tau^2)$ are inequivalent simple
$\tau^2$-twisted $M_k^0$-modules.

Let $0 \ne v \in T_{\chi'_j}$, $j=1,2$. Then
\begin{equation*}
(\tom^2)_1 (1 \ot v) = \frac{1}{15}(1 \ot v),\qquad
(\tom^2)_1 v^j = \frac{2}{3}v^j,
\end{equation*}
where $v^j = h'_1(-\frac{2}{3}) \ot v - (-1)^j \sqrt{-3}
h'_2(-\frac{1}{3})^2 \ot v$. Furthermore,
\begin{equation*}
K_2 (1 \ot v) = (-1)^j \frac{2}{9}(1 \ot v), \qquad
K_2 v^j = - (-1)^j \frac{52}{9}v^j.
\end{equation*}
Hence $V_L^{T_{\chi'_j}}(\tau^2) \cong (M_T^j(\tau^2) \ot M_t^j)
\oplus (W_T^j(\tau^2) \ot W_t^j)$ as $\tau^2$-twisted $M_k^0 \ot
M_t^0$-modules for $j=1,2$, where
\begin{align*}
M_T^j(\tau^2) &= \{ u \in V_L^{T_{\chi'_j}}(\tau^2) \,|\,
(\tom^2)_1 u = \frac{2}{3}u \},\\
W_T^j(\tau^2) &= \{ u \in V_L^{T_{\chi'_j}}(\tau^2) \,|\, (\tom^2)_1
u = \frac{1}{15}u \}, \quad j=1,2.
\end{align*}

As in the $\tau$-twisted case, $M_T^j(\tau^2)$, $j=0,1,2$ are
equivalent and $W_T^j(\tau^2)$, $j=0,1,2$ are equivalent. Set
$M_T(\tau^2) = M_T^0(\tau^2)$ and $W_T(\tau^2) = W_T^0(\tau^2)$.
Then
\begin{equation}\label{eq:Twisted-2}
V_L^{T_{\chi'_j}}(\tau^2) \cong
(M_T(\tau^2) \ot M_t^j) \oplus
(W_T(\tau^2) \ot W_t^j), \quad j=0,1,2
\end{equation}
as $\tau^2$-twisted $M_k^0 \ot M_t^0$-modules.

\begin{rmk}\label{rmk:correction}
The weight $3$ vector $K$ was denoted by different symbols in
previous papers, namely, $v_t$, $v^3$, and $q$ were used in
\cite{DLTYY}, \cite{KLY2}, and \cite{KMY}, respectively. They are
related as follows: $K = - 2\sqrt{2} v_t = - 2\sqrt{2} v^3 =
2\sqrt{2} q$. Thus, in the proof of \cite[Proposition 6.8]{KLY2}
$(v^3)_2$ should act on the top level of $V_L^{T_{\chi_j}}(\tau)$ as
a scalar multiple of $(-1)^j/9\sqrt{2}$ for $j=1,2$. Moreover,
(6.46) of \cite{KLY2} and the equation for $V_L^{T_{\chi_j}}(\tau)$
on page 265 of \cite{DLTYY} should be replaced with Equation
\eqref{eq:Twisted-1}. This correction does not affect the results in
\cite{DLTYY}. However, certain changes must be necessary in
\cite{KLY2} along the correction.
\end{rmk}

Note that
\begin{align*}
M_T(\tau^i)(\vep) &= \{ u \in M_T(\tau^i) \,|\,
\tau^i u = \xi^\vep u\},\\
W_T(\tau^i)(\vep) &= \{ u \in W_T(\tau^i) \,|\,
\tau^i u = \xi^\vep u \}
\end{align*}
for $i=1,2$, $\vep=0,1,2$. Another notation was used in
\cite{KLY2}, namely,
\begin{equation*}
M_T(\tau^i)^\vep = \bigoplus_{n \in 1/9 + \vep/3 + \Z}
(M_T(\tau^i))_n, \qquad
W_T(\tau^i)^\vep = \bigoplus_{n \in 2/45 + \vep/3 + \Z}
(W_T(\tau^i))_n,
\end{equation*}
where $U_n$ denotes the eigenspace of $U$ with
eigenvalue $n$ for $(\tom^1)_1$. They are related as follows.
\begin{equation}\label{eq:Rel-simple}
M_T(\tau^i)^\vep = M_T(\tau^i)(2\vep), \qquad
W_T(\tau^i)^\vep = W_T(\tau^i)(2\vep-1).
\end{equation}

Likewise,
\begin{equation*}
(V_L^{T_{\chi_j}}(\tau))^\vep =
\bigoplus_{n \in 1/9 + \vep/3 + \Z} (V_L^{T_{\chi_j}}(\tau))_n,
\qquad
(V_L^{T_{\chi'_j}}(\tau^2))^\vep =
\bigoplus_{n \in 1/9 + \vep/3 + \Z}
(V_L^{T_{\chi'_j}}(\tau^2))_n
\end{equation*}
of \cite{KLY2} are denoted by
\begin{equation}\label{eq:Rel-twisted}
(V_L^{T_{\chi_j}}(\tau))^\vep = V_L^{T_{\chi_j}}(\tau)(2\vep),
\qquad
(V_L^{T_{\chi'_j}}(\tau^2))^\vep =
V_L^{T_{\chi'_j}}(\tau^2)(2\vep)
\end{equation}
in our notation for $j=0,1,2$ and $\vep=0,1,2$, where $U_n$ is the
eigenspace of $U$ with eigenvalue $n$ for $\om_1$.

By \eqref{eq:MW-t}, the minimal eigenvalues of $(\tom^2)_1$ on
$M_t^0$ and $W_t^0$ are $0$ and $2/5$, respectively, while those on
$M_t^j$ and $W_t^j$, $j=1,2$, are $2/3$ and $1/15$, respectively.
Hence it follows from \eqref{eq:Twisted-1} that
\begin{equation}\label{eq:Twisted-1-vep}
\begin{split}
(V_L^{T_{\chi_0}}(\tau))^\vep & \cong (M_T(\tau)^\vep \ot M_t^0)
\oplus (W_T(\tau)^{\vep - 1} \ot W_t^0),\\
(V_L^{T_{\chi_j}}(\tau))^\vep & \cong (M_T(\tau)^{\vep +1} \ot
M_t^{3-j})
\oplus (W_T(\tau)^\vep \ot W_t^{3-j}), \quad j=1,2
\end{split}
\end{equation}
as $M^0$-modules for $\vep = 0,1,2$, where $M^0 = M(0) \ot M_t^0$.
Similarly,
\begin{equation}\label{eq:Twisted-2-vep}
\begin{split}
(V_L^{T_{\chi'_0}}(\tau^2))^\vep & \cong (M_T(\tau^2)^\vep
\ot M_t^0)
\oplus (W_T(\tau^2)^{\vep - 1} \ot W_t^0),\\
(V_L^{T_{\chi'_j}}(\tau^2))^\vep & \cong (M_T(\tau^2)^{\vep +1}
\ot M_t^j)
\oplus (W_T(\tau^2)^\vep \ot W_t^j), \quad j=1,2
\end{split}
\end{equation}
as $M^0$-modules for $\vep = 0,1,2$ (cf. \cite[(7.17)]{KLY2}).

The following fusion rules of simple $M(0)$-modules will be
necessary for the study of simple $V_L^\tau$-modules.
\begin{equation}\label{eq:W0-fusion}
\begin{split}
W(0) \times M_k^c &= W_k^c,\\
W(0) \times W_k^c &= M_k^c + W_k^c,\\
W(0) \times M(\vep) &= W(\vep),\\
W(0) \times W(\vep) &= M(\vep) + W(\vep),\\
W(0) \times M_T(\tau^i)(\vep) &= W_T(\tau^i)(\vep),\\
W(0) \times W_T(\tau^i)(\vep) &= M_T(\tau^i)(\vep) +
W_T(\tau^i)(\vep)
\end{split}
\end{equation}
for $i=1,2$ and $\vep = 0,1,2$. In fact, the first four fusion
rules, that is, fusion rules among simple $M(0)$-modules which
appear in untwisted simple $V_L$-modules can be found in
\cite{Tanabe}. The last two fusion rules involve simple
$M(0)$-modules which appear in $\tau^i$-twisted simple
$V_L$-modulesis. Their proofs are given in Appendix
\ref{app:fusion-rules}.

Fusion rules possess certain symmetries. Let $M^i$, $i=1,2,3$ be
modules for a vertex operator algebra $V$. Then by
\cite[Propositions 5.4.7 and 5.5.2]{FHL}
\begin{equation*}
\dim I_V \binom{M^3}{M^1 \  M^2} = \dim I_V \binom{M^3}{M^2 \ M^1}
= \dim I_V \binom{(M^2)'}{M^1 \  (M^3)'},
\end{equation*}
where $(M^i)'$ is the contragredient module of $M^i$. Recall that
the contragredient module $(U',Y_{U'})$ of a $V$-module $(U,Y_U)$ is
defined as follows. As a vector space $U' = \op_n (U_n)^*$ is the
restricted dual of $U$ and $Y_{U'}(\,\cdot\,,z)$ is determined by
\begin{equation*}
\la Y_{U'}(a,z)v, u \ra = \la v, Y_U(e^{zL(1)}(-z^{-2})^{L(0)}a,
z^{-1})u \ra
\end{equation*}
for $a \in V$, $u \in U$, and $v \in U'$.

In our case $M(0)$ is generated by the Virasoro element $\tom^1$ and
the weight $3$ vector $J$. Moreover, $\la L^1(0)v, u \ra = \la v,
L^1(0)u\ra$ and $\la J(0)v, u \ra = - \la v, J(0)u \ra$. Since the
$20$ simple $M(0)$-modules are distinguished by the action of
$L^1(0)$ and $J(0)$ on their top levels, we know from Tables 1, 3,
and 4 of \cite{DLTYY} that the contragredient modules of the simple
$M(0)$-modules are as follows.
\begin{gather*}
M(\vep)' \cong M(2\vep), \quad W(\vep)' \cong W(2\vep), \quad
\vep = 0,1,2,\\
(M_k^c)' \cong M_k^c, \quad (W_k^c)' \cong W_k^c,\\
M_T(\tau)(\vep)' \cong M_T(\tau^2)(\vep), \quad W_T(\tau)(\vep)'
\cong W_T(\tau^2)(\vep), \quad \vep = 0,1,2
\end{gather*}
(see also \cite[Lemma 3.7]{DLM1} and \cite[Section 4.2]{Tanabe}).

\section{Structure of simple modules}\label{sec:Structure}

Recall that $V_L^\tau = V_L(0) = M^0 \op W^0$ with $M^0 = M(0) \ot
M_t^0$ and $W^0 = W(0) \ot W_t^0$. In this section we study the
structure of the $30$ known simple $V_L^\tau$-modules listed in
Lemma \ref{lem:30SIMPLES}. We discuss decompositions of these simple
modules as modules for $M^0$. Those decompositions have been
obtained in \cite{KLY2}. We review them briefly. (Some corrections
are necessary in \cite{KLY2}, see Remark \ref{rmk:correction} for
details.)

A vector in a $V_L^\tau$-module is said to be of weight $h$ if it is
an eigenvector for $L(0)=\om_1$ with eigenvalue $h$. We calculate
the action of $(\tom^1)_1$, $(\tom^2)_1$, $J_2$, $K_2$, $P_1$, $(J_1
P)_2$, and $(K_1 P)_2$ on the top levels of the $30$ known simple
$V_L^\tau$-modules. Recall that the top level of a module means the
homogeneous subspace of the module of smallest weight. The
calculation is accomplished directly from the definition of
untwisted or twisted vertex operators associated with the lattice
$L$ and the automorphisms $\tau$ and $\tau^2$ (cf. \cite{DL2, FLM,
LL}). The results in this section will be used to determine the Zhu
algebra $A(V_L^\tau)$ of $V_L^\tau$ in Section
\ref{sec:Classification}.

The vectors $J_1 P$ and $K_1 P$ are of weight $3$. Their precise
form in terms of the lattice vertex operator algebra $V_L$ is as
follows.
\begin{align*}
J_1 P & = 2\beta_1(-1)^3 + 3\beta_1(-1)^2 \beta_2(-1)
- 3\beta_1(-1) \beta_2(-1)^2 - 2\beta_2(-1)^3 \\
& \qquad -4 \Big( (\beta_2 - \beta_0)(-1)x(\alpha_1) + (\beta_0 -
\beta_1)(-1)x(\alpha_2)
+ (\beta_1 - \beta_2)(-1)x(\alpha_0) \Big)\\
& = \frac{13}{9} \Big( 2\beta_1(-1)^3 + 3\beta_1(-1)^2 \beta_2(-1)
- 3\beta_1(-1) \beta_2(-1)^2 - 2\beta_2(-1)^3 \Big) - 4K,
\end{align*}
\begin{align*}
K_1 P & = 3 \Big( \beta_1(-2)\beta_2(-1) -
\beta_2(-2)\beta_1(-1) \Big)\\
& \qquad - \Big( (\beta_2 - \beta_0)(-1)y(\alpha_1) + (\beta_0 -
\beta_1)(-1)y(\alpha_2)
+ (\beta_1 - \beta_2)(-1)y(\alpha_0) \Big)\\
& = \frac{7}{2} \Big( \beta_1(-2)\beta_2(-1) -
\beta_2(-2)\beta_1(-1) \Big) + J.
\end{align*}

\subsection{The simple module $V_L(0)$}
$V_L(0) = M^0 \oplus W^0$ as $M^0$-modules. The top level of
$V_L(0)$ is $\C\1$. By a property of the vacuum vector, all of
$(\tom^1)_1$, $(\tom^2)_1$, $J_2$, $K_2$, $P_1$, $(J_1 P)_2$, and
$(K_1 P)_2$ act as $0$ on $\C\1$.

\subsection{The simple module $V_L(\vep), \vep=1,2$}
By  \eqref{eq:VLvep-dec}, $V_L(\vep) \cong ( M(\vep) \otimes M_t^0 )
\oplus ( W(\vep)\ot W_t^0)$ as $M^0$-modules for $\vep=1,2$. The top
level of $V_L(\vep)$ is $\C\bv^{2,\vep}$, where
$\bv^{2,\vep}=\al_1(-1) - \xi^\vep \al_2(-1) \in W(\vep)\ot W_t^0$.
The following hold.
\begin{gather*}
(\tom^1)_1 \bv^{2,\vep} = \frac{3}{5}\bv^{2,\vep}, \quad
(\tom^2)_1 \bv^{2,\vep} = \frac{2}{5}\bv^{2,\vep}, \quad
J_2 \bv^{2,\vep} = -(-1)^{\vep} 2\sqrt{-3}\bv^{2,\vep},\\
K_2 \bv^{2,\vep} = 0,\quad P_1 \bv^{2,\vep} = 0,\quad (J_1 P)_2
\bv^{2,\vep} = 0, \quad (K_1 P)_2 \bv^{2,\vep} = (-1)^\vep 12
\sqrt{-3} \bv^{2,\vep}.
\end{gather*}

\subsection{The simple module $V_{L^{(0,j)}}(0), j=1,2$}
For $j=1,2$, \eqref{eq:VLij-dec} implies that $V_{L^{(0,j)}}$ is a
direct sum of simple $M^0$-modules of the form $A \otimes B$, where
$A$ is a simple $M(0)$-module and $B$ is a simple $M_t^0$-module
isomorphic to $M_t^j$ or $W_t^j$. For convenience, set $U^j(\vep) =
V_{L^{(0,j)}}(\vep)$, $j=1,2$, $\vep=0,1,2$. Let
\begin{equation*}
\bv^{3,j}  =
e^{(-1)^j(\be_1 - \be_2)/3} + e^{(-1)^j(\be_2 - \be_0)/3} +
e^{(-1)^j(\be_0 - \be_1)/3}.
\end{equation*}
Then $\bv^{3,j} \in U^j(0)$. Moreover, $(\om^1)_1 \bv^{3,j}
= (\om^2)_1 \bv^{3,j} = 0$ and $(\tom^2)_1 \bv^{3,j} =
(2/3) \bv^{3,j}$. Hence $\bv^{3,j} \in M_t^j$ and $U^j(0)$
contains an $M_t^0$-submodule isomorphic to $M_t^j$.
By the fusion rule $M_t^j \times W_t^0 = W_t^j$
of $M_t^0$-modules and \cite[Proposition 11.9]{DL2},
$U^j(0)$ contains an $M_t^0$-submodule isomorphic to
$W_t^j$ also. Thus $U^j(0)$ contains simple
$M^0$-submodules of the form $A \otimes M_t^j$
and $A' \otimes W_t^j$ for some simple $M(0)$-modules $A$
and $A'$.

The minimal weight of $V_{L^{(0,j)}}$ is $2/3$. Its weight subspace
is of dimension $3$ and spanned by $e^{(-1)^j(\be_1 - \be_2)/3}$,
$e^{(-1)^j(\be_2 - \be_0)/3}$, and $e^{(-1)^j(\be_0 - \be_1)/3}$.
Thus the weight $2/3$ subspace of $U^j(0)$ is $\C \bv^{3,j}$. Since
$(\tom^1)_1 \bv^{3,j} = 0$ and since only $M(0)$ is the simple
$M(0)$-module whose minimal weight ($=$ eigenvalue of $(\tom^1)_1$)
is $0$ by \cite{DLTYY}, we conclude that $U^j(0)$ contains a simple
$M^0$-submodule isomorphic to $M(0) \otimes M_t^j$.

The minimal eigenvalue of $(\tom^2)_1$ in $W_t^j$ is $1/15$. Thus
the eigenvalues of $(\tom^1)_1$ on $A'$ must be of the form $3/5 +
n$, $n \in \Z$. By \cite{DLTYY}, only $W(0)$, $W(1)$, $W(2)$ are the
simple $M(0)$-modules whose weights are of this form. The minimal
weight of these simple modules are $8/5$, $3/5$ and $3/5$,
respectively. Since the weight $2/3$ subspace of $U^j(0)$ is one
dimensional, we see that $U^j(0)$ contains a simple $M^0$-submodule
isomorphic to $W(0) \otimes W_t^j$.

>From the fusion rules for $M_t^0$-modules, we obtain the fusion
rules
\begin{align*}
(M(\vep) \otimes M_t^0) \times (M(0) \otimes M_t^j)
&= M(\vep) \otimes M_t^j,\\
(W(\vep) \otimes W_t^0) \times (M(0) \otimes M_t^j)
&= W(\vep) \otimes W_t^j
\end{align*}
for $M^0$-modules. Hence $U^j(\vep) \cong (M(\vep) \otimes M_t^j)
\oplus (W(\vep) \otimes W_t^j)$ for $j=1,2$ and $\vep=0,1,2$ by
\eqref{eq:VLij-dec} and \eqref{eq:VLvep-dec}. In particular,
$V_{L^{(0,j)}}(0) \cong (M(0) \otimes M_t^j) \oplus (W(0) \otimes
W_t^j)$ as $M^0$-modules, $j=1,2$. The top level of
$V_{L^{(0,j)}}(0)$ is $\C \bv^{3,j} \subset M(0) \otimes M_t^j$. The
following hold.
\begin{gather*}
(\tom^1)_1 \bv^{3,j} = 0, \quad
(\tom^2)_1 \bv^{3,j} = \frac{2}{3} \bv^{3,j}, \quad
J_2 \bv^{3,j} = 0, \quad
K_2 \bv^{3,j} = -(-1)^j \frac{52}{9}\bv^{3,j}, \\
P_1 \bv^{3,j} = 0, \quad
(J_1 P)_2 \bv^{3,j} = 0, \quad (K_1 P)_2 \bv^{3,j} = 0.
\end{gather*}

\subsection{The simple module $V_{L^{(0,j)}}(\vep), j=1,2, \vep = 1,2$}
We have shown above that $V_{L^{(0,j)}}(\vep) \cong (M(\vep) \ot
M_t^j) \op (W(\vep)\ot W_t^j)$ as $M^0$-modules, $j=1,2$,
$\vep=1,2$. The top level of $V_{L^{(0,j)}}(\vep)$ is $\C
\bv^{4,j,\vep}$, where
\begin{equation*}
\bv^{4,j,\vep} =
e^{(-1)^j(\be_1 - \be_2)/3} +
\xi^{2\vep}e^{(-1)^j(\be_2 - \be_0)/3} +
\xi^{\vep}e^{(-1)^j(\be_0 - \be_1)/3} \in W(\vep) \ot W_t^j.
\end{equation*}

The following hold.
\begin{gather*}
(\tom^1)_1 \bv^{4,j,\vep} = \frac{3}{5} \bv^{4,j,\vep}, \quad
(\tom^2)_1 \bv^{4,j,\vep} = \frac{1}{15} \bv^{4,j,\vep}, \quad
J_2 \bv^{4,j,\vep} = -(-1)^{\vep}2\sqrt{-3} \bv^{4,j,\vep},\\
K_2 \bv^{4,j,\vep} = (-1)^j\frac{2}{9}\bv^{4,j,\vep},
\quad
P_1 \bv^{4,j,\vep} = -(-1)^{j+\vep}\sqrt{-3} \bv^{4,j,\vep},\\
(J_1 P)_2 \bv^{4,j,\vep} = - (-1)^j 24 \bv^{4,j,\vep}, \quad
(K_1 P)_2 \bv^{4,j,\vep} = - (-1)^\vep 2 \sqrt{-3} \bv^{4,j,\vep}.
\end{gather*}

\subsection{The simple module $V_{L^{(c,0)}}$}
By \eqref{eq:VLij-dec}, $V_{L^{(c,0)}} \cong (M_k^c \ot M_t^0) \op
(W_k^c \ot W_t^0)$ as $M^0$-modules. The top level of
$V_{L^{(c,0)}}$ is of dimension $2$ with basis $\{ \bv^{5,1},
\bv^{5,2} \}$, where $\bv^{5,1} = e^{\be_1/2} - e^{-\be_1/2} \in
M_k^c \ot M_t^0$, $\bv^{5,2} = e^{\be_1/2} + e^{-\be_1/2} \in W_k^c
\ot W_t^0$. The following hold.
\begin{gather*}
(\tom^1)_1 \bv^{5,1} = \frac{1}{2} \bv^{5,1}, \quad
(\tom^1)_1 \bv^{5,2} = \frac{1}{10} \bv^{5,2}, \quad
(\tom^2)_1 \bv^{5,1} = 0, \quad
(\tom^2)_1 \bv^{5,2} = \frac{2}{5} \bv^{5,2},\\
J_2 \bv^{5,j} = 0, \quad
K_2 \bv^{5,j} = 0, \,
j=1,2, \quad
P_1 \bv^{5,1} = - \bv^{5,2}, \quad
P_1 \bv^{5,2} = \bv^{5,1},\\
(J_1 P)_2 \bv^{5,j} = 0, \quad
(K_1 P)_2 \bv^{5,j} = 0, \, j =1,2.
\end{gather*}

\subsection{The simple module $V_{L^{(c,j)}}$, $j=1,2$}
By \eqref{eq:VLij-dec}, $V_{L^{(c,j)}} \cong (M_k^c \ot M_t^j) \op
(W_k^c \ot W_t^j)$ as $M^0$-modules, $j=1,2$. The top level of
$V_{L^{(c,j)}}$ is $\C \bv^{6,j}$, where $\bv^{6,j} =
e^{-(-1)^j(\be_2 - \be_0)/6} \in W_k^c \ot W_t^j$. The following
hold.
\begin{gather*}
(\tom^1)_1 \bv^{6,j} = \frac{1}{10} \bv^{6,j}, \quad
(\tom^2)_1 \bv^{6,j} = \frac{1}{15} \bv^{6,j}, \quad
J_2 \bv^{6,j} = 0, \quad
K_2 \bv^{6,j} = (-1)^j\frac{2}{9} \bv^{6,j}, \\
P_1 \bv^{6,j} = 0, \quad
(J_1 P)_2 \bv^{6,j} = (-1)^j 2 \bv^{6,j}, \quad
(K_1 P)_2 \bv^{6,j} = 0.
\end{gather*}

\subsection{The simple module $V_L^{T_{\chi_0}} (\tau)(0)$}
By \eqref{eq:Twisted-1-vep}, $V_L^{T_{\chi_0}} (\tau)(0) \cong
(M_T(\tau)(0) \ot M_t^0) \op (W_T(\tau)(0) \ot W_t^0)$ as
$M^0$-modules. The top level of $V_L^{T_{\chi_0}} (\tau)(0)$ is
$\C\bv^7$, where $\bv^7 = 1 \ot v \in M_T(\tau)(0) \ot M_t^0$ and $0
\ne v \in T_{\chi_0}$. The following hold.
\begin{gather*}
(\tom^1)_1 \bv^7 = \frac{1}{9} \bv^7, \quad
(\tom^2)_1 \bv^7 = 0, \quad
J_2 \bv^7 = \frac{14}{81}\sqrt{-3} \bv^7, \quad
K_2 \bv^7 = 0, \\
P_1 \bv^7 = 0, \quad
(J_1 P)_2 \bv^7 = 0, \quad (K_1 P)_2 \bv^7 = 0.
\end{gather*}

\subsection{The simple module $V_L^{T_{\chi_0}} (\tau)(1)$}
By \eqref{eq:Twisted-1-vep}, $V_L^{T_{\chi_0}} (\tau)(1) \cong
(M_T(\tau)(1) \ot M_t^0) \op (W_T(\tau)(1) \ot W_t^0)$ as
$M^0$-modules. The top level of $V_L^{T_{\chi_0}} (\tau)(1)$ is of
dimension $2$ with basis $\{ \bv^{8,1}, \bv^{8,2} \}$, where
$\bv^{8,1} = h_2(-1/3)^2 \ot v \in M_T(\tau)(1) \ot M_t^0$,
$\bv^{8,2} = h_1(-2/3) \ot v \in W_T(\tau)(1) \ot W_t^0$ and $0 \ne
v \in T_{\chi_0}$. The following hold.
\begin{gather*}
(\tom^1)_1 \bv^{8,1} = \big(\frac{1}{9} +
\frac{2}{3}\big) \bv^{8,1}, \quad
(\tom^2)_1 \bv^{8,1} = 0, \quad
J_2 \bv^{8,1} = - \frac{238}{81}\sqrt{-3} \bv^{8,1}, \quad
K_2 \bv^{8,1} = 0, \\
P_1 \bv^{8,1} = - \frac{4}{3}\bv^{8,2}, \quad
(J_1 P)_2 \bv^{8,1} = \frac{104}{9}\sqrt{-3}\bv^{8,2}, \quad
(K_1 P)_2 \bv^{8,1} = 0,
\end{gather*}
\begin{gather*}
(\tom^1)_1 \bv^{8,2} = \big(\frac{2}{45} +
\frac{1}{3}\big) \bv^{8,2}, \quad
(\tom^2)_1 \bv^{8,2} = \frac{2}{5} \bv^{8,2}, \quad
J_2 \bv^{8,2} = - \frac{22}{81}\sqrt{-3} \bv^{8,2}, \quad
K_2 \bv^{8,2} = 0, \\
P_1 \bv^{8,2} = 2 \bv^{8,1}, \quad
(J_1 P)_2 \bv^{8,2} = - \frac{52}{3}\sqrt{-3} \bv^{8,1}, \quad
(K_1 P)_2 \bv^{8,2} = - \frac{20}{3} \sqrt{-3} \bv^{8,2}.
\end{gather*}

\subsection{The simple module $V_L^{T_{\chi_0}} (\tau)(2)$}
By \eqref{eq:Twisted-1-vep}, $V_L^{T_{\chi_0}} (\tau)(2) \cong
(M_T(\tau)(2) \ot M_t^0) \op (W_T(\tau)(2) \ot W_t^0)$ as
$M^0$-modules. The top level of $V_L^{T_{\chi_0}} (\tau)(2)$ is
$\C\bv^9$, where $\bv^9 = h_2(-1/3) \ot v \in W_T(\tau)(2) \ot
W_t^0$ and $0 \ne v \in T_{\chi_0}$. The following hold.
\begin{gather*}
(\tom^1)_1 \bv^9 = \frac{2}{45} \bv^9, \quad
(\tom^2)_1 \bv^9 = \frac{2}{5} \bv^9, \quad
J_2 \bv^9 = - \frac{4}{81}\sqrt{-3} \bv^9, \quad
K_2 \bv^9 = 0, \\
P_1 \bv^9 = 0, \quad
(J_1 P)_2 \bv^9 = 0, \quad
(K_1 P)_2 \bv^9 = \frac{4}{3} \sqrt{-3} \bv^9.
\end{gather*}

\subsection{The simple module $V_L^{T_{\chi_j}} (\tau)(0)$, $j=1,2$}
By \eqref{eq:Twisted-1-vep}, $V_L^{T_{\chi_j}} (\tau)(0) \cong
(M_T(\tau)(2) \ot M_t^{3-j}) \op (W_T(\tau)(2) \ot W_t^{3-j})$ as
$M^0$-modules for $j=1,2$. The top level of $V_L^{T_{\chi_j}}
(\tau)(0)$ is $\C\bv^{10,j}$, where $\bv^{10,j} = 1 \ot v \in
W_T(\tau)(2) \ot W_t^{3-j}$ and $0 \ne v \in T_{\chi_j}$. The
following hold.
\begin{gather*}
(\tom^1)_1 \bv^{10,j} = \frac{2}{45} \bv^{10,j}, \quad
(\tom^2)_1 \bv^{10,j} = \frac{1}{15} \bv^{10,j}, \quad
J_2 \bv^{10,j} = - \frac{4}{81}\sqrt{-3} \bv^{10,j}, \\
K_2 \bv^{10,j} = -(-1)^j \frac{2}{9} \bv^{10,j}, \quad
P_1 \bv^{10,j} = (-1)^j \frac{1}{9}\sqrt{-3} \bv^{10,j},\\
(J_1 P)_2 \bv^{10,j} = (-1)^j \frac{8}{9} \bv^{10,j}, \quad
(K_1 P)_2 \bv^{10,j} = - \frac{2}{9} \sqrt{-3} \bv^{10,j}.
\end{gather*}

\subsection{The simple module $V_L^{T_{\chi_j}} (\tau)(1)$, $j=1,2$}
By \eqref{eq:Twisted-1-vep}, $V_L^{T_{\chi_j}} (\tau)(1) \cong
(M_T(\tau)(0) \ot M_t^{3-j}) \op (W_T(\tau)(0) \ot W_t^{3-j})$ as
$M^0$-modules for $j=1,2$. The top level of $V_L^{T_{\chi_j}}
(\tau)(1)$ is of dimension $2$ with basis $\{ \bv^{11,j,1},
\bv^{11,j,2}\}$, where
\begin{align*}
\bv^{11,j,1} & = h_1(-2/3) \ot v
- (-1)^j \sqrt{-3}h_2(-1/3)^2 \ot v
\in M_T(\tau)(0) \ot M_t^{3-j},\\
\bv^{11,j,2} & = 2h_1(-2/3) \ot v
+ (-1)^j \sqrt{-3}h_2(-1/3)^2 \ot v
\in W_T(\tau)(0) \ot W_t^{3-j}
\end{align*}
and $0 \ne v \in T_{\chi_j}$. The following hold.
\begin{gather*}
(\tom^1)_1 \bv^{11,j,1} = \frac{1}{9} \bv^{11,j,1}, \quad
(\tom^2)_1 \bv^{11,j,1} = \frac{2}{3} \bv^{11,j,1}, \quad
J_2 \bv^{11,j,1} = \frac{14}{81}\sqrt{-3} \bv^{11,j,1}, \\
K_2 \bv^{11,j,1} = (-1)^j \frac{52}{9}\bv^{11,j,1}, \quad
P_1 \bv^{11,j,1} = - (-1)^j \frac{4}{9} \sqrt{-3}\bv^{11,j,2},\\
(J_1 P)_2 \bv^{11,j,1} = (-1)^j \frac{52}{9} \bv^{11,j,2}, \\
\quad (K_1 P)_2 \bv^{11,j,1} = - \frac{28}{9} \sqrt{-3}
\bv^{11,j,2},
\end{gather*}
\begin{gather*}
(\tom^1)_1 \bv^{11,j,2} = \big(\frac{2}{45} +
\frac{2}{3}\big) \bv^{11,j,2}, \quad
(\tom^2)_1 \bv^{11,j,2} = \frac{1}{15} \bv^{11,j,2}, \quad
J_2 \bv^{11,j,2} = \frac{176}{81}\sqrt{-3} \bv^{11,j,2}, \\
K_2 \bv^{11,j,2} = - (-1)^j \frac{2}{9}\bv^{11,j,2}, \quad
P_1 \bv^{11,j,2} = - (-1)^j \frac{8}{9} \sqrt{-3}\bv^{11,j,1}
+ (-1)^j \frac{5}{9} \sqrt{-3} \bv^{11,j,2}, \\
(J_1 P)_2 \bv^{11,j,2} = (-1)^j \frac{104}{9} \bv^{11,j,1}
- (-1)^j \frac{200}{9} \bv^{11,j,2}, \\
\quad (K_1 P)_2 \bv^{11,j,2} = - \frac{56}{9}\sqrt{-3}
\bv^{11,j,1} - \frac{10}{9} \sqrt{-3} \bv^{11,j,2}.
\end{gather*}

\subsection{The simple module $V_L^{T_{\chi_j}} (\tau)(2)$, $j=1,2$}
By \eqref{eq:Twisted-1-vep}, $V_L^{T_{\chi_j}} (\tau)(2) \cong
(M_T(\tau)(1) \ot M_t^{3-j}) \op (W_T(\tau)(1) \ot W_t^{3-j})$  as
$M^0$-modules for $j=1,2$. The top level of $V_L^{T_{\chi_j}}
(\tau)(2)$ is $\C\bv^{12,j}$, where $\bv^{12,j} = h_2(-1/3) \ot v
\in W_T(\tau)(1) \ot W_t^{3-j}$ and $0 \ne v \in T_{\chi_j}$. The
following hold.
\begin{gather*}
(\tom^1)_1 \bv^{12,j} =
\big( \frac{2}{45} + \frac{1}{3}\big) \bv^{12,j}, \quad
(\tom^2)_1 \bv^{12,j} = \frac{1}{15} \bv^{12,j}, \quad
J_2 \bv^{12,j} = - \frac{22}{81}\sqrt{-3} \bv^{12,j}, \\
K_2 \bv^{12,j} = - (-1)^j \frac{2}{9} \bv^{12,j}, \quad
P_1 \bv^{12,j} = - (-1)^j \frac{5}{9}\sqrt{-3} \bv^{12,j},\\
(J_1 P)_2 \bv^{12,j} = (-1)^j \frac{8}{9} \bv^{12,j}, \quad
(K_1 P)_2 \bv^{12,j} = \frac{10}{9} \sqrt{-3} \bv^{12,j}.
\end{gather*}

\subsection{The simple module $V_L^{T_{\chi'_0}} (\tau^2)(0)$}
By \eqref{eq:Twisted-2-vep}, $V_L^{T_{\chi'_0}} (\tau^2)(0) \cong
(M_T(\tau^2)(0) \ot M_t^0) \op (W_T(\tau^2)(0) \ot W_t^0)$ as
$M^0$-modules. The top level of $V_L^{T_{\chi'_0}} (\tau^2)(0)$ is
$\C\bv^{13}$, where $\bv^{13} = 1 \ot v \in M_T(\tau^2)(0) \ot
M_t^0$ and $0 \ne v \in T_{\chi'_0}$. The following hold.
\begin{gather*}
(\tom^1)_1 \bv^{13} = \frac{1}{9} \bv^{13}, \quad
(\tom^2)_1 \bv^{13} = 0, \quad
J_2 \bv^{13} = - \frac{14}{81}\sqrt{-3} \bv^{13}, \quad
K_2 \bv^{13} = 0,\\
P_1 \bv^{13} = 0, \quad
(J_1 P)_2 \bv^{13} = 0, \quad (K_1 P)_2 \bv^{13} = 0.
\end{gather*}

\subsection{The simple module $V_L^{T_{\chi'_0}} (\tau^2)(1)$}
By \eqref{eq:Twisted-2-vep}, $V_L^{T_{\chi'_0}} (\tau^2)(1) \cong
(M_T(\tau^2)(1) \ot M_t^0) \op (W_T(\tau^2)(1) \ot W_t^0)$ as
$M^0$-modules. The top level of $V_L^{T_{\chi'_0}} (\tau^2)(1)$ is
of dimension $2$ with basis $\{ \bv^{14,1}, \bv^{14,2}\}$, where
$\bv^{14,1} = h_2'(-1/3)^2 \ot v \in M_T(\tau^2)(1) \ot M_t^0$,
$\bv^{14,2} = h_1'(-2/3) \ot v \in W_T(\tau^2)(1) \ot W_t^0$ and $0
\ne v \in T_{\chi_0'}$. The following hold.
\begin{gather*}
(\tom^1)_1 \bv^{14,1} = \big(\frac{1}{9} +
\frac{2}{3} \big) \bv^{14,1}, \quad
(\tom^2)_1 \bv^{14,1} = 0, \quad
J_2 \bv^{14,1} = \frac{238}{81}\sqrt{-3} \bv^{14,1},\quad
K_2 \bv^{14,1} = 0,\\
P_1 \bv^{14,1} = - \frac{4}{3}\bv^{14,2},\quad
(J_1 P)_2 \bv^{14,1} = - \frac{104}{9} \sqrt{-3} \bv^{14,2},
\quad (K_1 P)_2 \bv^{14,1} = 0,
\end{gather*}
\begin{gather*}
(\tom^1)_1 \bv^{14,2} = \big(\frac{2}{45} +
\frac{1}{3} \big) \bv^{14,2}, \quad
(\tom^2)_1 \bv^{14,2} = \frac{2}{5} \bv^{14,2}, \quad
J_2 \bv^{14,2} = \frac{22}{81}\sqrt{-3} \bv^{14,2}, \quad
K_2 \bv^{14,2} = 0,\\
P_1 \bv^{14,2} = 2 \bv^{14,1}, \quad
(J_1 P)_2 \bv^{14,2} = \frac{52}{3} \sqrt{-3} \bv^{14,1},
\quad (K_1 P)_2 \bv^{14,2} = \frac{20}{3} \sqrt{-3} \bv^{14,2}.
\end{gather*}

\subsection{The simple module $V_L^{T_{\chi'_0}} (\tau^2)(2)$}
By \eqref{eq:Twisted-2-vep}, $V_L^{T_{\chi_0'}} (\tau^2)(2) \cong
(M_T(\tau^2)(2) \ot M_t^0) \op (W_T(\tau^2)(2) \ot W_t^0)$ as
$M^0$-modules. The top level of $V_L^{T_{\chi'_0}} (\tau^2)(2)$ is
$\C\bv^{15}$, where $\bv^{15} = h_2'(-1/3) \ot v \in W_T(\tau^2)(2)
\ot W_t^0$ and $0 \ne v \in T_{\chi_0'}$. The following hold.
\begin{gather*}
(\tom^1)_1 \bv^{15} = \frac{2}{45} \bv^{15}, \quad
(\tom^2)_1 \bv^{15} = \frac{2}{5} \bv^{15}, \quad
J_2 \bv^{15} = \frac{4}{81}\sqrt{-3} \bv^{15}, \quad
K_2 \bv^{15} = 0,\\
P_1 \bv^{15} = 0, \quad
(J_1 P)_2 \bv^{15} = 0, \quad
(K_1 P)_2 \bv^{15} = - \frac{4}{3} \sqrt{-3} \bv^{15}.
\end{gather*}

\subsection{The simple module $V_L^{T_{\chi'_j}} (\tau^2)(0)$, $j=1,2$}
By \eqref{eq:Twisted-2-vep}, $V_L^{T_{\chi'_j}} (\tau^2)(0) \cong
(M_T(\tau^2)(2) \ot M_t^j) \op (W_T(\tau^2)(2) \ot W_t^j)$ as
$M^0$-modules for $j=1,2$. The top level of $V_L^{T_{\chi'_j}}
(\tau^2)(0)$ is $\C\bv^{16,j}$, where $\bv^{16,j} = 1 \ot v \in
W_T(\tau^2)(2) \ot W_t^j$ and $0 \ne v \in T_{\chi'_j}$. The
following hold.
\begin{gather*}
(\tom^1)_1 \bv^{16,j} = \frac{2}{45} \bv^{16,j}, \quad
(\tom^2)_1 \bv^{16,j} = \frac{1}{15} \bv^{16,j}, \quad
J_2 \bv^{16,j} = \frac{4}{81}\sqrt{-3} \bv^{16,j}, \\
K_2 \bv^{16,j} = (-1)^j \frac{2}{9} \bv^{16,j}, \quad
P_1 \bv^{16,j} = (-1)^j \frac{1}{9}\sqrt{-3} \bv^{16,j},\\
(J_1 P)_2 \bv^{16,j} = - (-1)^j \frac{8}{9} \bv^{16,j}, \quad
(K_1 P)_2 \bv^{16,j} = \frac{2}{9} \sqrt{-3} \bv^{16,j}.
\end{gather*}

\subsection{The simple module $V_L^{T_{\chi'_j}} (\tau^2)(1)$, $j=1,2$}
By \eqref{eq:Twisted-2-vep}, $V_L^{T_{\chi'_j}} (\tau^2)(1) \cong
(M_T(\tau^2)(0) \ot M_t^j) \op (W_T(\tau^2)(0) \ot W_t^j)$ as
$M^0$-modules for $j=1,2$. The top level of $V_L^{T_{\chi'_j}}
(\tau^2)(1)$ is of dimension $2$ with basis $\{ \bv^{17,j,1},
\bv^{17,j,2}\}$, where
\begin{align*}
\bv^{17,j,1} & = h_1'(-2/3) \ot v
- (-1)^j \sqrt{-3}h_2'(-1/3)^2 \ot v
\in M_T(\tau^2)(0) \ot M_t^j,\\
\bv^{17,j,2} & = 2h_1'(-2/3) \ot v
+ (-1)^j \sqrt{-3}h_2'(-1/3)^2 \ot v
\in W_T(\tau^2)(0) \ot W_t^j
\end{align*}
and $0 \ne v \in T_{\chi'_j}$. The following hold.
\begin{gather*}
(\tom^1)_1 \bv^{17,j,1} = \frac{1}{9} \bv^{17,j,1}, \quad
(\tom^2)_1 \bv^{17,j,1} = \frac{2}{3} \bv^{17,j,1}, \quad
J_2 \bv^{17,j,1} = - \frac{14}{81}\sqrt{-3} \bv^{17,j,1}, \\
K_2 \bv^{17,j,1} = - (-1)^j \frac{52}{9} \bv^{17,j,1}, \quad
P_1 \bv^{17,j,1} = - (-1)^j \frac{4}{9} \sqrt{-3}\bv^{17,j,2},\\
(J_1 P)_2 \bv^{17,j,1} = - (-1)^j \frac{52}{9} \bv^{17,j,2},\\
(K_1 P)_2 \bv^{17,j,1} = \frac{28}{9} \sqrt{-3} \bv^{17,j,2},
\end{gather*}
\begin{gather*}
(\tom^1)_1 \bv^{17,j,2} = \big(\frac{2}{45} +
\frac{2}{3} \big) \bv^{17,j,2}, \quad
(\tom^2)_1 \bv^{17,j,2} = \frac{1}{15} \bv^{17,j,2}, \quad
J_2 \bv^{17,j,2} = - \frac{176}{81}\sqrt{-3} \bv^{17,j,2}, \\
K_2 \bv^{17,j,2} = (-1)^j \frac{2}{9} \bv^{17,j,2}, \quad
P_1 \bv^{17,j,2} = - (-1)^j \frac{8}{9} \sqrt{-3}\bv^{17,j,1}
+ (-1)^j \frac{5}{9} \sqrt{-3} \bv^{17,j,2},\\
(J_1 P)_2 \bv^{17,j,2} = - (-1)^j \frac{104}{9} \bv^{17,j,1} +
(-1)^j \frac{200}{9} \bv^{17,j,2},\\
(K_1 P)_2 \bv^{17,j,2} = \frac{56}{9} \sqrt{-3} \bv^{17,j,1} +
\frac{10}{9} \sqrt{-3} \bv^{17,j,2}.
\end{gather*}

\subsection{The simple module $V_L^{T_{\chi'_j}} (\tau^2)(2)$, $j=1,2$}
By \eqref{eq:Twisted-2-vep}, $V_L^{T_{\chi'_j}} (\tau^2)(2) \cong
(M_T(\tau^2)(1) \ot M_t^j) \op (W_T(\tau^2)(1) \ot W_t^j)$ as
$M^0$-modules for $j=1,2$. The top level of $V_L^{T_{\chi'_j}}
(\tau^2)(2)$ is $\C\bv^{18,j}$, where $\bv^{18,j} = h_2'(-1/3) \ot v
\in W_T(\tau^2)(1) \ot W_t^j$ and $0 \ne v \in T_{\chi'_j}$. The
following hold.
\begin{gather*}
(\tom^1)_1 \bv^{18,j} =
\big( \frac{2}{45} + \frac{1}{3} \big) \bv^{18,j}, \quad
(\tom^2)_1 \bv^{18,j} = \frac{1}{15} \bv^{18,j}, \quad
J_2 \bv^{18,j} = \frac{22}{81}\sqrt{-3} \bv^{18,j}, \\
K_2 \bv^{18,j} = (-1)^j \frac{2}{9} \bv^{18,j}, \quad
P_1 \bv^{18,j} = - (-1)^j \frac{5}{9}\sqrt{-3} \bv^{18,j}, \\
(J_1 P)_2 \bv^{18,j} = - (-1)^j \frac{8}{9} \bv^{18,j}, \quad
(K_1 P)_2 \bv^{18,j} = - \frac{10}{9} \sqrt{-3} \bv^{18,j}.
\end{gather*}

\subsection{Symmetries by $\sg$}
Let us consider the automorphisms $\sg$ and $\theta$ of $V_L$ which
are lifts of the isometries $\sg$ and $\theta$ of the lattice $L$
defined by \eqref{eq:3Aut}. Clearly, $\sg \tau \sg = \tau^2$, $\sg
\theta = \theta \sg$, and $\tau \theta = \theta \tau$. Thus $\sg$
and $\theta$ induce automorphisms of $V_L^\tau$ of order $2$. We
have $\sg J = -J$, $\sg K = -K$, $\sg P = P$, $\theta J = J$,
$\theta K = -K$, and $\theta P = -P$. Hence $\sg$ and $\theta$
induce the same automorphism of $M_t^0$ and $\theta$ is the identity
on $M(0)$. Note also that $\sg (J_1 P) = -J_1 P$ and $\sg (K_1 P) =
-K_1 P$.

>From the action of $\sg$ on the top level of the $30$ known simple
$V_L^\tau$-modules or the action of $J_2$, $K_2$, $(J_1 P)_2$, and
$(K_1 P)_2$, we know how $\sg$ permutes those simple
$V_L^\tau$-modules. In fact, $\sg$ transforms $V_{L^{(c,0)}}$ into
an equivalent simple $V_L^\tau$-module and interchanges the
remaining simple $V_L^\tau$-modules as follows.
\begin{align*}
V_L(1) & \leftrightarrow V_L(2),
& V_{L^{(0,1)}}(\vep) & \leftrightarrow
V_{L^{(0,2)}}(2\vep), \quad \vep = 0,1,2,\\
V_{L^{(c,1)}} & \leftrightarrow V_{L^{(c,2)}},
& V_L^{T_{\chi_j}} (\tau)(\vep) & \leftrightarrow
V_L^{T_{\chi'_j}} (\tau^2)(\vep), \quad j, \vep =0,1,2.
\end{align*}

Note that $\sg h_i = \xi^{3-i}h'_i$, $i=1,2$. The top level of
$V_L^{T_{\chi'_j}} (\tau^2)(\vep)$ can be obtained by replacing
$h_i(i/3 + n)$ with $h'_i(i/3 + n)$ in the top level of
$V_L^{T_{\chi_j}} (\tau)(\vep)$ for $j,\vep=0,1,2$. The
corresponding action of $\sg$ on the simple $M(0)$-modules was
discussed in \cite[Section 4.4]{DLTYY}.

\section{Classification of simple modules}
\label{sec:Classification} We keep the notation in the preceding
section. Thus $V_L^\tau = M^0 \op W^0$ with $M^0 = M(0) \ot M_t^0$
and $W^0 = W(0) \ot W_t^0$. In this section we show that any simple
$V_L^\tau$-module is equivalent to one of the $30$ simple
$V_L^\tau$-modules listed in Lemma \ref{lem:30SIMPLES}. The result
will be established by considering the Zhu algebra $A(V_L^\tau)$ of
$V_L^\tau$.

First, we review some notation and basic formulas for the Zhu
algebra $A(V)$ of a vertex operator algebra $(V,Y, \1, \om)$. Define
two binary operations
\begin{equation}\label{eq:Zhu-Operation}
u \ast v = \sum_{i=0}^{\infty} \binom{\wt u}{i}u_{i-1}v, \quad u
\circ v = \sum_{i=0}^{\infty} \binom{\wt u}{i}u_{i-2}v
\end{equation}
for $u,v \in V$ with $u$ being homogeneous and extend $\ast$ and
$\circ$ for arbitrary $u \in V$ by linearity. Let $O(V)$ be the
subspace of $V$ spanned by all $u \circ v$ for $u,v \in V$. Set
$A(V) = V/O(V)$. By \cite[Theorem 2.1.1]{Z}, $O(V)$ is a two-sided
ideal with respect to the operation $\ast$. Thus $\ast$ induces an
operation in $A(V)$. Denote by $[v]$ the image of $v \in V$ in
$A(V)$. Then $[u]\ast [v] = [u\ast v]$ and $A(V)$ is an associative
algebra by this operation. Moreover, $[\1]$ is the identity and
$[\om]$ is in the center of $A(V)$. For $u,v \in V$, we write $u
\sim v$ if $[u] = [v]$. For $\varphi,\psi \in \End V$, we write
$\varphi \sim \psi$ if $\varphi v \sim \psi v$ for all $v \in V$. We
need some basic formulas (cf. \cite{Z}).
\begin{equation}\label{eq:vu}
v \ast u \sim \sum_{i=0}^{\infty} \binom{\wt (u) - 1}{i}u_{i-1}v,
\end{equation}
\begin{equation}\label{gen}
\sum_{i=0}^{\infty} \binom{\wt (u) + m}{i}u_{i-n-2}v \in O(V),
\quad n \ge m \ge 0.
\end{equation}

Moreover (cf. \cite{Wang1}),
\begin{equation}\label{eq:virn}
L(-n) \sim (-1)^n \big\{ (n-1)\big( L(-2)+L(-1)\big) + L(0)\big\},
\quad n \ge 1,
\end{equation}
\begin{equation}\label{eq:virp}
[\om] \ast [u] = [(L(-2)+L(-1))u],
\end{equation}
where $L(n) = \om_{n+1}$. From \eqref{eq:virn} and \eqref{eq:virp}
we have
\begin{equation}\label{f1}
[L(-n)u] = (-1)^n (n-1) [\om]\ast [u] + (-1)^n [L(0)u], \quad n
\ge 1.
\end{equation}

If $u \in V$ is of weight $2$, then $u(-n-3) + 2u(-n-2) + u(-n-1)
\sim 0$ by \eqref{gen}, where $u(n) = u_{n+1}$. Hence
\begin{equation}\label{eq:Weight2-1}
u(-n) \sim (-1)^n \big( (n-1)u(-2) + (n-2)u(-1)\big)
\end{equation}
for $n \ge 1$. Then it follows from \eqref{eq:Zhu-Operation} and
\eqref{eq:vu} that
\begin{equation}\label{eq:Weight2-2}
u(-n) w \sim (-1)^n \big( -u \ast w + nw \ast u + u(0)w \big)
\end{equation}
for $n \ge 1$, $w \in V$. Likewise, if $u$ is of weight $3$ and
$u(n) = u_{n+2}$, then
\begin{equation}\label{eq:Weight3-1}
\begin{split}
u(-n) \sim (-1)^{n+1} & \Big( \frac{1}{2}(n-1)(n-2)u(-3)
+ (n-1)(n-3)u(-2) \\
& \quad + \frac{1}{2}(n-2)(n-3)u(-1)\Big),
\end{split}
\end{equation}
\begin{equation}\label{eq:Weight3-2}
\begin{split}
u(-n)w \sim (-1)^{n+1} & \Big( nu(-1)w + (n-1) u(0)w \\
& \quad -(n-1)u \ast w + \frac{1}{2}n(n-1) w \ast u \Big)
\end{split}
\end{equation}
for $n \ge 1$, $w \in V$.

For a homogeneous vector $u \in V$, $o(u)=u_{\wt (u) - 1}$ is the
weight zero component operator of $Y(u,z)$. Extend $o(u)$ for
arbitrary $u \in V$ by linearity. Note that we call a module in the
sense of \cite{Z} an $\N$-graded weak module here. If $M =
\oplus_{n=0}^{\infty} M_{(n)}$ is an $\N$-graded weak $V$-module
with $M_{(0)} \ne 0,$ then $o(u)$ acts on its top level $M_{(0)}$.
Zhu's theory \cite{Z} says: (1) $o(u)o(v) = o(u \ast v)$ as
operators on the top level $M_{(0)}$ and $o(u)$ acts as $0$ on
$M_{(0)}$ if $u \in O(V)$. Thus $M_{(0)}$ is an $A(V)$-module, where
$[u]$ acts on $M_{(0)}$ as $o(u)$. (2) The map $M \mapsto M_{(0)}$
is a bijection between the set of isomorphism classes of simple
$\N$-graded weak $V$-modules and the set of isomorphism classes of
simple $A(V)$-modules.

Let us return to our $V_L^\tau$. As in Section \ref{sec:VL-tau}, we
write $L^i(n) = (\tom^i)_{n+1}$, $i=1,2$, $J(n) = J_{n+2}$, and
$K(n) = K_{n+2}$. The Zhu algebras $A(M(0))$ and $A(M_t^0)$ were
determined in \cite{DLTYY} and \cite{KMY}, respectively. Since
$O(M^0) \subset O(V_L^\tau)$, the image of $M(0)$ (resp. $M_t^0$) in
$A(V_L^\tau)$ is a homomorphic image of $A(M(0))$ (resp.
$A(M_t^0)$). It is generated by $[\tom^1]$ and $[J]$ (resp.
$[\tom^2]$ and $[K]$).

By a direct calculation, we have
\begin{equation}\label{eq:PnP}
\begin{split}
P_1 P &= -16\tom^1 - 6\tom^2,\\
P_0 P &= -8(\tom^1)_{-2}\1 - 3(\tom^2)_{-2}\1,\\
P_{-1} P &= \frac{5}{273}J_1 K_1 P - \frac{12}{7}(\tom^1)_{-3}\1
-\frac{18}{13}(\tom^2)_{-3}\1\\
& \qquad -\frac{36}{7}(\tom^1)_{-1}(\tom^1)_{-1}\1 - \frac{9}{13}
(\tom^2)_{-1}(\tom^2)_{-1}\1 - 16(\tom^1)_{-1}(\tom^2)_{-1}\1,\\
P_{-2} P &= \frac{1}{84}J_0 K_1 P + \frac{1}{156}
J_1 K_0 P\\
& \qquad -\frac{8}{7}(\tom^1)_{-4}\1 -
\frac{12}{13}(\tom^2)_{-4}\1
-\frac{36}{7}(\tom^1)_{-2}(\tom^1)_{-1}\1 -
\frac{9}{13}(\tom^2)_{-2}(\tom^2)_{-1}\1\\
& \qquad -8(\tom^1)_{-2}(\tom^2)_{-1}\1 -
8(\tom^1)_{-1}(\tom^2)_{-2}\1.
\end{split}
\end{equation}
Moreover, $J_2P = K_2P = 0$. Then, using the formulas
\eqref{eq:virn}--\eqref{eq:Weight3-2}, we obtain
\begin{equation}\label{eq:PP}
\begin{split}
[P] \ast [P] &=  \frac{5}{273}[J_1 K_1 P] - \frac{36}{7}[\tom^1]
\ast [\tom^1] - \frac{9}{13}[\tom^2] \ast [\tom^2]\\
& \qquad - 16[\tom^1]\ast [\tom^2] + \frac{4}{7}[\tom^1] +
\frac{6}{13}[\tom^2],
\end{split}
\end{equation}
\begin{equation}\label{eq:P-circ-P}
\begin{split}
[P \circ P ] &= \frac{1}{84}[J]\ast [K_1 P]
- \frac{1}{84}[K_1 P] \ast [J] + \frac{1}{156}[K]\ast [J_1 P]
- \frac{1}{156}[J_1 P]\ast [K]\\
&= 0.
\end{split}
\end{equation}

It turns out that $A(V_L^\tau)$ is generated by $[\tom^1]$,
$[\tom^2]$, $[J]$, $[K]$, and $[P]$ (cf. Corollary
\ref{cor:generator}). However, we first prove the following
intermediate assertion.

\begin{prop}
The Zhu algebra $A(V_L^\tau)$ is generated by $[\tom^1]$,
$[\tom^2]$, $[J]$, $[K]$, $[P]$, $[J_1 P]$, and $[K_1 P]$.
\end{prop}

\begin{proof}
Recall that $L^i(n)P = 0$ for $i=1,2$, $n \ge 1$, $L^1(0)P =
(8/5)P$, $L^2(0)P = (2/5)P$, and $J(n)P=K(n)P=0$ for $n \ge 0$. Thus
from the commutation relations \eqref{eq:L1L1}--\eqref{eq:JJ} and
\eqref{eq:L2L2}--\eqref{eq:KK} we see that $W^0$ is spanned by the
vectors of the form
\begin{equation}\label{eq:Normal-form}
L^1(-j_1) \cdots L^1(-j_r) L^2(-k_1) \cdots L^2(-k_s)
J(-m_1)\cdots J(-m_p) K(-n_1) \cdots K(-n_q)P
\end{equation}
with $j_1 \ge \cdots \ge j_r \ge 1$, $k_1 \ge \cdots \ge k_s \ge 1$,
$m_1 \ge \cdots \ge m_p \ge 1$, $n_1 \ge \cdots \ge n_q \ge 1$. Let
$v$ be a vector of this form. Its weight is
\begin{equation*}
j_1 + \cdots +j_r + k_1 + \cdots +k_s +
m_1 + \cdots +m_p + n_1 + \cdots +n_q + 2.
\end{equation*}

Since $V_L^\tau = M^0 \op W^0$ and since the image of $M^0$ in
$A(V_L^\tau)$ is generated by $[\tom^1]$, $[\tom^2]$, $[J]$, and
$[K]$, it suffices to show that the image $[v]$ of $v$ in
$A(V_L^\tau)$ is contained in the subalgebra generated by
$[\tom^1]$, $[\tom^2]$, $[J]$, $[K]$, $[P]$, $[J_1 P]$, and $[K_1
P]$. We proceed by induction on the weight of $v$. By the formula
\eqref{eq:Weight2-2} with $u=\tom^i$, $i=1,2$ and the induction on
the weight, we may assume that $r=s=0$, that is,
\begin{equation*}
v = J(-m_1)\cdots J(-m_p) K(-n_1) \cdots K(-n_q)P.
\end{equation*}

Moreover, by the formula \eqref{eq:Weight3-2} with $u=J$, we may
assume that $m_1 = \cdots = m_p = 1$. Since $J(m)$ and $K(n)$
commute, we may also assume that $n_1 = \cdots = n_q = 1$ by a
similar argument. Then $v = J(-1)^p K(-1)^q P$.

Next, we reduce $v$ to the case $p \le 1$. For this purpose, we use
a singular vector
\begin{equation}\label{eq:W0-singular}
5J(-1)^2P + 2496L^1(-2)P - 195L^1(-1)^2 P = 0.
\end{equation}
in $W(0)$. Suppose $ p \ge 2$. Then, since $K(-1)$ commute with
$J(m)$ and $L^1(n)$, \eqref{eq:W0-singular} implies that $v =
J(-1)^p K(-1)^q P$ is a linear combination of
$J(-1)^{p-2}L^1(-2)K(-1)^q P$ and $J(-1)^{p-2}L^1(-1)^2 K(-1)^q P$.
By the commutation relation \eqref{eq:L1J}, these two vectors can be
written in the form $L^1(-2)HK(-1)^q P$ and $L^1(-1)^2 H'K(-1)^q P$,
where $H$ (resp. $H'$) is a polynomial in $J(-1)$ and $J(-3)$ (resp.
$J(-1)$, $J(-2)$, and $J(-3)$). Then by the formula
\eqref{eq:Weight2-2} with $u = \tom^1$ and the induction on the
weight, the assertion holds for $v$. Hence we may assume that $p \le
1$.

There is a singular vector
\begin{equation}\label{eq:Wt0-singular}
K(-1)^2 P - 210L^2(-2)P = 0
\end{equation}
in $W_t^0$. Thus, by a similar argument as above, we may assume that
$q \le 1$. Finally, it follows from \eqref{eq:PP} that $[J(-1) K(-1)
P]$ can be written by $[\tom^1]$, $[\tom^2]$, and $[P]$ in
$A(V_L^\tau)$. The proof is complete.
\end{proof}

Let us classify the simple $V_L^\tau$-modules. Our argument is based
on the knowledge of simple modules for $M(0)$ and $M_t^0$ together
with fusion rules \eqref{eq:W0-fusion} and \eqref{eq:Ternary-fusion}
among them. Set
\begin{gather*}
\M_1 = \{ M(\vep), M_k^c, M_T(\tau^i)(\vep) \,|\,
i=1,2, \vep=0,1,2 \},\\
\W_1 = \{ W(\vep), W_k^c, W_T(\tau^i)(\vep) \,|\,
i=1,2, \vep=0,1,2 \},\\
\M_2 = \{ M_t^j \,|\, j = 0,1,2\}, \qquad
\W_2 = \{ W_t^j \,|\, j = 0,1,2\}.
\end{gather*}

Then $\M_1 \cup \W_1$ (resp. $\M_2 \cup \W_2$) is a complete set of
representatives of isomorphism classes of simple $M(0)$-modules
(resp. simple $M_t^0$-modules). A main point is that the fusion
rules of the following form hold.
\begin{equation}\label{eq:MW-fusion}
\begin{split}
W(0) \times M^1 = W^1 &, \qquad
W(0) \times W^1 = M^1 + W^1,\\
W_t^0 \times M^2 = W^2 &, \qquad W_t^0 \times W^2 = M^2 + W^2,
\end{split}
\end{equation}
where $M^i \in \M_i$, $i=1,2$, and $W^i \in \W_i$ is determined by
$M^i$ through the fusion rule $W(0) \times M^1 = W^1$ or $W_t^0
\times M^2 = W^2$.

Recall that $M^0$ is rational, $C_2$-cofinite, and of CFT type. Thus
every $\N$-graded weak $M^0$-module is a direct sum of simple
$M^0$-modules. As a result, every $\N$-graded weak $V_L^\tau$-module
is decomposed into a direct sum of simple $M^0$-modules, and in
particular $L(0) = \om_1$ acts semisimply on it. Each weight
subspace, that is, each eigenspace for $L(0)$ is not necessarily a
finite dimensional space. However, any simple weak $V_L^\tau$-module
is a simple ordinary $V_L^\tau$-module by \cite[Corollary 5.8]{ABD},
since $V_L^\tau$ is $C_2$-cofinite and of CFT type.

We note that
\begin{equation}\label{eq:W-cdot-W}
W^0 \cdot W^0 = V_L^\tau.
\end{equation}
Indeed, $W^0 \cdot W^0 = \spn \{ a_n b\,|\,a,b \in W^0, n \in \Z\}$
is an $M^0$-submodule of $V_L^\tau$ by \eqref{eq:Formula1}. Since
$P, J_1K_1P \in W^0$ and $\tom^1$, $\tom^2 \in M^0$, \eqref{eq:PnP}
implies that $W^0 \cdot W^0 = M^0 \op W^0$.

Each simple $M^0$-module is isomorphic to a tensor product $A \ot B$
of a simple $M(0)$-module $A$ and a simple $M_t^0$-module $B$. We
show that only restricted simple $M^0$-modules can appear in
$\N$-graded weak $V_L^\tau$-modules.

\begin{lem}\label{lem:MM-or-WW}
Let $U$ be an $\N$-graded weak $V_L^\tau$-module. Then any simple
$M^0$-submodule of $U$ is isomorphic to $M^1 \ot M^2$ or $W^1 \ot
W^2$ for some $M^i \in \M_i$ and $W^i \in \W_i$, $i=1,2$.
\end{lem}

\begin{proof}
Suppose $U$ contains a simple $M^0$-submodule $S^0 \cong M^1 \ot
W^2$ with $M^1 \in \M_1$ and $W^2 \in \W_2$. Let $S = V_L^\tau \cdot
S^0 = \spn\{ a_n w \,|\, a \in V_L^\tau, w \in S^0, n \in \Z\}$.
Then \eqref{eq:Formula1} implies that $S$ is the $\N$-graded weak
$V_L^\tau$-submodule of $U$ generated by $S^0$. By the construction
of $S$, the difference of any two eigenvalues of $L(0)$ in $S$ is an
integer. In fact, $S$ is an ordinary $V_L^\tau$-module by Remark
\ref{rmk:ordinary-submodule}.

If $v$ is a nonzero vector in $V_L^\tau$, then $v_n S^0 \ne 0$ for
some $n \in \Z$. Indeed, Lemma \ref{lem:Formula2} implies that the
set $\{ v \in V_L^\tau\,|\, v_n S^0 = 0 \mbox{ for all } n \in \Z\}$
is an ideal of $V_L^\tau$. It is in fact $0$, since $V_L^\tau$ is a
simple vertex operator algebra and $S^0$ is a simple $M^0$-module.
Then by the fusion rules \eqref{eq:MW-fusion}, a simple $M^0$-module
isomorphic to $W^1 \ot M^2$ or $W^1 \ot W^2$ must appear in $S$.
However, the difference of the minimal eigenvalues of $L(0)$ in $M^1
\ot W^2$ and $W^1 \ot M^2$, or in $M^1 \ot W^2$ and $W^1 \ot W^2$ is
not an integer. This is a contradiction. Thus $U$ does not contain a
simple $M^0$-submodule isomorphic to $M^1 \ot W^2$. By a similar
argument, we can also show that there is no simple $M^0$-submodule
isomorphic to $W^1 \ot M^2$ in $U$. Hence the assertion holds.
\end{proof}

Set $\M = \{ M^1 \ot M^2\,|\,M^i \in \M_i, i=1,2\}$ and $\W = \{
W^1\ot W^2\,|\,W^i \in \W_i, i=1,2\}$. Then each of $\M$ and $\W$
consists of $30$ inequivalent simple $M^0$-modules. The top level of
every simple $M^0$-module is of dimension one.

\begin{lem}\label{lem:one-dim-top-level}
If $U$ is a simple $\N$-graded weak $V_L^\tau$-module whose top
level is of dimension one, then $U$ is isomorphic to one of the $23$
known simple $V_L^\tau$-modules with one dimensional top level,
namely, $V_{L^{(0,j)}}(\vep)$, $j=0,1,2$, $\vep=0,1,2$,
$V_{L^{(c,j)}}$, $j=1,2$, $V_L^{T_{\chi_j}}(\tau)(\vep)$, $j=0,1,2$,
$\vep=0,2$, and $V_L^{T_{\chi'_j}}(\tau^2)(\vep)$, $j=0,1,2$,
$\vep=0,2$.
\end{lem}

\begin{proof}
Since $U$ is a direct sum of simple $M^0$-modules and since the top
level, say $U_\lambda$ of $U$ is assumed to be of dimension one, it
follows from Lemma \ref{lem:MM-or-WW} that $U_\lambda$ is isomorphic
to the top level of $M^1 \ot M^2$ or the top level of $W^1 \ot W^2$
as an $A(M^0)$-module for some $M^i \in \M_i$, $W^i \in \W_i$,
$i=1,2$. The Zhu algebra $A(M^0) \cong A(M(0)) \times A(M_t^0)$ is
commutative and the action of $A(M^0)$ on the top level of $M^1 \ot
M^2$ and the top level of $W^1 \ot W^2$ are known. Indeed, we know
all possible action of the elements $[\tom^1]$, $[\tom^2]$, $[J]$,
and $[K]$ of $A(V_L^\tau)$ on $U_\lambda$. Let $[\tom^1]$,
$[\tom^2]$, $[J]$, and $[K]$ act on $U_\lambda$ as scalars $a_1$,
$a_2$, $b_1$, and $b_2$, respectively. There are $60$ possible such
quadruplets $(a_1, a_2, b_1, b_2)$.

Let $[P]$, $[J_1 P]$, and $[K_1 P]$ act on $U_\lambda$ as scalars
$x_1$, $x_2$, and $x_3$, respectively. Then it follows from
\eqref{eq:PP} that $[J_1K_1 P]$ acts on $U_\lambda$ as a scalar
\begin{equation}\label{eq:JKP-scalar}
(273/5)x_1^2 + (1404/5)a_1^2+(189/5)a_2^2 \\
 + (4368/5)a_1 a_2 -
(156/5)a_1 - (126/5)a_2.
\end{equation}

>From Appendix \ref{app:Vectors-for-Zhu} and formulas
\eqref{eq:virn}--\eqref{eq:Weight3-2} we see that $[P \circ (J_1
P)] = 0$ gives
\begin{equation}\label{eq:P-circ-JP-scalar}
15 b_2 x_1 + 5 a_2 x_3 - 2 x_3 = 0.
\end{equation}
Likewise,
\begin{equation}\label{eq:P-circ-KP-scalar}
(15a_2 - 1) x_2 = 0,
\end{equation}
since $[P \circ (K_1 P)] = 0$.

Using \eqref{eq:JKP-scalar}, we can calculate $[(J_1 P) \ast (J_1
P)]$, $[(K_1 P) \ast (K_1 P)]$, and $[(J_1 P) \ast (K_1 P)]$ in a
similar way and verify that the following equations hold.
\begin{equation}\label{eq:JPJP-scalar}
\begin{split}
x_2^2 =& \Big( (229164/575)a_1 - (37856/425)a_2
+ 1669382/48875 \Big)x_1^2 \\
& - (56/85)b_2 x_2 - (4056/115)b_1 x_3
+ (348994464/107525)a_1^3 \\
& + (137149584/9775)a_1^2 a_2 - (1030224/1375)a_1^2
+ (7064876/9775)a_1 a_2^2 \\
& - (40788488/48875)a_1 a_2
+ (16160456/537625)a_1
- (419184/9775)a_2^3 \\
& - (200994/48875)a_2^2
+ (1065516/48875)a_2 - (3042/187)b_1^2,
\end{split}
\end{equation}
\begin{equation}\label{eq:KPKP-scalar}
\begin{split}
x_3^2 =& \Big(- (37044/575)a_1 - (5684/85)a_2
+ 741713/97750 \Big)x_1^2 \\
& + (28/221)b_2 x_2 + (216/115)b_1 x_3
- (54559344/107525)a_1^3 \\
& - (28217448/9775)a_1^2 a_2 + (254982/1375)a_1^2
- (25042724/9775)a_1 a_2^2 \\
& + (26308184/48875)a_1 a_2
- (8127098/537625)a_1
- (4775148/25415)a_2^3 \\
& + (188338017/1270750)a_2^2 - (9722139/635375)a_2 -
(180/187)b_1^2,
\end{split}
\end{equation}
\begin{equation}\label{eq:JPKP-scalar}
\begin{split}
x_2 x_3 =& \Big( - (864/5)a_1^2 + (1248/25)a_1 a_2
+ (1152/5)a_2^2 + (5904/125)a_1 \\
& + (184176/125)a_2 - 62112/625 \Big)x_1 - 36b_1 b_2.
\end{split}
\end{equation}

We have obtained a system of equations
\eqref{eq:P-circ-JP-scalar}--\eqref{eq:JPKP-scalar} for $x_1$,
$x_2$, $x_3$. We can solve this system of equations with respect to
the $60$ possible quadruplets $(a_1,a_2,b_1,b_2)$. Actually, there
is no solution for $37$ quadruplets of $(a_1,a_2,b_1,b_2)$. For each
of the remaining $23$ quadruplets $(a_1,a_2,b_1,b_2)$, the system of
equations possesses a unique solution $(x_1, x_2, x_3)$.
Furthermore, the $23$ sets $(a_1, a_2, b_1, b_2, x_1, x_2,x_3)$ of
values determined in this way coincide with the action of
$[\tom^1]$, $[\tom^2]$, $[J]$, $[K]$, $[P]$, $[J_1 P]$, and $[K_1
P]$ on the top level of the $23$ known simple $V_L^\tau$-modules
with one dimensional top level described in Section
\ref{sec:Structure}. Since $A(V_L^\tau)$ is generated by these seven
elements, this implies that $U_\lambda$ is isomorphic to the top
level of one of the $23$ simple $V_L^\tau$-modules listed in the
assertion as an $A(V_L^\tau)$-module. Thus the lemma holds by Zhu's
theorem.
\end{proof}

\begin{rmk}
We also obtain some equations for $x_1x_2$ and $x_1x_3$ from $[P
\ast (J_1P)]$ and $[P \ast (K_1P)]$. However, they are not
sufficient to determine $x_1$, $x_2$, and $x_3$.
\end{rmk}

\begin{lem}\label{lem:MM}
Every $\N$-graded weak $V_L^\tau$-module contains a simple
$M^0$-submodule isomorphic to a member of $\M$.
\end{lem}

\begin{proof}
Suppose false and let $U$ be an $\N$-graded weak $V_L^\tau$-module
which contains no simple $M^0$-submodule isomorphic to a member of
$\M$. Then by Lemma \ref{lem:MM-or-WW}, there is a simple
$M^0$-submodule $W$ in $U$ such that $W \cong W^1 \ot W^2$ for some
$W^i \in \W_i$, $i=1,2$. The top level of $W$, say $W_\lambda$ for
some $\lambda \in \Q$, is a one dimensional space. Take $0 \ne w \in
W_\lambda$ and let $S = V_L^\tau\cdot w = \{a_n w\,|\,a \in
V_L^\tau, n \in \Z\}$, which is an ordinary $V_L^\tau$-module by
\eqref{eq:Formula1} and Remark \ref{rmk:ordinary-submodule}. Since
$V_L^\tau = M^0 \op W^0$, it follows from our assumption and the
fusion rules \eqref{eq:MW-fusion} that $S$ is isomorphic to a direct
sum of finite number of copies of $W$ as an $M^0$-module. Thus
$[\tom^1]$, $[\tom^2]$, $[J]$, and $[K]$ act on the top level
$S_\lambda$ of $S$ as scalars, say $a_1$, $a_2$, $b_1$, and $b_2$,
respectively. Then by a similar calculation as in the proof of Lemma
\ref{lem:one-dim-top-level}, we see that $[P \circ (K_1 P)] = 0$
implies
\begin{equation}\label{eq:P-circ-KP-MM}
(15 a_2 - 1) o(J_1 P) = 0
\end{equation}
as an operator on the top level $S_\lambda$. Recall that $[u] \in
A(V_L^\tau)$ acts on $S_\lambda$ as $o(u) = u_{\wt (u) -1}$ for a
homogeneous vector $u$ of $V_L^\tau$. Furthermore, we can calculate
that
\begin{equation}
\begin{split}
o(J_1 P)o(P) - o(P)o(J_1 P) &= 0,\\
o(K_1 P)o(P) - o(P)o(K_1 P) &= \frac{2}{13}(15 a_2 - 1) o(J_1 P),\\
o(J_1 P)o(K_1 P) - o(K_1 P)o(J_1 P) &= \frac{96}{125}(15 a_2 -
1)(65a_1 + 100a_2 +441)o(P)
\end{split}
\end{equation}
as operators on $S_\lambda$.

By \eqref{eq:P-circ-KP-MM}, $15 a_2 - 1 = 0$ or $o(J_1 P) = 0$ and
so $o(P)$, $o(J_1 P)$, and $o(K_1 P)$ commute each other. Thus the
action of $A(V_L^\tau)$ on $S_\lambda$ is commutative. Hence we can
choose a one dimensional $A(V_L^\tau)$-submodule $T$ of $S_\lambda$.
Zhu's theory tells us that there is a simple $\N$-graded weak
$V_L^\tau$-module $R$ whose top level $R_\lambda$ is isomorphic to
$T$ as an $A(V_L^\tau)$-module. Since $\dim R_\lambda = 1$, $R$ is
isomorphic to one of the $23$ simple $V_L^\tau$-modules listed in
Lemma \ref{lem:one-dim-top-level}. In particular, $R$ contains a
simple $M^0$-submodule $M$ isomorphic to a member of $\M$. Now,
consider the $V_L^\tau$-submodule $V_L^\tau \cdot T$ of $S$
generated by $T$. By Lemma \ref{lem:surjective-LU}, there is a
surjective homomorphism of $V_L^\tau$-modules from $V_L^\tau \cdot
T$ onto $R$. Then $V_L^\tau \cdot T$ must contain a simple
$M^0$-submodule isomorphic to $M$. This contradicts our assumption.
The proof is complete.
\end{proof}

\begin{lem}\label{lem:VM}
Let $U$ be an $\N$-graded weak $V_L^\tau$-module and $M$ be a simple
$M^0$-submodule of $U$ such that $M \cong M^1 \ot M^2$ as
$M^0$-modules for some $M^i \in \M_i$, $i=1,2$. Then $V_L^\tau \cdot
M = \spn \{ a_n u \,|\, a \in V_L^\tau, u \in M, n \in \Z\}$ is a
simple $V_L^\tau$-module. Moreover, $V_L^\tau \cdot M = M \op W$,
where $W$ is a simple $M^0$-module isomorphic to $W^1 \ot W^2$ and
$W^i$, $i=1,2$ are determined from $M^i$ by the fusion rules $W(0)
\times M^1 = W^1$ and $W_t^0 \times M^2 = W^2$ of
\eqref{eq:MW-fusion}.
\end{lem}

\begin{proof}
By Remark \ref{rmk:ordinary-submodule}, $V_L^\tau \cdot M$ is an
ordinary $V_L^\tau$-module. Note that $V_L^\tau \cdot M = (M^0 +
W^0)\cdot M = M + W^0 \cdot M$. We see that $W^0 \cdot M \ne 0$ by a
similar argument as in the proof of Lemma \ref{lem:MM-or-WW}.
Actually, $W^0 \cdot(W^0 \cdot M) \supset (W^0 \cdot W^0) \cdot M =
V_L^\tau \cdot M$ (cf. Lemma \ref{lem:Formula2} and
\eqref{eq:W-cdot-W}) implies $W^0 \cdot M \ne 0$ also. Moreover,
$W^0 \cdot M$ is an $M^0$-module by \eqref{eq:Formula1}. Since $M^0$
is rational, $W^0 \cdot M$ is decomposed into a direct sum of simple
$M^0$-modules, say $W^0 \cdot M = \op_{\gamma \in \Gamma} S^\gamma$.
Let $W = W^1 \ot W^2$, where $W^i \in \W_i$, $i=1,2$ are determined
by the fusion rules $W(0) \times M^1 = W^1$ and $W_t^0 \times M^2 =
W^2$. The space $I_{M^0}\binom{W}{W^0\ M}$ of intertwining operators
of type $\binom{W}{W^0\ M}$ is of dimension one and each $S^\gamma$
is isomorphic to $W$.

We want to show that $| \Gamma | = 1$. Suppose $\Gamma$ contains at
least two elements and take $\gamma_1, \gamma_2 \in \Gamma$,
$\gamma_1 \ne \gamma_2$. Let $\psi : S^{\gamma_2} \rightarrow
S^{\gamma_1}$ be an isomorphism of $M^0$-modules and $p_\gamma : W^0
\cdot M \rightarrow S^\gamma$ be a projection. For $a \in W^0$ and
$u \in M$, set
\begin{equation*}
\G_{\gamma_1}(a,z)u = p_{\gamma_1}Y_U(a,z)u, \qquad
\G_{\gamma_2}(a,z)u = \psi p_{\gamma_2}Y_U(a,z)u,
\end{equation*}
where $Y_U(a,z)$ is the vertex operator of the $\N$-graded weak
$V_L^\tau$-module $U$. Then $\G_{\gamma_i}(\,\cdot\,,z)$, $i=1,2$
are nonzero members in the one dimensional space
$I_{M^0}\binom{W}{W^0\ M}$, so that $\mu \G_{\gamma_1}(\,\cdot\,,z)
= \G_{\gamma_2}(\,\cdot\,,z)$ for some $0 \ne \mu \in \C$. Let $0
\ne v \in S^{\gamma_1}$. Then $v \in W^0 \cdot M$ and so $v = \sum_j
(a^j)_{n_j} u^j$ for some $a^j \in W^0$, $u^j \in M$, $n_j \in \Z$.
Take the coefficients of $z^{-n_j-1}$ in both sides of $\mu
\G_{\gamma_1}(a^j,z)u^j = \G_{\gamma_2}(a^j,z)u^j$. Then $\mu
p_{\gamma_1}((a^j)_{n_j} u^j) = \psi p_{\gamma_2}((a^j)_{n_j} u^j)$.
Summing up both sides of the equation with respect to $j$, we have
$\mu p_{\gamma_1} v = \psi p_{\gamma_2} v$. However, $v \in
S^{\gamma_1}$ implies that $p_{\gamma_1} v = v$ and $p_{\gamma_2} v
= 0$. This is a contradiction since $\mu \ne 0$ and $v \ne 0$. Thus
$| \Gamma| = 1$ and $W^0 \cdot M \cong W$ as required.

If $V_L^\tau \cdot M$ is not a simple $V_L^\tau$-module, then there
is a proper $V_L^\tau$-submodule $N$ of $V_L^\tau \cdot M$. Since
$M$ and $W$ are simple $M^0$-modules, $N$ must be isomorphic to $M$
or $W$ as an $M^0$-module. Then the top level of $N$ is of dimension
one. The simple $V_L^\tau$-modules with one dimensional top level
are classified in Lemma \ref{lem:one-dim-top-level}. Each of them is
a direct sum of two simple $M^0$-modules. However, $N$ is not of
such a form. Thus $V_L^\tau \cdot M$ is a simple $V_L^\tau$-module.
\end{proof}

\begin{lem}\label{lem:Unique}
Let $U = M \op W$ be an $M^0$-module such that $M \cong M^1 \ot M^2$
and $W \cong W^1 \ot W^2$ for some $M^i \in \M_i$ and $W^i \in
\W_i$, $i=1,2$. Then $U$ admits at most one simple $V_L^\tau$-module
structure.
\end{lem}

\begin{proof}
Assume that $(U,Y_1)$ and $(U,Y_2)$ are simple $V_L^\tau$-modules
such that $Y_i(a,z) = Y(a,z)$ for all $a \in M^0$, $i=1,2$, where
$(U, Y)$ is the given $M^0$-mudule structure. We denote the vertex
operator of $V_L^\tau$ by $\ty(v,z)$ for $v \in V_L^\tau$. Let
$p_{M^0}:V_L^\tau \rightarrow M^0$ and $p_{W^0}:V_L^\tau \rightarrow
W^0$ be projections and define $\I(\,\cdot\,,z)$ and
$\J(\,\cdot\,,z)$ by
\begin{equation*}
\I(a,z)b = p_{M^0} \ty(a,z)b, \qquad \J(a,z)b = p_{W^0} \ty(a,z)b
\end{equation*}
for $a, b \in W^0$. Then by \eqref{eq:W-cdot-W}, $\I(\,\cdot\,,z)$
and $\J(\,\cdot\,,z)$ are nonzero intertwining operators of type
$\binom{M^0}{W^0\ W^0}$ and $\binom{W^0}{W^0\ W^0}$, respectively.
By the fusion rules \eqref{eq:MW-fusion}, the space
$I_{M^0}\binom{M^0}{W^0\ W^0}$ of $M^0$-intertwining operators of
type $\binom{M^0}{W^0\ W^0}$ is of dimension one. Likewise, $\dim
I_{M^0}\binom{W^0}{W^0\ W^0} = 1$. Note that $W^0 \cdot M^0 \subset
W^0$ and that $\I(a,z)b + \J(a,z)b = \ty(a,z)b$.

Let $p_M : U \rightarrow M$ and $p_W : U \rightarrow W$ be
projections. Define $\F_i^M(\,\cdot\,,z)$ and
$\F_i^W(\,\cdot\,,z)$, $i=1,2$ by
\begin{equation*}
\F^M_i(a,z)w = p_M Y_i(a,z)w, \qquad \F^W_i(a,z)w = p_W Y_i(a,z)w
\end{equation*}
for $a \in W^0$ and $w \in W$. Then $\F_i^M(\,\cdot\,,z)$ and
$\F_i^W(\,\cdot\,,z)$ are intertwining operators of type
$\binom{M}{W^0\ W}$ and $\binom{W}{W^0\ W}$, respectively. Clearly,
$\F_i^M(a,z)w + \F_i^W(a,z)w = Y_i(a,z)w$. If $\F_i^M(\,\cdot\,,z) =
0$, then $W^0 \cdot W \subset W$ and so $V_L^\tau \cdot W = M^0
\cdot W + W^0 \cdot W \subset W$. This is a contradiction, since $U$
is a simple $V_L^\tau$-module. Hence $\F_i^M(\,\cdot\,,z) \ne 0$.
Let
\begin{equation*}
\G_i^W(a,z)v = Y_i(a,z)v
\end{equation*}
for $a \in W^0$, $v \in M$. Then $\G_i^W(\,\cdot\,,z)$ is a nonzero
intertwining operator of type $\binom{W}{W^0\ M}$ by
\eqref{eq:MW-fusion}. The space of $M^0$-intertwining operators
$I_{M^0}\binom{W}{W^0\ M}$ of type $\binom{W}{W^0\ M}$ is of
dimension one by \eqref{eq:MW-fusion}. Similarly, $\dim
I_{M^0}\binom{M}{W^0\ W} = \dim I_{M^0}\binom{W}{W^0\ W} = 1$.
Therefore, $\F_2^M(\,\cdot\,,z) = \lambda \F_1^M(\,\cdot\,,z)$,
$\F_2^W(\,\cdot\,,z) = \mu \F_1^W(\,\cdot\,,z)$, and
$\G_2^W(\,\cdot\,,z) = \gamma \G_1^W(\,\cdot\,,z)$ for some
$\lambda, \mu , \gamma \in \C$ with $\lambda \ne 0$ and $\gamma \ne
0$.

Now,
\begin{align*}
Y_i(a,z_1)Y_i(b,z_2)v &= \big( \F_i^M(a,z_1) + \F_i^W(a,z_1)\big)
\G_i^W(b,z_2)v,\\
Y_i(b,z_2)Y_i(a,z_1)v &= \big( \F_i^M(b,z_2) + \F_i^W(b,z_2)\big)
\G_i^W(a,z_1)v,\\
Y_i(\ty(a,z_0)b,z_2)v &= Y_i(\I(a,z_0)b,z_2)v +
\G_i^W(\J(a,z_0)b,z_2)v
\end{align*}
for $a, b \in W^0$ and $v \in M$. Taking the image of both sides of
the Jacobi identity
\begin{equation}\label{eq:Jacobi-1}
\begin{split}
& z_0^{-1}\delta\Big(\frac{z_1-z_2}{z_0}\Big)Y_i(a,z_1)Y_i(b,z_2)v
- z_0^{-1}\delta\Big(\frac{z_2-z_1}{-z_0}\Big)Y_i(b,z_2)Y_i(a,z_1)v\\
& \qquad =
z_2^{-1}\delta\Big(\frac{z_1-z_0}{z_2}\Big)Y_i(\ty(a,z_0)b,z_2)v
\end{split}
\end{equation}
under the projection $p_M$, we obtain
\begin{equation}\label{eq:Jacobi-2}
\begin{split}
&
z_0^{-1}\delta\Big(\frac{z_1-z_2}{z_0}\Big)\F_i^M(a,z_1)\G_i^W(b,z_2)v
- z_0^{-1}\delta\Big(\frac{z_2-z_1}{-z_0}\Big)\F_i^M(b,z_2)\G_i^W(a,z_1)v\\
& \qquad =
z_2^{-1}\delta\Big(\frac{z_1-z_0}{z_2}\Big)Y_i(\I(a,z_0)b,z_2)v.
\end{split}
\end{equation}

Likewise, if we take the image of both sides of
\eqref{eq:Jacobi-1} under the projection $p_W$, then
\begin{equation}\label{eq:Jacobi-3}
\begin{split}
&
z_0^{-1}\delta\Big(\frac{z_1-z_2}{z_0}\Big)\F_i^W(a,z_1)\G_i^W(b,z_2)v
- z_0^{-1}\delta\Big(\frac{z_2-z_1}{-z_0}\Big)\F_i^W(b,z_2)\G_i^W(a,z_1)v\\
& \qquad =
z_2^{-1}\delta\Big(\frac{z_1-z_0}{z_2}\Big)\G_i^W(\J(a,z_0)b,z_2)v.
\end{split}
\end{equation}
Comparing Equation \eqref{eq:Jacobi-3} for $i=1$ and $i=2$, we have
\begin{equation*}
\gamma(\mu - 1)
z_2^{-1}\delta\Big(\frac{z_1-z_0}{z_2}\Big)\G_1^W(\J(a,z_0)b,z_2)v
= 0,
\end{equation*}
since $\F_2^M(\,\cdot\,,z) = \lambda \F_1^M(\,\cdot\,,z)$,
$\F_2^W(\,\cdot\,,z) = \mu \F_1^W(\,\cdot\,,z)$, and
$\G_2^W(\,\cdot\,,z) = \gamma \G_1^W(\,\cdot\,,z)$. Now,
$z_2^{-1}\delta\big(\frac{z_1-z_0}{z_2}\big) =
z_1^{-1}\delta\big(\frac{z_2+z_0}{z_1}\big)$ by \cite[Proposition
8.8.5]{FLM} and so the above equation is equivalent to the following
assertion.
\begin{equation*}
\gamma(\mu - 1) (z_2+z_0)^k \G_1^W(\J(a,z_0)b,z_2)v = 0 \quad
\mbox{for all} \quad k \in \Z.
\end{equation*}
This implies that
\begin{equation*}
\gamma(\mu - 1) \G_1^W(\J(a,z_0)b,z_2)v = 0,
\end{equation*}
since $\G_1^W(\J(a,z_0)b,z_2)v \in W((z_0))[[z_2, z_2^{-1}]]$.
Then since $\J(\,\cdot\,,z) \ne 0$ and $\G_1^W(\,\cdot\,,z) \ne
0$, we conclude that $\mu = 1$.

Next, we use Equation \eqref{eq:Jacobi-2}. Since $\I(a,z_0)b \in
M^0((z_0))$, we have $Y_1(\I(a,z_0)b, z_2)v = Y_2(\I(a,z_0)b, z_2)v$
by our assumption. Then it follows from \eqref{eq:Jacobi-2} for
$i=1,2$ that
\begin{equation*}
(\lambda\gamma - 1)
z_2^{-1}\delta\Big(\frac{z_1-z_0}{z_2}\Big)Y_1(\I(a,z_0)b,z_2)v =
0.
\end{equation*}
Since $\I(\,\cdot\,,z) \ne 0$ and $M$ is a simple
$(M^0,Y_1)$-module, a similar argument as above gives that
$\lambda\gamma = 1$.

For $a \in M^0$, $b \in W^0$, $v \in M$, and $w \in W$,
\begin{align*}
& Y_i(a+b,z)(v+w)\\
& \qquad = Y_i(a,z)v + Y_i(a,z)w + \G_i^W(b,z)v + \big( \F_i^M(b,z)
+ \F_i^W(b,z)\big) w.
\end{align*}
Note that $Y_i(a,z)v$, $\F_i^M(b,z)w \in M((z))$ and $Y_i(a,z)w$,
$\G_i^W(b,z)v$, $\F_i^W(b,z)w \in W((z))$. Define $\varphi : U
\rightarrow U$ by $\varphi(u) = \lambda u$ if $u \in M$ and
$\varphi(u) = u$ if $u \in W$. Since $\mu = 1$ and $\lambda\gamma =
1$, we can verify that
\begin{equation*}
Y_2(a+b,z)\varphi(v+w) = \varphi\big( Y_1(a+b,z)(v+w)\big).
\end{equation*}
Thus $\varphi$ is an isomorphism of $V_L^\tau$-modules from
$(U,Y_1)$ onto $(U,Y_2)$. This completes the proof.
\end{proof}

\begin{rmk}
The proof of the above lemma is essentially the same as that of
\cite[Lemma C.3]{LYY}. Consider the Jacobi identity for $a,b \in
W^0$ and $w \in W$ and take the images of both sides of the identity
under the projections $p_M$ and $p_W$, respectively. Then
\begin{align*}
&
z_0^{-1}\delta\Big(\frac{z_1-z_2}{z_0}\Big)\F_i^M(a,z_1)\F_i^W(b,z_2)w
- z_0^{-1}\delta\Big(\frac{z_2-z_1}{-z_0}\Big)\F_i^M(b,z_2)\F_i^W(a,z_1)w\\
& \qquad =
z_2^{-1}\delta\Big(\frac{z_1-z_0}{z_2}\Big)\F_i^M(\J(a,z_0)b,z_2)w,
\end{align*}
\begin{align*}
& z_0^{-1}\delta\Big(\frac{z_1-z_2}{z_0}\Big) \big(
\G_i^W(a,z_1)\F_i^M(b,z_2) + \F_i^W(a,z_1)\F_i^W(b,z_2) \big)w\\
& \qquad\qquad - z_0^{-1}\delta\Big(\frac{z_2-z_1}{-z_0}\Big)
\big(
\G_i^W(b,z_2)\F_i^M(a,z_1) + \F_i^W(b,z_2)\F_i^W(a,z_1) \big)w\\
& \qquad = z_2^{-1}\delta\Big(\frac{z_1-z_0}{z_2}\Big) \big(
Y_i(\I(a,z_0)b,z_2) + \F_i^W(\J(a,z_0)b,z_2)\big)w.
\end{align*}

Each of these two equations gives the identical equations in case
of $i=1$ and $i=2$ provided that $\mu =1$ and $\lambda\gamma = 1$.
\end{rmk}

\begin{thm}
There are exactly $30$ inequivalent simple $V_L^\tau$-modules. They
are represented by the $30$ simple $V_L^\tau$-modules listed in
Lemma \ref{lem:30SIMPLES}.
\end{thm}

\begin{proof}
Let $U$ be a simple $V_L^\tau$-module. Then by Lemma \ref{lem:MM},
$U$ contains a simple $M^0$-submodule $M$ isomorphic to a member of
$\M$. Since $U$ is a simple $V_L^\tau$-module, Lemma \ref{lem:VM}
implies that $U = M \op W$ for some simple $M^0$-submodule $W$
isomorphic to a member of $\W$. In fact, the isomorphism class of
$W$ is uniquely determined by $M$. By Lemma \ref{lem:Unique}, $U$
admits a unique $V_L^\tau$-module structure. Since $\M$ consists of
$30$ members, it follows that there are at most $30$ inequivalent
simple $V_L^\tau$-module. Hence the assertion holds.
\end{proof}

\begin{thm}
$V_L^\tau$ is a rational vertex operator algebra.
\end{thm}

\begin{proof}
It is sufficient to show that every $\N$-graded weak
$V_L^\tau$-module $U$ is a sum of simple $V_L^\tau$-modules. Since
$M^0$ is rational, $U$ is a direct sum of simple $M^0$-modules. Thus
by Lemma \ref{lem:MM-or-WW}, we may assume that $U =
\big(\op_{\gamma \in \Gamma} S^\gamma \big) \op \big( \op_{\lambda
\in \Lambda} S^\lambda \big)$, where $S^\gamma$ is isomorphic to a
member of $\M$ and $S^\lambda$ is isomorphic to a member of $\W$. We
know that $V_L^\tau \cdot S^\gamma$ is a simple $V_L^\tau$-module by
Lemma \ref{lem:VM}. Set $N = \sum_{\gamma \in \Gamma} V_L^\tau \cdot
S^\gamma$. Since $U/N$ has no simple $M^0$-submodule isomorphic to a
member of $\M$, it follows from Lemma \ref{lem:MM} that $U = N$ and
the proof is complete.
\end{proof}

\begin{cor}\label{cor:generator}
The Zhu algebra $A(V_L^\tau)$ of $V_L^\tau$ is a $51$ dimensional
semisimple associative algebra isomorphic to a direct sum of $23$
copies of the one dimensional algebra $\C$ and $7$ copies of the
algebra $\mathrm{Mat}_2(\C)$ of $2 \times 2$ matrices. Moreover,
$A(V_L^\tau)$ is generated by $[\tom^1]$, $[\tom^2]$, $[J]$, $[K]$,
and $[P]$.
\end{cor}

\begin{proof}
Since $V_L^\tau$ is rational, $A(V_L^\tau)$ is a finite dimensional
semisimple associative algebra (cf. \cite[Theorem 8.1]{DLM1},
\cite[Theorem 2.2.3]{Z}). We know all the simple $V_L^\tau$-modules
and the action of $[\tom^1]$, $[\tom^2]$, $[J]$, $[K]$, and $[P]$ on
their top levels in Section \ref{sec:Structure}. Hence we can
determine the structure of $A(V_L^\tau)$ as in the assertion.
\end{proof}

\appendix
\section{Some fusion rules for  $M(0)$}\label{app:fusion-rules}
We give a proof of the fusion rules
\begin{align*}
W(0) \times M_T(\tau^i)(\vep) &= W_T(\tau^i)(\vep),\\
W(0) \times W_T(\tau^i)(\vep) &= M_T(\tau^i)(\vep) +
W_T(\tau^i)(\vep),
\end{align*}
$i=1,2$, $\vep = 0,1,2$ of simple $M(0)$-modules in
\eqref{eq:W0-fusion}.

Recall that $V_L^\tau \cong M^0 \op W^0$, where $M^0 = M(0) \ot
M_t^0$ and $W^0 = W(0)\ot W_t^0$. Set $\hM_T(\tau^i)(\vep) =
M_T(\tau^i)(\vep) \ot M_t^0$ and $\hW_T(\tau^i)(\vep) =
W_T(\tau^i)(\vep) \ot W_t^0$, which are simple $M^0$-modules. Then
\begin{align*}
V_L^{T_{\chi_0}}(\tau)(\vep) &\cong \hM_T(\tau)(\vep) \op
\hW_T(\tau)(\vep),\\
V_L^{T_{\chi'_0}}(\tau^2)(\vep) &\cong \hM_T(\tau^2)(\vep) \op
\hW_T(\tau^2)(\vep)
\end{align*}
as $M^0$-modules by \eqref{eq:Twisted-1-vep} and
\eqref{eq:Twisted-2-vep}. Denote by $Y_1(\, \cdot \,,z)$ (resp.
$Y_2(\, \cdot \,,z)$) the vertex operator of the simple
$V_L^\tau$-module $V_L^{T_{\chi_0}}(\tau)(\vep)$ (resp.
$V_L^{T_{\chi'_0}}(\tau^2)(\vep)$). Let $p_M :
V_L^{T_{\chi_0}}(\tau)(\vep) \rightarrow \hM_T(\tau)(\vep)$ and $p_W
: V_L^{T_{\chi_0}}(\tau)(\vep) \rightarrow \hW_T(\tau)(\vep)$ be
projections. We also use the same symbol $p_M$ or $p_W$ to denote a
projection from $V_L^{T_{\chi'_0}}(\tau^2)(\vep)$ onto
$\hM_T(\tau^2)(\vep)$ or onto $\hW_T(\tau^2)(\vep)$. We fix $i =
1,2$ and $\vep = 0,1,2$. For simplicity of notation, set $\hM =
\hM_T(\tau^i)(\vep)$ and $\hW = \hW_T(\tau^i)(\vep)$.

Let $\F_i^M(a,z)w = p_M Y_i(a,z)w$ and $\F_i^W(a,z)w = p_W
Y_i(a,z)w$ for $a \in W^0$ and $w \in \hW$. Then $\F_i^M(\,\cdot
\,,z)$ and $\F_i^W(\,\cdot \,,z)$ are intertwining operators of type
$\binom{\hM}{W^0 \ \hW}$ and $\binom{\hW}{W^0 \ \hW}$, respectively.
Likewise, let $\G_i^W(a,z)v = Y_i(a,z)v$ for $a \in W^0$ and $v \in
\hM$. Then $\G_i^W(\,\cdot\,,z)$ is an intertwining operator of type
$\binom{\hW}{W^0 \ \hM}$, since the fusion rule $W_t^0 \times M_t^0
= W_t^0$ of $M_t^0$-modules implies that $W^0 \cdot \hM = \spn \{a_n
\hM\,|\, a \in W^0, n \in \Z\}$ is contained in $\hW$. If
$\G_i^W(\,\cdot\,,z)=0$, then $V_L^\tau \cdot \hM = (M^0 + W^0)\cdot
\hM \subset \hM$. This is a contradiction, since
$V_L^{T_{\chi_0}}(\tau)(\vep)$ and $V_L^{T_{\chi'_0}}(\tau^2)(\vep)$
are simple $V_L^\tau$-modules. Thus $\G_i^W(\,\cdot\,,z) \ne 0$.
Similarly, $\F_i^M(\,\cdot\,,z) \ne 0$. Indeed, if
$\F_i^M(\,\cdot\,,z) = 0$, then $V_L^\tau \cdot \hW \subset \hW$,
which is a contradiction. Assume that $\F_i^W(\,\cdot\,,z) = 0$.
Then $W^0 \cdot \hW \subset \hM$ and so $W^0 \cdot (W^0 \cdot \hW)
\subset \hW$. However, $W^0 \cdot (W^0 \cdot \hW) \supset (W^0 \cdot
W^0) \cdot \hW = V_L^\tau \cdot \hW$ by  Lemma \ref{lem:Formula2}
and \eqref{eq:W-cdot-W}. This contradiction implies that
$\F_i^W(\,\cdot\,,z) \ne 0$.

Restricting the three nonzero intertwining operators
$\F_i^M(\,\cdot\,,z)$, $\F_i^W(\,\cdot\,,z)$, and
$\G_i^W(\,\cdot\,,z)$ to the first component of each of the tensor
products $W^0 = W(0) \ot W_t^0$, $\hM = M_T(\tau^i)(\vep) \ot
M_t^0$, and $\hW = W_T(\tau^i)(\vep) \ot W_t^0$, we obtain nonzero
intertwining operators of type $\binom{M_T(\tau^i)(\vep)}{W(0) \
W_T(\tau^i)(\vep)}$, $\binom{W_T(\tau^i)(\vep)}{W(0) \
W_T(\tau^i)(\vep)}$, and $\binom{W_T(\tau^i)(\vep)}{W(0) \
M_T(\tau^i)(\vep)}$ for $M(0)$-modules, respectively.

Let $N^2$ be one of $M_T(\tau^i)(\vep)$, $W_T(\tau^i)(\vep)$,
$i=1,2$, $\vep = 0,1,2$ and let $N^3$ be any of the $20$ simple
$M(0)$-modules. Then the top level $N^j_{(0)}$ of $N^j$ is of
dimension one. By \cite{DLTYY}, the Zhu algebra $A(M(0))$ of
$M(0)$ is generated by $[\tom^1]$ and $[J]$. Moreover, we know the
action of $o(\tom^1)$ and $o(J)$ on $N^j_{(0)}$. Thus by a similar
argument as in pages 192 and 193 of \cite{Tanabe}, we can
calculate that the dimension of
\begin{equation*}
\Hom_{A(M(0))}( A(W(0)) \ot_{A(M(0))} N^2_{(0)}, N^3_{(0)})
\end{equation*}
is at most one and it is equal to one if and only if the pair
$(N^2, N^3)$ is one of
\begin{equation*}
(M_T(\tau^i)(\vep), W_T(\tau^i)(\vep)),\quad (W_T(\tau^i)(\vep),
M_T(\tau^i)(\vep)),\quad \text{or } (W_T(\tau^i)(\vep),
W_T(\tau^i)(\vep))
\end{equation*}
for $i=1,2$, $\vep = 0,1,2$. Note that $W(0)$ was denoted by
$W_k^{0(0)}$ in \cite{Tanabe}. Now, the desired fusion rules are
obtained by \cite[Proposition 2.10 and Corollary 2.13]{Li2}.

\section{Some vectors in $V_L^\tau$}\label{app:Vectors-for-Zhu}
We describe certain vectors in $V_L^\tau$ which are used for
determining the Zhu algebra of $V_L^\tau$ in Section
\ref{sec:Classification}. The calculation was done by a computer
algebra system Risa/Asir.

{\footnotesize
\begin{align*}
P_1(J_1 P) = &
-(312/7)J_{-1}\1
-(80/7)K_1P, \\
P_0(J_1 P) = &
-(104/7)J_{-2}\1
-(15/7)(\tom^1)_0K_1P
-(40/17)(\tom^2)_0K_1P
-(40/17)K_0P, \\
P_{-1}(J_1 P)
= &
-(1404/119)(\tom^1)_{-1}J_{-1}\1
-(312/7)(\tom^2)_{-1}J_{-1}\1
-(156/119)J_{-3}\1
-(208/35)(\tom^1)_{-1}K_1P \\
&
+(12/7)(\tom^1)_0(\tom^1)_0K_1P
-(15/34)(\tom^1)_0(\tom^2)_0K_1P
-(832/357)(\tom^2)_{-1}K_1P \\
&
-(15/34)(\tom^1)_0K_0P
-(16/17)K_{-1}P, \\
P_{-2}(J_1 P)
= &
-(156/17)(\tom^1)_{-2}J_{-1}\1
-(156/7)(\tom^2)_{-2}J_{-1}\1
-(208/119)(\tom^1)_{-1}J_{-2}\1 \\
&
-(104/7)(\tom^2)_{-1}J_{-2}\1
+(78/119)J_{-4}\1
-(52/5)(\tom^1)_{-2}K_1P
+(286/35)(\tom^1)_{-1}(\tom^1)_0K_1P \\
&
-(9/7)(\tom^1)_0(\tom^1)_0(\tom^1)_0K_1P
-(104/85)(\tom^1)_{-1}(\tom^2)_0K_1P
+(6/17)(\tom^1)_0(\tom^1)_0(\tom^2)_0K_1P \\
&
-(52/119)(\tom^1)_0(\tom^2)_{-1}K_1P
+(132/259)(\tom^2)_{-2}K_1P
-(164/259)(\tom^2)_{-1}(\tom^2)_0K_1P \\
&
-(104/85)(\tom^1)_{-1}K_0P
+(6/17)(\tom^1)_0(\tom^1)_0K_0P
-(164/259)(\tom^2)_{-1}K_0P \\
&
-(3/17)(\tom^1)_0K_{-1}P
-(492/259)(\tom^2)_0K_{-1}P
+(1152/259)K_{-2}P,
\end{align*}
}
{\footnotesize
\begin{align*}
(J_1 P)_1P = &
-(312/7)J_{-1}\1
-(80/7)K_1P, \\
(J_1 P)_0P = &
-(208/7)J_{-2}\1
-(65/7)(\tom^1)_0K_1P
-(40/17)(\tom^2)_0K_1P
-(40/17)K_0P, \\
(J_1 P)_{-1}P
= &
-(1404/119)(\tom^1)_{-1}J_{-1}\1
-(312/7)(\tom^2)_{-1}J_{-1}\1
-(1924/119)J_{-3}\1
-(208/35)(\tom^1)_{-1}K_1P \\
&
-(13/7)(\tom^1)_0(\tom^1)_0K_1P
-(65/34)(\tom^1)_0(\tom^2)_0K_1P
-(832/357)(\tom^2)_{-1}K_1P \\
&
-(65/34)(\tom^1)_0K_0P
-(16/17)K_{-1}P,\\
(J_1 P)_{-2}P
= &
-(312/119)(\tom^1)_{-2}J_{-1}\1
-(156/7)(\tom^2)_{-2}J_{-1}\1
-(1196/119)(\tom^1)_{-1}J_{-2}\1 \\
&
-(208/7)(\tom^2)_{-1}J_{-2}\1
-(78/17)J_{-4}\1
+(156/35)(\tom^1)_{-2}K_1P
-(494/35)(\tom^1)_{-1}(\tom^1)_0K_1P \\
&
+(13/6)(\tom^1)_0(\tom^1)_0(\tom^1)_0K_1P
-(104/85)(\tom^1)_{-1}(\tom^2)_0K_1P
-(13/34)(\tom^1)_0(\tom^1)_0(\tom^2)_0K_1P \\
&
-(676/357)(\tom^1)_0(\tom^2)_{-1}K_1P
+(132/259)(\tom^2)_{-2}K_1P
-(164/259)(\tom^2)_{-1}(\tom^2)_0K_1P \\
&
-(104/85)(\tom^1)_{-1}K_0P
-(13/34)(\tom^1)_0(\tom^1)_0K_0P
-(164/259)(\tom^2)_{-1}K_0P \\
&
-(13/17)(\tom^1)_0K_{-1}P
-(492/259)(\tom^2)_0K_{-1}P
+(1152/259)K_{-2}P,
\end{align*}
}
{\footnotesize
\begin{align*}
P_1(K_1 P) = &
-(63/13)K_{-1}\1
+(20/13)J_1P, \\
P_0(K_1 P) = &
-(21/13)K_{-2}\1
+(10/23)(\tom^1)_0J_1P
-(15/13)(\tom^2)_0J_1P
+(10/23)J_0P, \\
P_{-1}(K_1 P)
= &
-(168/13)(\tom^1)_{-1}K_{-1}\1
-(189/299)(\tom^2)_{-1}K_{-1}\1
-(126/299)K_{-3}\1 \\
&
+(1176/1495)(\tom^1)_{-1}J_1P
+(42/253)(\tom^1)_0(\tom^1)_0J_1P
-(15/46)(\tom^1)_0(\tom^2)_0J_1P \\
&
-(2/13)(\tom^2)_{-1}J_1P
+(84/253)(\tom^1)_0J_0P
-(15/46)(\tom^2)_0J_0P
-(144/115)J_{-1}P, \\
P_{-2}(K_1 P)
= &
-(84/13)(\tom^1)_{-2}K_{-1}\1
-(9/23)(\tom^2)_{-2}K_{-1}\1
-(56/13)(\tom^1)_{-1}K_{-2}\1 \\
&
-(48/299)(\tom^2)_{-1}K_{-2}\1
-(27/299)K_{-4}\1
-(24/65)(\tom^1)_{-2}J_1P
+(36/65)(\tom^1)_{-1}(\tom^1)_0J_1P \\
&
-(347/3354)(\tom^1)_0(\tom^1)_0(\tom^1)_0J_1P
-(882/1495)(\tom^1)_{-1}(\tom^2)_0J_1P \\
&
-(63/506)(\tom^1)_0(\tom^1)_0(\tom^2)_0J_1P
-(1/23)(\tom^1)_0(\tom^2)_{-1}J_1P
+(18/13)(\tom^2)_{-2}J_1P \\
&
-(17/13)(\tom^2)_{-1}(\tom^2)_0J_1P
+(36/65)(\tom^1)_{-1}J_0P
-(347/1118)(\tom^1)_0(\tom^1)_0J_0P \\
&
-(63/253)(\tom^1)_0(\tom^2)_0J_0P
-(1/23)(\tom^2)_{-1}J_0P
+(108/65)(\tom^1)_0J_{-1}P \\
&
+(108/115)(\tom^2)_0J_{-1}P
-(228/65)J_{-2}P,
\end{align*}
}
{\footnotesize
\begin{align*}
(K_1 P)_1 P = &
-(63/13)K_{-1}\1
+(20/13)J_1P, \\
(K_1 P)_0 P = &
-(42/13)K_{-2}\1
+(10/23)(\tom^1)_0J_1P
+(35/13)(\tom^2)_0J_1P
+(10/23)J_0P, \\
(K_1 P)_{-1}P
= &
-(168/13)(\tom^1)_{-1}K_{-1}\1
-(189/299)(\tom^2)_{-1}K_{-1}\1
-(609/299)K_{-3}\1 \\
&
+(1176/1495)(\tom^1)_{-1}J_1P
+(42/253)(\tom^1)_0(\tom^1)_0J_1P
+(35/46)(\tom^1)_0(\tom^2)_0J_1P \\
&
+(28/13)(\tom^2)_{-1}J_1P
+(84/253)(\tom^1)_0J_0P
+(35/46)(\tom^2)_0J_0P
-(144/115)J_{-1}P, \\
(K_1 P)_{-2}P
= &
-(84/13)(\tom^1)_{-2}K_{-1}\1
-(72/299)(\tom^2)_{-2}K_{-1}\1
-(112/13)(\tom^1)_{-1}K_{-2}\1 \\
&
-(141/299)(\tom^2)_{-1}K_{-2}\1
-(27/23)K_{-4}\1
-(24/65)(\tom^1)_{-2}J_1P
+(36/65)(\tom^1)_{-1}(\tom^1)_0J_1P \\
&
-(347/3354)(\tom^1)_0(\tom^1)_0(\tom^1)_0J_1P
+(2058/1495)(\tom^1)_{-1}(\tom^2)_0J_1P \\
&
+(147/506)(\tom^1)_0(\tom^1)_0(\tom^2)_0J_1P
+(14/23)(\tom^1)_0(\tom^2)_{-1}J_1P
-(7/13)(\tom^2)_{-2}J_1P \\
&
+(28/13)(\tom^2)_{-1}(\tom^2)_0J_1P
+(36/65)(\tom^1)_{-1}J_0P
-(347/1118)(\tom^1)_0(\tom^1)_0J_0P \\
&
+(147/253)(\tom^1)_0(\tom^2)_0J_0P
+(14/23)(\tom^2)_{-1}J_0P
+(108/65)(\tom^1)_0J_{-1}P \\
&
-(252/115)(\tom^2)_0J_{-1}P
-(228/65)J_{-2}P,
\end{align*}
}
{\footnotesize
\begin{align*}
(J_1 P)_2(J_1 P) = &
8112(\tom^1)_{-2}\1
+1872(\tom^2)_{-2}\1,\\
(J_1 P)_1(J_1 P) = &
8112(\tom^1)_{-1}(\tom^1)_{-1}\1
+16224(\tom^1)_{-1}(\tom^2)_{-1}\1
+864(\tom^2)_{-3}\1
+432(\tom^2)_{-1}(\tom^2)_{-1}\1
-8J_1K_1P,\\
(J_1 P)_0(J_1 P)
= &
8112(\tom^1)_{-2}(\tom^1)_{-1}\1
+8112(\tom^1)_{-2}(\tom^2)_{-1}\1
+8112(\tom^1)_{-1}(\tom^2)_{-2}\1
+576(\tom^2)_{-4}\1 \\
&
+432(\tom^2)_{-2}(\tom^2)_{-1}\1
-(52/23)(\tom^1)_0J_1K_1P
-(28/17)(\tom^2)_0J_1K_1P
-(28/17)J_1K_0P \\
&
-(52/23)J_0K_1P, \\
(J_1 P)_{-1}(J_1 P)
= &
-(146016/77)(\tom^1)_{-5}\1
-(1014000/1309)(\tom^1)_{-3}(\tom^1)_{-1}\1
+(2314962/1309)(\tom^1)_{-2}(\tom^1)_{-2}\1 \\
&
+(1565616/1309)(\tom^1)_{-1}(\tom^1)_{-1}(\tom^1)_{-1}\1
+4056(\tom^1)_{-2}(\tom^2)_{-2}\1 \\
&
+8112(\tom^1)_{-1}(\tom^1)_{-1}(\tom^2)_{-1}\1
+3744(\tom^1)_{-1}(\tom^2)_{-3}\1
+1872(\tom^1)_{-1}(\tom^2)_{-1}(\tom^2)_{-1}\1 \\
&
+(7776/23)(\tom^2)_{-5}\1
+(5472/23)(\tom^2)_{-3}(\tom^2)_{-1}\1
+(2304/23)(\tom^2)_{-2}(\tom^2)_{-2}\1 \\
&
+(432/23)(\tom^2)_{-1}(\tom^2)_{-1}(\tom^2)_{-1}\1
-(3042/187)J_{-1}J_{-1}\1
+(5876/805)(\tom^1)_{-1}J_1K_1P \\
&
+(5590/1771)(\tom^1)_0(\tom^1)_0J_1K_1P
-(182/391)(\tom^1)_0(\tom^2)_0J_1K_1P
-(416/255)(\tom^2)_{-1}J_1K_1P \\
&
+(11180/1771)(\tom^1)_0J_0K_1P
-(182/391)(\tom^2)_0J_0K_1P
-(182/391)(\tom^1)_0J_1K_0P \\
&
-(4056/115)J_{-1}K_1P
-(182/391)J_0K_0P
-(56/85)J_1K_{-1}P,
\end{align*}
}
{\footnotesize
\begin{align*}
(K_1 P)_2(K_1 P) = &
-672(\tom^1)_{-2}\1
-882(\tom^2)_{-2}\1,\\
(K_1 P)_1(K_1 P)
= &
-144(\tom^1)_{-3}\1
-432(\tom^1)_{-1}(\tom^1)_{-1}\1
-4704(\tom^1)_{-1}(\tom^2)_{-1}\1 \\
&
-(2646/13)(\tom^2)_{-3}\1
-(7056/13)(\tom^2)_{-1}(\tom^2)_{-1}\1
-(30/13)J_1K_1P,\\
(K_1 P)_0(K_1 P)
= &
-96(\tom^1)_{-4}\1
-432(\tom^1)_{-2}(\tom^1)_{-1}\1
-2352(\tom^1)_{-2}(\tom^2)_{-1}\1
-2352(\tom^1)_{-1}(\tom^2)_{-2}\1 \\
&
-(1764/13)(\tom^2)_{-4}\1
-(7056/13)(\tom^2)_{-2}(\tom^2)_{-1}\1
-(15/23)(\tom^1)_0J_1K_1P \\
&
-(105/221)(\tom^2)_0J_1K_1P
-(15/23)J_0K_1P
-(105/221)J_1K_0P,\\
(K_1 P)_{-1}(K_1 P)
= &
(2688/11)(\tom^1)_{-5}\1
+(15792/187)(\tom^1)_{-3}(\tom^1)_{-1}\1
-(50316/187)(\tom^1)_{-2}(\tom^1)_{-2}\1 \\
&
-(32928/187)(\tom^1)_{-1}(\tom^1)_{-1}(\tom^1)_{-1}\1
-504(\tom^1)_{-3}(\tom^2)_{-1}\1
-1176(\tom^1)_{-2}(\tom^2)_{-2}\1 \\
&
-1512(\tom^1)_{-1}(\tom^1)_{-1}(\tom^2)_{-1}\1
-(7056/13)(\tom^1)_{-1}(\tom^2)_{-3}\1 \\
&
-(18816/13)(\tom^1)_{-1}(\tom^2)_{-1}(\tom^2)_{-1}\1
+(49392/299)(\tom^2)_{-5}\1
+(31752/299)(\tom^2)_{-3}(\tom^2)_{-1}\1 \\
&
-(73647/299)(\tom^2)_{-2}(\tom^2)_{-2}\1
-(42336/299)(\tom^2)_{-1}(\tom^2)_{-1}(\tom^2)_{-1}\1
-(180/187)J_{-1}J_{-1}\1 \\
&
-(1764/1495)(\tom^1)_{-1}J_1K_1P
-(63/253)(\tom^1)_0(\tom^1)_0J_1K_1P
-(105/782)(\tom^1)_0(\tom^2)_0J_1K_1P\\
&
-(812/663)(\tom^2)_{-1}J_1K_1P
-(126/253)(\tom^1)_0J_0K_1P
-(105/782)(\tom^2)_0J_0K_1P \\
&
-(105/782)(\tom^1)_0J_1K_0P
+(216/115)J_{-1}K_1P
-(105/782)J_0K_0P
+(28/221)J_1K_{-1}P,
\end{align*}
}
{\footnotesize
\begin{align*}
(J_1 P)_2(K_1 P) = &
-780(\tom^1)_0P
+720(\tom^2)_0P,\\
(J_1 P)_1(K_1 P) = &
-(2496/5)(\tom^1)_{-1}P
-156(\tom^1)_0(\tom^1)_0P
+585(\tom^1)_0(\tom^2)_0P
+96(\tom^2)_{-1}P,\\
(J_1 P)_0(K_1 P)
= &
(1872/5)(\tom^1)_{-2}P
-(5928/5)(\tom^1)_{-1}(\tom^1)_0P
182(\tom^1)_0(\tom^1)_0(\tom^1)_0P
+(1872/5)(\tom^1)_{-1}(\tom^2)_0P \\
&
+117(\tom^1)_0(\tom^1)_0(\tom^2)_0P
+78(\tom^1)_0(\tom^2)_{-1}P
-864(\tom^2)_{-2}P
+816(\tom^2)_{-1}(\tom^2)_0P,\\
(J_1 P)_{-1}(K_1 P)
= &
-36J_{-1}K_{-1}\1
+(3936/5)(\tom^1)_{-3}P
-(2064/5)(\tom^1)_{-2}(\tom^1)_0P
-(864/5)(\tom^1)_{-1}(\tom^1)_{-1}P \\
&
-(2352/5)(\tom^1)_{-1}(\tom^1)_0(\tom^1)_0P
+98(\tom^1)_0(\tom^1)_0(\tom^1)_0(\tom^1)_0P
-(1404/5)(\tom^1)_{-2}(\tom^2)_0P \\
&
+(4446/5)(\tom^1)_{-1}(\tom^1)_0(\tom^2)_0P
-(273/2)(\tom^1)_0(\tom^1)_0(\tom^1)_0(\tom^2)_0P \\
&
+(1248/25)(\tom^1)_{-1}(\tom^2)_{-1}P
+(78/5)(\tom^1)_0(\tom^1)_0(\tom^2)_{-1}P
-702(\tom^1)_0(\tom^2)_{-2}P \\
&
+663(\tom^1)_0(\tom^2)_{-1}(\tom^2)_0P
-576(\tom^2)_{-3}P
+144(\tom^2)_{-2}(\tom^2)_0P
+(1152/5)(\tom^2)_{-1}(\tom^2)_{-1}P,
\end{align*}
}
{\footnotesize
\begin{align*}
(K_1 P)_2(J_1 P) = &
-180(\tom^1)_0P
-1680(\tom^2)_0P,\\
(K_1 P)_1(J_1 P) = &
-(2496/5)(\tom^1)_{-1}P
144(\tom^1)_0(\tom^1)_0P
-315(\tom^1)_0(\tom^2)_0P
-1344(\tom^2)_{-1}P,\\
(K_1 P)_0(J_1 P)
= &
-(4368/5)(\tom^1)_{-2}P
+(3432/5)(\tom^1)_{-1}(\tom^1)_0P
-108(\tom^1)_0(\tom^1)_0(\tom^1)_0P \\
&
-(4368/5)(\tom^1)_{-1}(\tom^2)_0P
+252(\tom^1)_0(\tom^1)_0(\tom^2)_0P
-252(\tom^1)_0(\tom^2)_{-1}P
+336(\tom^2)_{-2}P \\
&
-1344(\tom^2)_{-1}(\tom^2)_0P,\\
(K_1 P)_{-1}(J_1 P)
= &
-36J_{-1}K_{-1}\1
-(2304/5)(\tom^1)_{-3}P
-(504/5)(\tom^1)_{-2}(\tom^1)_0P
-(864/5)(\tom^1)_{-1}(\tom^1)_{-1}P \\
&
+(2328/5)(\tom^1)_{-1}(\tom^1)_0(\tom^1)_0P
-72(\tom^1)_0(\tom^1)_0(\tom^1)_0(\tom^1)_0P
-(7644/5)(\tom^1)_{-2}(\tom^2)_0P \\
&
+(6006/5)(\tom^1)_{-1}(\tom^1)_0(\tom^2)_0P
-189(\tom^1)_0(\tom^1)_0(\tom^1)_0(\tom^2)_0P
-(17472/25)(\tom^1)_{-1}(\tom^2)_{-1}P\\
&
+(1008/5)(\tom^1)_0(\tom^1)_0(\tom^2)_{-1}P
+63(\tom^1)_0(\tom^2)_{-2}P
-252(\tom^1)_0(\tom^2)_{-1}(\tom^2)_0P \\
&
+864(\tom^2)_{-3}P
-96(\tom^2)_{-2}(\tom^2)_0P
-(4608/5)(\tom^2)_{-1}(\tom^2)_{-1}P.
\end{align*}
}


\begin{thebibliography}{99}
\bibitem{ABD}
T. Abe, G. Buhl and C. Dong, Rationality, regularity, and
$C_2$-cofiniteness, \emph{Trans. Amer. Math. Soc.} {\bfseries 356}
(2004), 3391--3402.

\bibitem{AD}
T. Abe and C. Dong, Classification of irreducible modules for the
vertex operator algebra $V_L^+$: General case, \emph{J. Algebra}
{\bfseries 273} (2004), 657--685.

\bibitem{ADL}
T. Abe, C. Dong and H.S. Li, Fusion rules for the vertex operator
algebra $M(1)$ and $V_L^+$, \emph{Comm. Math. Phys.} {\bfseries 253}
(2005), 171--219.

\bibitem{BMP}
P. Bouwknegt, J. McCarthy and K. Pilch,
\emph{The $\mathcal{W}_3$ Algebra}, Lecture Notes in Physics,
{\bfseries m42}, Springer, Berlin 1996.

\bibitem{Buhl} G. Buhl, A spanning set for VOA modules,
\emph{J. Algebra} {\bfseries 254} (2002), 125--151.

\bibitem{D1}
C. Dong, Vertex algebras associated with even lattices,
\emph{J. Algebra} {\bfseries 161} (1993), 245--265.

\bibitem{D2}
C. Dong, Twisted modules for vertex algebras associated with even
lattices,
\emph{J. Algebra} {\bfseries 165} (1994), 91--112.

\bibitem{DLTYY}
C. Dong, C.H. Lam, K. Tanabe, H. Yamada and K. Yokoyama,
$\Z_3$ symmetry and $W_3$ algebra in lattice vertex
operator algebras,
\emph{Pacific J. Math.} {\bfseries 215} (2004), 245--296.

\bibitem{DL2}
C. Dong and J. Lepowsky, The algebraic structure of
relative twisted vertex operators,
\emph{J. Pure and Applied Algebra} {\bfseries 110}(1996),
259--295.

\bibitem{DLM0}
C. Dong, H.S. Li and  G. Mason, Regularity of rational vertex
operator algebras, \emph{Adv. Math.} {\bfseries 132} (1997),
148--166.

\bibitem{DLM1}
C. Dong, H.S. Li and  G. Mason, Twisted representations of vertex
operator algebras, \emph{Math. Ann.} {\bfseries 310} (1998),
571--600.

\bibitem{DLM2}
C. Dong, H.S. Li and  G. Mason, Modular-invariance of trace
functions in orbifold theory and generalized moonshine,
\emph{Comm. Math. Phys.} {\bfseries 214} (2000), 1--56.

\bibitem{DLMN}
C. Dong, H.S. Li, G. Mason and S. P. Norton,
Associative subalgebras of the Griess algebra and related topics,
in: \emph{Proc. of the Conference on the Monster and Lie
algebras at The Ohio State University, May 1996,}
ed. by J. Ferrar and K. Harada,
Walter de Gruyter, Berlin-New York, 1998, 27--42.

\bibitem{DM}
C. Dong and G. Mason, On quantum Galois theory,
\emph{Duke Math. J.} {\bfseries 86} (1997), 305--321.

\bibitem{DY}
C. Dong and G. Yamskulna, Vertex operator algebras, generalized
doubles and dual pairs, \emph{Math. Z.} {\bfseries 241} (2002),
397--423.

\bibitem{FHL}  I. B. Frenkel, Y. Huang and J. Lepowsky,
\emph{On axiomatic approaches to vertex operator algebras and
modules}, Mem. Amer. Math. Soc. {\bfseries 104}, 1993.

\bibitem{FLM}
I. B. Frenkel, J. Lepowsky and A. Meurman, \emph{Vertex Operator
Algebras and the Monster}, Pure and Applied Math., Vol. {\bfseries
134}, Academic Press, 1988.

\bibitem{GN}
M. R. Gaberdiel and A. Neitzke, Rationality, quasirationality and
finite $W$-algebras, \emph{Comm. Math. Phys.} {\bfseries 238}
(2003), 305--331.

\bibitem{KLY1}
K. Kitazume, C.H. Lam and H. Yamada, Decomposition of the moonshine
vertex operator algebra as Virasoro modules, \emph{J. Algebra}
{\bfseries 226} (2000), 893--919.

\bibitem{KLY2}
K. Kitazume, C.H. Lam and H. Yamada, $3$-state Potts model,
moonshine vertex operator algebra and $3A$ elements of the monster
group, \emph{Internat. Math. Res. Notices}, No. {\bfseries 23}
(2003), 1269--1303.

\bibitem{KMY}
M. Kitazume, M. Miyamoto and H. Yamada, Ternary codes and vertex
operator algebras, \emph{J. Algebra} {\bfseries 223} (2000),
379--395.

\bibitem{LY}
C.H. Lam and H. Yamada, $\mathbb{Z}_{2}\times
\mathbb{Z}_{2}$ codes and vertex operator algebras,
\emph{J. Algebra} {\bfseries 224} (2000), 268--291.

\bibitem{LYY}
C.H. Lam, H. Yamada and H. Yamauchi, McKay's observation and vertex
operator algebras generated by two conformal vectors of central
charge $1/2$, \emph{Internat. Math. Res. Papers} {\bfseries 2005:3}
(2005), 117--181.

\bibitem{L}
J. Lepowsky, Calculus of twisted vertex operators,
\emph{Proc. Natl. Acad. Sci. USA} {\bfseries 82} (1985),
8295--8299.

\bibitem{LL}
J. Lepowsky and H.S. Li, \emph{Introduction to Vertex Operator
Algebras and their Representations}, Progress in Mathematics,
\textbf{227}, Birkhauser Boston, Inc., Boston, MA, 2004.

\bibitem{Li1}
H.S. Li, Some finiteness properties of regular vertex operator
algebras, \emph{J. Algebra} {\bfseries 212} (1999), 495--514.

\bibitem{Li2}
H.S. Li, Determining fusion rules by $A(V)$-modules and bimodules,
\emph{J. Algebra} {\bfseries 212} (1999), 515--556.

\bibitem{Li3}
H.S. Li, The regular representation, Zhu's $A(V)$-theory, and
induced modules, \emph{J. Algebra} {\bfseries 238} (2001), 159--193.

\bibitem{M2}
M. Miyamoto, $3$-State Potts model and automorphism of vertex
operator algebra of order $3$, \emph{J. Algebra} {\bfseries 239}
(2001), 56--76.

\bibitem{MT}
M. Miyamoto and K. Tanabe, Uniform product of $A_{g,n}(V)$ for an
orbifold model $V$ and $G$-twisted Zhu algebra, \emph{J. Algebra}
{\bfseries 274} (2004), 80--96.

\bibitem{NT}
K. Nagatomo and A. Tsuchiya, Conformal field theories associated to
regular chiral vertex operator algebras. I. Theories over the
projective line, \emph{Duke Math. J.} {\bfseries 128} (2005),
393--471.

\bibitem{Tanabe}
K. Tanabe, On intertwining operators and finite automorphism groups
of vertex operator algebras, \emph{J. Algebra} {\bfseries 287}
(2005), 174--198.

\bibitem{Wang1}
W. Wang, Rationality of Virasoro vertex operator algebras,
\emph{Duke Math. J.} {\bfseries 71}, \emph{Internat. Math. Res.
Notice} (1993), 197--211.

\bibitem{Yamauchi}
H. Yamauchi, Modularity on vertex operator algebras arising from
semisimple primary vectors, \emph{Internat. J. Math.} {\bfseries
15} (2004), 87--109.

\bibitem{Yamskulna}
G. Yamskulna, $C_2$-cofiniteness of vertex operator algebra $V_L^+$
when $L$ is a rank one lattice, \emph{Comm. Algebra} {\bfseries 32}
(2004), 927--954.

\bibitem{Z}
Y. Zhu, Modular invariance of characters of vertex operator
algebras,
\emph{J. Amer. Math. Soc.} {\bfseries 9} (1996),
237--302.

\end{thebibliography}
\end{document}